\input amstex
\let\myfrac=\frac
\input eplain
\let\frac=\myfrac
\input epsf




\loadeufm \loadmsam \loadmsbm
\message{symbol names}\UseAMSsymbols\message{,}

\font\myfontdefault=cmr10

\font\mytdmchapfont=cmb10 at 14pt
\font\mytdmheadfont=cmb10 at 10pt
\font\mytdmsubheadfont=cmr10

\magnification 1200
\newif\ifinappendices
\newif\ifundefinedreferences
\newif\ifchangedreferences
\newif\ifloadreferences
\newif\ifmakebiblio
\newif\ifmaketdm

\undefinedreferencestrue
\changedreferencesfalse


\loadreferencestrue
\makebibliofalse
\maketdmfalse

\def\headpenalty{-400}     
\def\proclaimpenalty{-200} 

%
%

\def\alphanum#1{\ifcase #1 _\or A\or B\or C\or D\or E\or F\or G\or H\or I\or J\or K\or L\or M\or N\or O\or P\or Q\or R\or S\or T\or U\or V\or W\or X\or Y\or Z\fi}
\def\gobbleeight#1#2#3#4#5#6#7#8{}

\newwrite\references
\newwrite\tdm
\newwrite\biblio

\newcount\chapno
\newcount\headno
\newcount\subheadno
\newcount\procno
\newcount\figno
\newcount\citationno

\def\setcatcodes{%
\catcode`\!=0 \catcode`\\=11}%

\ifloadreferences
    {\catcode`\@=11 \catcode`\_=11%
    \global\def\_@citation@Almgren{1}
\global\def\_@citation@AndBarbBegZegh{2}
\global\def\_@citation@Cabezas{3}
\global\def\_@citation@CaffNirSprI{4}
\global\def\_@citation@CaffNirSprII{5}
\global\def\_@citation@CaffNirSprIII{6}
\global\def\_@citation@CaffNirSprV{7}
\global\def\_@citation@ChenYau{8}
\global\def\_@citation@GallKapMard{9}
\global\def\_@citation@Guan{10}
\global\def\_@citation@GuanSpruckI{11}
\global\def\_@citation@GuanSpruckII{12}
\global\def\_@citation@GuanSpruckIII{13}
\global\def\_@citation@Gut{14}
\global\def\_@citation@HarvLaws{15}
\global\def\_@citation@Huisken{16}
\global\def\_@citation@LabI{17}
\global\def\_@citation@LabII{18}
\global\def\_@citation@LoftI{19}
\global\def\_@citation@LoftII{20}
\global\def\_@citation@MazzPac{21}
\global\def\_@citation@RosSpruck{22}
\global\def\_@citation@SchKra{23}
\global\def\_@citation@SmiFCS{24}
\global\def\_@citation@SmiSLC{25}
\global\def\_@citation@SmiGKM{26}
\global\def\_@citation@SmiHPP{27}
\global\def\_@citation@SmiPKS{28}
\global\def\_@citation@SmiPHD{29}
\global\def\_@proc@TheoremDirichletI{1.1}
\global\def\_@proc@TheoremDirichletIA{1.2}
\global\def\_@proc@TheoremDirichletIB{1.3}
\global\def\_@proc@TheoremDirichletII{1.4}
\global\def\_@proc@TheoremDirichletIII{1.5}
\global\def\_@proc@TheoremDirichletIIIA{1.6}
\global\def\_@head@SpecialLagrangianCurvature{2}
\global\def\_@proc@LemmaConcavityOfSLCurvature{2.1}
\global\def\_@proc@LemmaSLCIsLowestEval{2.2}
\global\def\_@proc@LemmaIntersectionIsStillEpsilonConvex{2.3}
\global\def\_@proc@LemmaIntersectionStillHasCurvatureBoundedBelow{2.4}
\global\def\_@proc@LemmaRemoveableSingularity{2.5}
\global\def\_@proc@LemmaGeomMaxPrinc{2.6}
\global\def\_@proc@LemmaConcavityPreserved{2.8}
\global\def\_@proc@CorollaryConcaveFunction{2.9}
\global\def\_@proc@CorollaryApproximateHessianOfDistance{2.11}
\global\def\_@proc@LemmaSmoothingConvexSets{2.13}
\global\def\_@proc@LemmaCompactnessOfConvexGraphs{3.1}
\global\def\_@proc@LemmaUniformModulusOfContinuity{3.2}
\global\def\_@proc@PropositionConvexityCondition{3.3}
\global\def\_@proc@PropositionSubspacesCloseI{3.5}
\global\def\_@proc@PropositionSecondOrderVariation{3.6}
\global\def\_@proc@CorollarySubspacesCloseII{3.8}
\global\def\_@proc@PropositionControlledDerivative{3.9}
\global\def\_@proc@PropositionSecondOrderVariationII{3.10}
\global\def\_@proc@LemmaBoundaryLowerNormalBounds{3.11}
\global\def\_@proc@CorollaryLowerNormalBound{3.12}
\global\def\_@proc@LemmaConcaveFunction{3.13}
\global\def\_@proc@LemmaLowerBoundOfLaplacian{3.14}
\global\def\_@proc@LemmaNormalDerivative{3.15}
\global\def\_@proc@LemmaBoundedLaplacian{3.16}
\global\def\_@proc@LemmaSecondOrderBoundaryEstimates{3.17}
\global\def\_@proc@LemmaLaplacianOfMeanCurvature{3.18}
\global\def\_@proc@LemmaMeanCurvatureIsSubharmonic{3.19}
\global\def\_@proc@PropositionInteriorSecondOrderBounds{3.20}
\global\def\_@proc@PropositionFormulaForSecondFF{3.21}
\global\def\_@proc@ConvexityOfDistanceFunction{4.4}
\global\def\_@proc@PropositionConvexityInTheLimit{4.5}
\global\def\_@proc@LemmaCompactnessOfConvexSubsetsOfHyperbolicEnds{4.6}
\global\def\_@proc@LemmaConvexGraphsAreFat{4.7}
\global\def\_@proc@PropositionAdaptedDiscs{4.8}
\global\def\_@proc@PropositionConvexCurvesLift{5.2}
\global\def\_@proc@PropositionUniformityOfEpsilonConvexity{5.3}
\global\def\_@proc@LemmaCompactnessOfConvexPseudoImmersions{5.4}
    }%
\else
    \openout\references=references.tex
\fi

\newcount\newchapflag 
\newcount\showpagenumflag 

\global\chapno = -1 
\global\citationno=0
\global\headno = 0
\global\subheadno = 0
\global\procno = 0
\global\figno = 0

\def\resetcounters{%
\global\headno = 0%
\global\subheadno = 0%
\global\procno = 0%
\global\figno = 0%
}

\global\newchapflag=0 
\global\showpagenumflag=0 

\def\chinfo{\ifinappendices\alphanum\chapno\else\the\chapno\fi}%
\def\headinfo{\ifinappendices\alphanum\headno\else\the\headno\fi}%
\def\subheadinfo{\headinfo.\the\subheadno}
\def\procinfo{\headinfo.\the\procno}
\def\figinfo{\the\figno}        
\def\citationinfo{\the\citationno}%
\def\nextheadno{\global\advance\headno by 1 \global\subheadno = 0 \global\procno = 0}
\def\nextsubheadno{\global\advance\subheadno by 1}
\def\nextprocno{\global\advance\procno by 1 \procinfo}
\def\nextfigno{\global\advance\figno by 1 \figinfo}

{\global\let\noe=\noexpand%
%
%
\catcode`\@=11%
\catcode`\_=11%
\setcatcodes%
!global!def!_@@internal@@makeref#1{%
!global!expandafter!def!csname #1ref!endcsname##1{%
!csname _@#1@##1!endcsname%
!expandafter!ifx!csname _@#1@##1!endcsname!relax%
    !write16{#1 ##1 not defined - run saving references}%
    !undefinedreferencestrue%
!fi}}%
!global!def!_@@internal@@makelabel#1{%
!global!expandafter!def!csname #1label!endcsname##1{%
!edef!temptoken{!csname #1info!endcsname}%
!ifloadreferences%
    !expandafter!ifx!csname _@#1@##1!endcsname!relax%
        !write16{#1 ##1 not hitherto defined - rerun saving references}%
        !changedreferencestrue%
    !else%
        !expandafter!ifx!csname _@#1@##1!endcsname!temptoken%
        !else
            !write16{#1 ##1 reference has changed - rerun saving references}%
            !changedreferencestrue%
        !fi%
    !fi%
!else%
    !expandafter!edef!csname _@#1@##1!endcsname{!temptoken}%
    !edef!textoutput{!write!references{\global\def\_@#1@##1{!temptoken}}}%
    !textoutput%
!fi}}%
!global!def!makecounter#1{!_@@internal@@makelabel{#1}!_@@internal@@makeref{#1}}%
!unsetcatcodes%
}
\makecounter{ch}%
\makecounter{head}%
\makecounter{subhead}%
\makecounter{proc}%
\makecounter{fig}%
\makecounter{citation}%
\def\newref#1#2{%
\def\temptext{#2}%
\edef\bibliotextoutput{\expandafter\gobbleeight\meaning\temptext}%
\global\advance\citationno by 1\citationlabel{#1}%
\ifmakebiblio%
    \edef\fileoutput{\write\biblio{\noindent\hbox to 0pt{\hss$[\the\citationno]$}\hskip 0.2em\bibliotextoutput\medskip}}%
    \fileoutput%
\fi}%
\def\cite#1{%
$[\citationref{#1}]$%
\ifmakebiblio%
    \edef\fileoutput{\write\biblio{#1}}%
    \fileoutput%
\fi%
}%
%
%
%

\let\mypar=\par


\def\raggedleft{\leftskip=0pt plus 1fil \parfillskip=0pt}


\font\lettrinefont=cmr10 at 28pt
\def\lettrine #1[#2][#3]#4%
{\hangafter -#1 \hangindent #2
\noindent\hskip -#2 \vtop to 0pt{
\kern #3 \hbox to #2 {\lettrinefont #4\hss}\vss}}

\font\mylettrinefont=cmr10 at 28pt
\def\mylettrine #1[#2][#3][#4]#5%
{\hangafter -#1 \hangindent #2
\noindent\hskip -#2 \vtop to 0pt{
\kern #3 \hbox to #2 {\mylettrinefont #5\hss}\vss}}


\edef\Pagetitle={Blank}

\headline={\hfil\Pagetitle\hfil}

\footline={\hfil\myfontdefault\folio\hfil}

\def\nextoddpage
{
\newpage%
\ifodd\pageno%
\else%
    \global\showpagenumflag = 0%
    \null%
    \vfil%
    \eject%
    \global\showpagenumflag = 1%
\fi%
}


\def\newchap#1#2%
{%
%
%
\global\advance\chapno by 1%
\resetcounters%
%
%
\newpage%
\ifodd\pageno%
\else%
    \global\showpagenumflag = 0%
    \null%
    \vfil%
    \eject%
    \global\showpagenumflag = 1%
\fi%
\global\newchapflag = 1%
\global\showpagenumflag = 1%
%
%
{\font\chapfontA=cmsl10 at 30pt%
\font\chapfontB=cmsl10 at 25pt%
\null\vskip 5cm%
{\chapfontA\raggedleft\hfil%
{%
\ifnum\chapno=0
    \phantom{%
    \ifinappendices%
        Annexe \alphanum\chapno%
    \else%
        \the\chapno%
    \fi}%
\else%
    \ifinappendices%
        Annexe \alphanum\chapno%
    \else%
        \the\chapno%
    \fi%
\fi%
}%
\par}%
\vskip 2cm%
{\chapfontB\raggedleft%
\lineskiplimit=0pt%
\lineskip=0.8ex%
\hfil #1\par}%
\vskip 2cm%
}%
\edef\Pagetitle{#2}%
%
%
\ifmaketdm%
    \def\temp{#2}%
    \def\tempbis{\nobreak}%
    \edef\chaptitle{\expandafter\gobbleeight\meaning\temp}%
    \edef\mynobreak{\expandafter\gobbleeight\meaning\tempbis}%
    \edef\textoutput{\write\tdm{\bigskip{\noexpand\mytdmchapfont\noindent\chinfo\ - \chaptitle\hfill\noexpand\folio}\par\mynobreak}}%
\fi%
\textoutput%
}


\def\newhead#1%
{%
\ifhmode%
    \mypar%
\fi%
\ifnum\headno=0%
\ifinappendices
    \nobreak\vskip -\lastskip%
    \nobreak\vskip .5cm%
\fi
\else%
    \nobreak\vskip -\lastskip%
    \nobreak\vskip .5cm%
\fi%
\nextheadno%
\ifmaketdm%
    \def\temp{#1}%
    \edef\sectiontitle{\expandafter\gobbleeight\meaning\temp}%
    \edef\textoutput{\write\tdm{\noindent{\noexpand\mytdmheadfont\quad\headinfo\ - \sectiontitle\hfill\noexpand\folio}\par}}%
    \textoutput%
\fi%
\font\headfontA=cmbx10 at 14pt%
{\headfontA\noindent #1.\hfil}%
\nobreak\vskip .5cm%
}%


\def\newsubhead#1%
{%
\ifhmode%
    \mypar%
\fi%
\ifnum\subheadno=0%
\else%
    \penalty\headpenalty\vskip .4cm%
\fi%
\nextsubheadno%
\ifmaketdm%
    \def\temp{#1}%
    \edef\subsectiontitle{\expandafter\gobbleeight\meaning\temp}%
    \edef\textoutput{\write\tdm{\noindent{\noexpand\mytdmsubheadfont\quad\quad\subheadinfo\ - \subsectiontitle\hfill\noexpand\folio}\par}}%
    \textoutput%
\fi%
\font\subheadfontA=cmsl10 at 12pt
{\subheadfontA\noindent\subheadinfo\ #1.\hfil}%
\nobreak\vskip .25cm %
}%

%
%


\font\mathromanten=cmr10
\font\mathromanseven=cmr7
\font\mathromanfive=cmr5
\newfam\mathromanfam
\textfont\mathromanfam=\mathromanten
\scriptfont\mathromanfam=\mathromanseven
\scriptscriptfont\mathromanfam=\mathromanfive
\def\mathroman{\fam\mathromanfam}


\font\sf=cmss12

\font\sansseriften=cmss10
\font\sansserifseven=cmss7
\font\sansseriffive=cmss5
\newfam\sansseriffam
\textfont\sansseriffam=\sansseriften
\scriptfont\sansseriffam=\sansserifseven
\scriptscriptfont\sansseriffam=\sansseriffive
\def\mathsf{\fam\sansseriffam}


\font\bftwelve=cmb12

\font\boldten=cmb10
\font\boldseven=cmb7
\font\boldfive=cmb5
\newfam\mathboldfam
\textfont\mathboldfam=\boldten
\scriptfont\mathboldfam=\boldseven
\scriptscriptfont\mathboldfam=\boldfive
\def\mathbf{\fam\mathboldfam}


\font\mycmmiten=cmmi10
\font\mycmmiseven=cmmi7
\font\mycmmifive=cmmi5
\newfam\mycmmifam
\textfont\mycmmifam=\mycmmiten
\scriptfont\mycmmifam=\mycmmiseven
\scriptscriptfont\mycmmifam=\mycmmifive

\def\hexa#1{\ifcase #1 0\or 1\or 2\or 3\or 4\or 5\or 6\or 7\or 8\or 9\or A\or B\or C\or D\or E\or F\fi}
\mathchardef\mathi="7\hexa\mycmmifam7B
\mathchardef\mathj="7\hexa\mycmmifam7C


\font\mymsbmten=msbm10 at 8pt
\font\mymsbmseven=msbm7 at 5.6pt
\font\mymsbmfive=msbm5 at 4pt
\newfam\mymsbmfam
\textfont\mymsbmfam=\mymsbmten
\scriptfont\mymsbmfam=\mymsbmseven
\scriptscriptfont\mymsbmfam=\mymsbmfive

\mathchardef\mybeth="7\hexa\mymsbmfam69
\mathchardef\mygimmel="7\hexa\mymsbmfam6A
\mathchardef\mydaleth="7\hexa\mymsbmfam6B


\def\placelabel[#1][#2]#3{{%
\setbox10=\hbox{\raise #2cm \hbox{\hskip #1cm #3}}%
\ht10=0pt%
\dp10=0pt%
\wd10=0pt%
\box10}}%


\newif\ifinproclaim%
\global\inproclaimfalse%
\def\proclaim#1{%
\medskip%
%
%
\bgroup%
\inproclaimtrue%
\setbox10=\vbox\bgroup\leftskip=0.8em\noindent{\bftwelve #1}\sf%
}

\def\endproclaim{%
\egroup%
\setbox11=\vtop{\noindent\vrule height \ht10 depth \dp10 width 0.1em}%
\wd11=0pt%
\setbox12=\hbox{\copy11\kern 0.3em\copy11\kern 0.3em}%
\wd12=0pt%
\setbox13=\hbox{\noindent\box12\box10}%
\noindent\unhbox13%
\egroup%
\medskip\ignorespaces%
}

\def\proclaim#1{%
\medskip%
\bgroup%
\inproclaimtrue%
\noindent{\bftwelve #1}%
\nobreak\medskip%
\sf%
}

\def\endproclaim{%
\mypar\egroup\penalty\proclaimpenalty\medskip\ignorespaces%
}

\def\noskipproclaim#1{%
\medskip%
\bgroup%
\inproclaimtrue%
\noindent{\bf #1}\nobreak\sl%
}

\def\endnoskipproclaim{%
\mypar\egroup\penalty\proclaimpenalty\medskip\ignorespaces%
}


\def\ninn{{n\in\Bbb{N}}}

\def\proof{{\noindent\bf Proof:\ }}

\def\remark{{\noindent\sl Remark:\ }}

\def\msup{\mathop{{\mathroman Sup}}}
\def\minf{\mathop{{\mathroman Inf}}}
\def\msf#1{{\mathsf #1}}

\def\qed{~$\square$}
\def\munion{\mathop{\cup}}
\def\minter{\mathop{\cap}}
\def\myitem#1{%
    \noindent\hbox to .5cm{\hfill#1\hss}
}

\catcode`\@=11
\def\Eqalign#1{\null\,\vcenter{\openup\jot\m@th\ialign{%
\strut\hfil$\displaystyle{##}$&$\displaystyle{{}##}$\hfil%
&&\quad\strut\hfil$\displaystyle{##}$&$\displaystyle{{}##}$%
\hfil\crcr #1\crcr}}\,}
\catcode`\@=12

\def\makeop#1{%
\global\expandafter\def\csname op#1\endcsname{{\mathroman #1}}}%

\def\makeopsmall#1{%
\global\expandafter\def\csname op#1\endcsname{{\mathroman{\lowercase{#1}}}}}%

\makeopsmall{ArcTan}%
\makeopsmall{ArcCos}%
\makeop{Arg}%
\makeop{Det}%
\makeop{Log}%
\makeop{Re}%
\makeop{Im}%
\makeop{Dim}%
\makeopsmall{Tan}%
\makeop{Ker}%
\makeopsmall{Cos}%
\makeopsmall{Sin}%
\makeop{Exp}%
\makeopsmall{Tanh}%
\makeop{Tr}%
\makeop{End}%
\makeop{Long}%
\makeop{Ch}%
\makeop{Exp}%
\makeop{Eval}%
\makeop{Lift}%
\makeop{Int}%
\makeop{Ext}%
\makeop{Aire}%
\makeop{Im}%
\makeop{Conf}%
\makeop{Exp}%
\makeop{Mod}%
\makeop{Log}%
\makeop{Ext}%
\makeop{Int}%
\makeop{Dist}%
\makeop{Aut}%
\makeop{Id}%
\makeop{GL}%
\makeop{SO}%
\makeop{Homeo}%
\makeop{Vol}%
\makeop{Ric}%
\makeop{Hess}%
\makeop{Euc}%
\makeop{Isom}%
\makeop{Max}%
\makeop{SW}%
\makeop{SL}%
\makeop{Long}%
\makeop{Fix}%
\makeop{Wind}%
\makeop{Diag}%
\makeop{dVol}%
\makeop{Symm}%
\makeop{Ad}%
\makeop{Diam}%
\makeop{loc}%
\makeopsmall{Sinh}%
\makeop{Len}%
\makeop{Length}%
\makeop{Min}%
\makeop{Area}%
\font\mycirclefont=cmsy7
\def\textcircle{{\raise 0.3ex \hbox{\mycirclefont\char'015}}}

\let\emph=\bf

\hyphenation{quasi-con-formal}

%
%

\ifmakebiblio%
    \openout\biblio=biblio.tex %
    {%
        \edef\fileoutput{\write\biblio{\bgroup\leftskip=2em}}%
        \fileoutput
    }%
\fi%

\newref{Almgren}{Almgren F. J. Jr., {\sl Plateau's problem. An invitation to varifold geometry.}, Student Mathematical Library, {\bf 13}, American Mathematical Society, Providence, RI, (2001)}
\newref{AndBarbBegZegh}{Andersson L., Barbot T., B\'e guin F., Zeghib A., Cosmological time versus CMC time in spacetimes of constant curvature}
\newref{Cabezas}{Cabezas-Rivas E., Miquel V., Volume-preserving mean curvature flow in the hyperbolic cpace, {\sl Indiana Univ. Math. J.} {\bf 56} (2007), no.5, 2061--2086}
\newref{CaffNirSprI}{Caffarelli L., Nirenberg L., Spruck J., The Dirichlet problem for nonlinear second-Order elliptic equations. I. Monge Amp\`e re equation, {\sl Comm. Pure Appl. Math.} {\bf 37} (1984), no. 3, 369--402}
\newref{CaffNirSprII}{Caffarelli L., Kohn J. J., Nirenberg L., Spruck J., The Dirichlet problem for nonlinear second-order elliptic equations. II. Complex Monge Amp\`e re, and uniformly elliptic, equations, {\sl Comm. Pure Appl. Math.} {\bf 38} (1985), no. 2, 209--252}
\newref{CaffNirSprIII}{Caffarelli L., Nirenberg L., Spruck J., The Dirichlet problem for nonlinear second-order elliptic equations. III. Functions of the eigenvalues of the Hessian, {\sl Acta Math.} {\bf 155} (1985), no. 3-4, 261--301}
\newref{CaffNirSprV}{Caffarelli L., Nirenberg L., Spruck J., Nonlinear second-order elliptic equations. V. The Dirichlet problem for Weingarten hypersurfaces, {\sl Comm. Pure Appl. Math.} {\bf 41} (1988), no. 1, 47--70}
\newref{ChenYau}{Cheng S. Y., Yau S. T., On the regularity of the Monge-Ampère equation $\opDet(\partial^2u\partial x_i\partial x_j)=F(x,u)$, \sl{Comm. Pure Appl. Math.} {\bf 30} (1977), no. 1, 41--68}
\newref{GallKapMard}{Gallo D., Kapovich M., Marden A., The monodromy groups of Schwarzian equations on closed Riemann surfaces, {\sl
Ann. of Math.} {\bf 2} 151 (2000), no. 2, 625--704}
\newref{Guan}{Guan B., The Dirichlet problem for Monge-Ampère equations in non-convex domains and spacelike hypersurfaces of constant Gauss curvature, {\sl Trans. Amer. Math. Soc.} {\bf 350} (1998), 4955--4971}
\newref{GuanSpruckI}{Guan B., Spruck J., The existence of hypersurfaces of constant Gauss curvature with prescribed boundary, {\sl J. Differential Geom.} {\bf 62} (2002), no. 2, 259--287}
\newref{GuanSpruckII}{Guan B., Spruck J., Szapiel M., Hypersurfaces of constant curvature in Hyperbolic space I, arXiv:0810.1779}
\newref{GuanSpruckIII}{Guan B., Spruck J., Hypersurfaces of constant curvature in Hyperbolic space II, arXiv:0810.1781}
\newref{Gut}{Guti\'errez C., {\sl The Monge-Amp\`e re equation}, Progress in Nonlinear Differential Equations and Their Applications, {\bf 44}, Birkh\"a user, Boston, (2001)}
\newref{HarvLaws}{Harvey R., Lawson H. B. Jr., Calibrated geometries, {\sl Acta Math.} {\bf 148} (1982), 47--157}
\newref{Huisken}{Huisken G., Contracting convex hypersurfaces in Riemannian manifolds by their mean curvature, {\sl Invent. Math.} {\bf 84}  (1986), no. 3, 463--480}
\newref{LabI}{Labourie F., Un lemme de Morse pour les surfaces convexes (French), {\sl Invent. Math.} {\bf 141} (2000), no. 2, 239--297}
\newref{LabII}{Labourie F., Probl\`e me de Minkowski et surfaces \`a courbure constante dans les vari\'e t\'e s hyperboliques (French), {\sl Bull. Soc. Math. France} {\bf 119} (1991), no. 3, 307--325}
\newref{LoftI}{Loftin J. C., Affine spheres and convex $\Bbb{RP}\sp n$-manifolds, {\sl Amer. J. Math.} {\bf 123} (2001), no. 2, 255--274}
\newref{LoftII}{Loftin J. C., Riemannian metrics on locally projectively flat manifolds, {\sl Amer. J. Math.} {\bf 124} (2002), no. 3, 595--609}
\newref{MazzPac}{Mazzeo R, Pacard P., Constant curvature foliations on asymptotically hyperbolic spaces, arXiv:0710.2298}
\newref{RosSpruck}{Rosenberg H., Spruck J., On the existence of convex hypersurfaces of constant Gauss curvature in hyperbolic space, {\sl J. Differential Geom.} {\bf 40} (1994), no. 2, 379--409}
\newref{SchKra}{Schlenker J. M., Krasnov K., The Weil-Petersson metric and the renormalized volume of hyperbolic $3$-manifolds, arXiv:0907.2590}
\newref{SmiFCS}{Smith G., Moduli of Flat Conformal Structures of Hyperbolic Type, arXiv:0804.0744}
\newref{SmiSLC}{Smith G., Special Lagrangian curvature, arXiv:math/0506230}
\newref{SmiGKM}{Smith G., Equivariant plateau problems, {\sl Geom. Dedicata} {\bf 140} (2009), 95--135}
\newref{SmiHPP}{Smith G., Hyperbolic Plateau problems, arXiv:math/0506231}
\newref{SmiPKS}{Smith G., Pointed $k$-surfaces, {\sl Bull. Soc. Math. France} {\bf 134} (2006), no. 4, 509--557}
\newref{SmiPHD}{Smith G., Probl\`e mes elliptiques pour des sous-vari\'e t\'e s Riemanniennes, Th\`ese, Orsay, 2004}

\ifmakebiblio%
    {\edef\fileoutput{\write\biblio{\egroup}}%
    \fileoutput}%
\fi%

%
%
%
\document
\myfontdefault
\global\chapno=1
\global\showpagenumflag=1
\def\Pagetitle{}
\null
\vfill
\def\centre{\rightskip=0pt plus 1fil \leftskip=0pt plus 1fil \spaceskip=.3333em \xspaceskip=.5em \parfillskip=0em \parindent=0em}%
\def\textmonth#1{\ifcase#1\or January\or Febuary\or March\or April\or May\or June\or July\or August\or September\or October\or November\or December\fi}
\font\abstracttitlefont=cmr10 at 14pt
{\abstracttitlefont\centre The Non-Linear Dirichlet Problem in Hadamard Manifolds\par}
\bigskip
{\centre Graham Smith\par}
\bigskip
{\centre \the\day\ \textmonth\month\ \the\year\par}
\bigskip
{\centre Centre de Recerca Matem\`atica,\par
Facultat de Ci\`encies, Edifici C,\par
Universitat Aut\`onoma de Barcelona,\par
08193 Bellaterra,\par
Barcelona,\par
SPAIN\par}
\bigskip
\noindent{\emph Abstract:\ } We proof existence theorems for the
Dirichlet problem for hypersurfaces of constant special Lagrangian
curvature in Hadamard manifolds. The first results are obtained
using the continuity method and approximation and then refined using
two iterations of the Perron method. The a-priori estimates used in
the continuity method are valid in any ambient manifold.
\bigskip
\noindent{\emph Key Words:\ } Dirichlet problem, special Lagrangian curvature, non-linear elliptic PDEs, Hadamard manifolds.
\bigskip
\noindent{\emph AMS Subject Classification:\ } 58E12 (35J25, 35J60, 53A10, 53C21, 53C42)
%
%
\par
\vfill
\nextoddpage
\global\pageno=1
\def\Pagetitle{\sl The Non-Linear Dirichlet Problem in Hadamard Manifolds}
\newhead{Introduction}
\noindent This paper treats the problem of finding immersed hypersurfaces with prescribed boundaries and curvature conditions in manifolds of strictly negative sectional curvature.
\medskip
\noindent This is an old geometric problem. The simplest version is Plateau's problem (see, for example, \cite{Almgren}), which requires minimal hypersurfaces with specified boundary. In this case, the curvature condition (minimality) is linear in terms of the shape operator of the immersion. A more general linear problem is that of finding hypersurfaces of constant mean curvature with specified boundary, for which a substantial litterature exists.
\medskip
\noindent The next interesting problem concerns hypersurfaces of constant Gaussian curvature. This is much harder, since Gaussian curvature is a non-linear function of the shape operator. Various results exist, using various different techniques (the following list is not exhaustive): constant Gaussian curvature surfaces which are graphs over hyperplanes in $\Bbb{R}^n$ are obtained using the continuity method by Caffarelli, Nirenberg and Spruck in \cite{CaffNirSprV} and by Guan in \cite{Guan}; hypersurfaces whose boundary is the boundary of a given convex set in $\Bbb{R}^n$ are the obtained by Spruck and Guan using the Perron method in \cite{GuanSpruckI}; graphs over horospheres and open subsets of the ideal boundary of $\Bbb{H}^n$ are obtained using the continuity method again by Rosenberg and Spruck in \cite{RosSpruck} while Guan and Spruck obtain more general results, again in $\Bbb{H}^n$ using a mixture of the continuity method and the Perron method in \cite{GuanSpruckII} and \cite{GuanSpruckIII}; a slightly different species of local existence results is obtained using the Implicit Function Theorem by Mazzeo and Pacard in \cite{MazzPac}; and finally more general graphs are obtained in $3$-dimensional Hadamard manifolds by Labourie in \cite{LabI} using the theory of pseudo-holomorphic curves and these results are further developed in the case of $\Bbb{H}^3$ by the author in \cite{SmiHPP} and \cite{SmiPKS}.
\medskip
\noindent Gaussian curvature constitutes the simplest non-linear curvature, but there are many other interesting examples. Existence results for constant curvature hypersurfaces for these different notions of curvature are interesting for various reasons and have varied applications. In general they yield dimensional reductions of geometric problems, since the space of constant curvature hypersurfaces is typically a finite dimensional manifold in contrast to, for example, the space of convex immersions, which is much more complicated. The following are a few applications of these types of results: in \cite{LabII}, Labourie uses constant Gaussian curvature surfaces to study the structure of $3$-dimensional hyperbolic manifolds; in \cite{SchKra}, Schlenker and Krasnov use constant mean curvature surfaces to study the relationship between the Teichmueller space of a surface and its moduli space of hyperbolic metrics; in \cite{SmiFCS}, the author uses constant curvature hypersurfaces to obtain geometric results concerning the structure of hyperbolic ends; in \cite{AndBarbBegZegh}, Andersson, Barbot, Beguin and Zeghib use constant mean curvature hypersurfaces to study flat, de-Sitter and anti de-Sitter spacetimes; and, in a similar vein, in \cite{LoftI} and \cite{LoftII}, Loftin uses the existence results \cite{ChenYau} of Cheng and Yau for solutions of the Monge-Amp\`ere equation to construct affine structures on convex projective manifolds.
\medskip
\noindent This paper studies hypersurfaces of constant special Lagrangian (SL) curvature. SL curvature was introduced by the author in \cite{SmiPHD} and \cite{SmiSLC}, and is defined in section $2$. As its name might suggest, it is closely related to the special Legendrian structure of the unitary bundle of the ambient manifold, and is derived from the theory of Calibrated Geometries developed by Harvey and Lawson in \cite{HarvLaws}. The interest of SL curvature is twofold. Firstly, like Gaussian curvature, it is intimately related to convexity, and thus provides a natural tool in the study of convex problems (which is employed in the case of hyperbolic ends and flat conformal structures in \cite{SmiFCS}), and secondly it is regular, in the sense that a sequence of hypersurfaces of constant SL curvature only degenerates in one simple way, which can often then be excluded by geometric considerations. This simple property makes it much easier to obtain existence results than for other non linear curvatures. Crucially, Gaussian curvature exhibits this property only when the dimension of the ambient manifold is equal to $3$.
\medskip
\noindent SL curvature is only defined for convex immersed hypersurfaces and depends on an angle parameter, $\theta\in[0,n\pi/2[$. We thus denote it by $R_\theta$. We only concern ourselves with the case where $\theta\in[(n-1)\pi/2,n\pi/2[$, since, in this case, $R_\theta$ interacts wells with convexity (more precisely, it vanishes along the boundary of the cone of positive definite symmetric two forms and is positive in its interior). The case where $\theta>(n-1)\pi/2$ is more regular and, in general, results will be obtained for this case and then extended to the case where $\theta=(n-1)\pi/2$ by compactness. However, when $\theta=(n-1)\pi/2$, $R_\theta$ has the simplest form and the most interesting geometric properties. For example, when $n=2$:
$$
R_\theta = K^{1/2},
$$
\noindent where $K$ is the Gaussian curvature, and when $n=3$:
$$
R_\theta = (K/H)^{1/2},
$$
\noindent where $H$ is the mean curvature. In higher dimensions, $R_{(n-1)\pi/2}$ has a more complicated expression, but still exhibits the same properties: specifically, sequences of constant curvature hypersurfaces degenerate in exactly the same way as constant Gaussian curvature surfaces do in $3$-dimensional ambient manifolds (see \cite{LabI}). It is for this reason amongst others that SL curvature should be considered as an alternative higher dimensional generalisation of Gaussian curvature, in analogy to the way in which the symplectic structure can be considered as an alternative higher dimensional generalisation of $2$-dimensional volume.
\medskip
\noindent Throughout the rest of the introduction, we shall use the rescaled SL curvature $\hat{R}_\theta$ given by:
$$
\hat{R}_\theta = \opTan(\theta/n)R_\theta.
$$
\noindent This is chosen so that the rescaled SL curvature of a horosphere in hyperbolic space equals $1$, which is as it should be. Nevertheless, to save on multiplicative factors, throughout the rest of the paper, $R_\theta$ will be used.
\medskip
\noindent The first result we obtain is an existence theorem for graphs of given SL curvature over hypersurfaces. Let $M$ be an $(n+1)$-dimensional Hadamard manifold. Choose $\theta\in[(n-1)\pi/2,n\pi/2[$. Let $H\subseteq M$ be a smooth, convex hypersurface such that:
$$
\hat{R}_\theta(H) = R_0,
$$
\noindent Where $R_0$ is constant. Let $\Omega\subseteq H$ be a bounded open subset. Let $\hat{\Sigma}\subseteq M$ be a convex immersed hypersurface such that $\partial\hat{\Sigma}=\partial\Omega=:\Gamma$ and such that $\hat{R}_\theta(\hat{\Sigma})\geqslant R_1$ in the weak sense, where:
$$
R_1 \leqslant 1.
$$
\noindent Using the continuity method, we obtain:
\proclaim{Theorem \nextprocno}
\noindent Suppose that $\hat{\Sigma}$ is a graph over $\Omega$ and that $\Gamma$ is strictly convex as a subset of $M$ with respect to the outward pointing normal to $\Gamma$ in $\hat{\Sigma}$. If $\theta>(n-1)\pi/2$, then, for all $r\in[R_0,R_1]$, there exists an immersed hypersurface $\Sigma_r\subseteq M$ such that:
\medskip
\myitem{(i)} $\Sigma_r$ is $C^0$ and $C^\infty$ in its interior;
\medskip
\myitem{(ii)} $\partial\Sigma_r=\Gamma$;
\medskip
\myitem{(iii)} $\Sigma_r$ is a graph over $\Omega$ lying below $\hat{\Sigma}$; and
\medskip
\myitem{(iv)} $\hat{R}_\theta(\Sigma_r)=r$.
\medskip
\noindent Moreover, the same result holds for $\theta=(n-1)\pi/2$ provided that, in addition $\hat{\Sigma}$ is $\epsilon$-convex, for some $\epsilon>0$.
\endproclaim
\proclabel{TheoremDirichletI}
\remark The hypotheses of this theorem are satisfied when the norm of the second fundamental form of $H$ is small with respect to that of $\hat{\Sigma}$ and the normal of $\hat{\Sigma}$ is sufficiently bounded away from $TH$ along $\Gamma$. Explicitely, if $\hat{\Sigma}$ is $\epsilon$-convex, if the norm of the second fundamental form of $H$ is bounded above by $\delta$ and if the angle between the normal to $\hat{\Sigma}$ and $TH$ is bounded below by $\theta$ along $\Gamma$, then the hypotheses are satisfied provided that:
$$
\epsilon\opSin(\theta) - \delta > 0.
$$
\remark The proof of this theorem closely follows the approach towards the non-linear Dirichlet problem for functions over open subsets of $\Bbb{R}^n$ used in \cite{CaffNirSprI}, \cite{CaffNirSprIII} and \cite{CaffNirSprV}. The key innovation here is the use of geometric tools for the construction of appropriate barrier functions. Special Lagrangian curvature turns out to be particularly amenable to this analysis. However, given the generality of the results of \cite{CaffNirSprI}, \cite{CaffNirSprIII} and \cite{CaffNirSprV}, one would hope that these geometric constructions (especially Lemma \procref{LemmaBoundaryLowerNormalBounds}) may also be of use in the study of other curvatures.
\medskip
\remark The continuity method consists of two key stages: local deformation and compactness. Only the local deformation stage requires the ambient manifold to have strictly negative sectional curvature. Importantly, the compactness stage follows from a-priori estimates that are valid in any manifold.
\medskip
\noindent The special case where $M$ is $(n+1)$-dimensional hyperbolic space, $\Bbb{H}^{n+1}$, and $H$ is a totally geodesic hypersurface is interesting in itself. In particular, the hypersurfaces thus obtained are unique:
\goodbreak
\proclaim{Theorem \nextprocno}
\noindent Choose $\theta\in[(n-1)\pi/2,n\pi/2]$. Let $H\subseteq\Bbb{H}^{n+1}$ be a totally geodesic hypersurface. Let $\Omega\subset H$ be a bounded open subset. Let $\hat{\Sigma}\subseteq\Bbb{H}^{n+1}$ be a convex hypersurface which is a graph over $\Omega$ such that $\partial\hat{\Sigma}=\partial\Omega$ and:
$$
\hat{R}_\theta(\hat{\Sigma}) \geqslant R_1,
$$
\noindent in the weak sense, where $R_1\leqslant 1$.  If $\theta>(n-1)\pi/2$, then, for all $r\in[0,R_1]$, there exists a unique immersed hypersurface $\Sigma_r\subseteq M$ such that:
\medskip
\myitem{(i)} $\Sigma_r$ is $C^0$ and $C^\infty$ in its interior;
\medskip
\myitem{(ii)} $\partial\Sigma_r=\Gamma$;
\medskip
\myitem{(iii)} $\Sigma_r$ is a graph over $\Omega$ lying below $\hat{\Sigma}$; and
\medskip
\myitem{(iv)} $\hat{R}_\theta(\Sigma_r)=r$.
\medskip
\noindent Moreover, the same result holds for $\theta=(n-1)\pi/2$ provided that, in addition, $\hat{\Sigma}$ is $\epsilon$-convex, for some $\epsilon>0$.
\endproclaim
\proclabel{TheoremDirichletIA}
\noindent The following particular case is also of interest:
\proclaim{Theorem \nextprocno}
\noindent Let $H\subseteq\Bbb{H}^{n+1}$ be a totally geodesic hypersurface. Let $\Omega\subset H$ be a bounded open subset. If $\partial\Omega$ is $1$-convex, then, for all $\theta\in[(n-1)\pi/2,n\pi/2]$ and for all $r\in[0,1]$, there exists a unique immersed hypersurface $\Sigma_r\subseteq M$ such that:
\medskip
\myitem{(i)} $\Sigma_r$ is $C^0$ and $C^\infty$ in its interior;
\medskip
\myitem{(ii)} $\partial\Sigma_r=\Gamma$;
\medskip
\myitem{(iii)} $\Sigma_r$ is a graph over $\Omega$ lying below $\hat{\Sigma}$; and
\medskip
\myitem{(iv)} $\hat{R}_\theta(\Sigma_r)=r$.
\endproclaim
\proclabel{TheoremDirichletIB}
\remark Theorem \procref{TheoremDirichletIB} illustrates a general feature of hyperbolic space: that the curvature of horospheres, which is equal to $1$, provides a threshold for geometric results. This becomes particularly evident in the study of curvature flows (c.f. \cite{Huisken} and \cite{Cabezas}), and constitutes an important distinction between hyperbolic space and Euclidean space. In both spaces, the curvature of totally geodesic hypersurfaces, which is equal to $0$, forms one threshold, but, in Euclidean space, horospheres coincide with totally geodesic hypersurfaces, and so the horospheric threshold is absorbed into the totally geodesic one.
\medskip
\noindent We use the Perron method to generalise Theorem \procref{TheoremDirichletI}. Let $M$ be an $(n+1)$-dimensional Hadamard manifold of sectional curvature bounded above by $-1$. Let $\Sigma$ be a smooth, convex immersed hypersurface in $M$. Let $\Omega\subseteq\Sigma$ be an open subset and let $\hat{\Sigma}$ be a convex immersed hypersurface in $M$ which is a graph over the extended normal of $\Sigma$:
\goodbreak
\proclaim{Theorem \nextprocno}
\noindent Choose $\theta\in[(n-1)\pi/2,n\pi/2[$. Choose $0\leqslant R_0<R_1\leqslant 1$. Suppose that $R_\theta(\Sigma)\leqslant R_0$ and $R_\theta(\hat{\Sigma})\geqslant R_1$ in the weak sense. If $\theta>(n-1)\pi/2$, then, for all $r\in[R_0,R_1]$, there exists an immersed hypersurface $\Sigma_r$ in $M$ such that:
\medskip
\myitem{(i)} $\Sigma_r$ is a graph over $\hat{N}\Omega$;
\medskip
\myitem{(ii)} $\Sigma_r$ lies below $\hat{\Sigma}$ as a graph over $\hat{N}\Omega$;
\medskip
\myitem{(iii)} $\Sigma_r$ is smooth away from the boundary; and
\medskip
\myitem{(iv)} $\hat{R}_\theta(\Sigma_r) = r$.
\medskip
\noindent If $\theta=(n-1)\pi/2$, then the same result holds provided that, in addition, $\hat{\Sigma}$ is $\epsilon$-convex for some $\epsilon>0$.
\endproclaim
\proclabel{TheoremDirichletII}
\remark The notation is explained in section $4$.
\medskip
\remark The use of the Perron method in the proof of Theorem \procref{TheoremDirichletII} is inspired by the work \cite{GuanSpruckI} of Guan and Spruck on hypersurfaces of constant Gaussian curvature in $\Bbb{R}^n$.
\medskip
\remark A simple modification of Theorem \procref{TheoremDirichletII} yields results in the case where $\Sigma$ is a closed hypersurface in a negatively curved manifold. In particular, this yields the existence results of \cite{SmiFCS} in much greater generality than is required in that paper.
\medskip
\noindent The Perron method may iterated, this time using Theorem \procref{TheoremDirichletII} to provide the local existence results, and we obtain:
\proclaim{Theorem \nextprocno}
\noindent Let $M$ be an $(n-1)$-dimensional manifold of negative sectional curvature bounded above by $-1$. Let $N\subseteq M$ be a compact, convex immersed submanifold. Suppose that the diameter of immersions homotopic to $N$ is bounded below by $\epsilon>0$. Choose $\theta\in[(n-1)\pi/2[$ and $r\in]0,1[$. Then:
\medskip
\myitem{(i)} If $\theta>(n-1)\pi/2$, then there exists a smooth, convex, immersed submanifold $N_{r,\theta}\in M$, isotopic to $N$ such that:
$$
R_\theta(N_{r,\theta}) = r.
$$
\myitem{(ii)} If $\theta=(n-1)\pi/2$, then the same result holds provided that, in addition, $N$ is not homeomorphic to the sphere bundle $S^{n-1}\times S^1$.
\endproclaim
\proclabel{TheoremDirichletIII}
\remark the hypotheses are satisfied if, for example, $N$ is homotopically non-trivial and $M$ is compact or geometrically finite without cusps.
\medskip
\noindent An interesting application of this result concerns previous work of the author in \cite{SmiGKM}. Using the notation of that paper, we obtain:
\goodbreak
\proclaim{Theorem \nextprocno}
\noindent Let $M$ be a compact, three-dimensional manifold of negative sectional curvature bounded above by $-1$. Let $\Sigma$ be a compact Riemann surface of hyperbolic type. Let $\alpha:\pi_1(\Sigma)\rightarrow\pi_1(M)$ be a homomorphism such that:
\medskip
\myitem{(i)} $\alpha(\pi_1(\Sigma))$ is non-elementary; and
\medskip
\myitem{(ii)} the second Stiefel-Whitney class of $\alpha$ vanishes.
\medskip
\noindent Then, for all $k\in]0,1[$, there exists a smooth immersion $i:\Sigma\rightarrow M$ of constant Gaussian curvature equal to $k$ such that $i_*=\alpha$.
\endproclaim
\proclabel{TheoremDirichletIIIA}
\noindent We thus obtain algebraic conditions for the existence of constant Gaussian curvature hypersurfaces immersed in a compact three dimensional manifold of strictly negative curvature. Importantly, the construction of \cite{SmiGKM} (which is closely related to the work \cite{GallKapMard} of Gallo, Kapovich and Marden) yields infinitely many isotopy inequivalent convex immersions. Thus, by Theorem \procref{TheoremDirichletIII}, we obtain infinitely many families of distinct constant Gaussian curvature hypersurfaces in any homotopy class satisfying the algebraic hypotheses of Theorem \procref{TheoremDirichletIIIA}. There is no a-priori reason not to expect such degeneration even in higher dimensions which leads one to speculate on the structure of the space of solutions.
\medskip
\noindent The paper is arranged as follows:
\medskip
\myitem{(i)} various background concepts are introduced and studied in section $2$. Special Lagrangian curvature is introduced and it is shown that it is defined in terms of a homogeneous, concave function. The properties of convex subsets of Riemannian manifolds are studied in detail, and it is shown how mollifiers may be used to produce convex sets with certain desired properties;
\medskip
\myitem{(ii)} in section $3$, which is the most innovative part of the paper, geometric constructions are used to produce the barrier functions required to determine a-priori $C^2$ bounds on hypersurfaces of constant special Lagrangian curvature which are graphs for the case when $\theta>(n-1)\pi/2$. Using these estimates and the Continuity method, we prove Theorem \procref{TheoremDirichletI} for this case, and the case when $\theta=(n-1)\pi/2$ is proven by taking limits. This section concludes with the proofs of Theorems \procref{TheoremDirichletIA} and \procref{TheoremDirichletIB};
\medskip
\myitem{(iii)} in section $4$, inspired by \cite{LabI}, we introduce extended normals and graphs over extended normals. We then use the Perron method along with Theorem \procref{TheoremDirichletI} to prove Theorem \procref{TheoremDirichletII}; and
\medskip
\myitem{(iv)} in section $5$, we introduce the concept of pseudo-immersions as a compactification of the space of convex immersions. These are used as an important tool in the proof of Theorem \procref{TheoremDirichletIII}, which is carried out using the Perron method and Theorem \procref{TheoremDirichletII}.
\medskip
\noindent The author would like to thank Harold Rosenberg for drawing his attention to \cite{GuanSpruckI}, which has been a signficant impetus for the current paper. The author is also grateful to the CRM in Barcelona for providing the conditions necessary for carrying out this research.
\goodbreak
\newhead{Preliminaries}
\newsubhead{Immersed Submanifolds and Special Lagrangian Curvature}
\noindent Let $M$ be a smooth Riemannian manifold. An {\bf immersed submanifold} is a pair $\Sigma=(S,i)$ where $S$ is a smooth manifold and $i:S\rightarrow M$ is a smooth immersion. An {\bf immersed hypersurface} is an immersed submanifold of codimension $1$. We give $S$ the unique Riemannian metric $i^*g$ which makes $i$ into an isometry. We say that $\Sigma$ is {\bf complete} if and only if the Riemannian manifold $(S,i^*g)$ is.
\medskip
\noindent Let $UM$ be the unitary bundle of $M$ (i.e the bundle of unit vectors in $TM$. In the cooriented case (for example, when $I$ is convex), there exists a unique exterior normal vector field $\msf{N}$ over $i$. We denote $\hat{\mathi}=\msf{N}$ and call it the {\bf Gauss lift} of $i$. Likewise, we call the manifold $\hat{\Sigma}=(S,\hat{\mathi})$ the {\bf Gauss lift} of $\Sigma$.
\medskip
\noindent The special Lagrangian curvature, which is only defined for strictly convex immersed hypersurfaces, is defined as follows. Denote by $\opSymm(\Bbb{R}^n)$ the space of symmetric matrices over $\Bbb{R}^n$. We define $\Phi:\opSymm(\Bbb{R}^n)\rightarrow\Bbb{C}^*$ by:
\headlabel{SpecialLagrangianCurvature}
$$
\Phi(A) = \opDet(I+iA).
$$
\noindent Since $\Phi$ never vanishes and $\opSymm(\Bbb{R}^n)$ is simply connected, there exists a unique analytic function $\tilde{\Phi}:\opSymm(\Bbb{R}^n)\rightarrow\Bbb{C}$
such that:
$$
\tilde{\Phi}(I) = 0,\qquad e^{\tilde{\Phi}(A)} = \Phi(A)\qquad\forall A\in\opSymm(\Bbb{R}^n).
$$
\noindent We define the function $\opArcTan:\opSymm(\Bbb{R}^n)\rightarrow(-n\pi/2,n\pi/2)$ by:
$$
\arctan(A) = \opIm(\tilde{\Phi}(A)).
$$
\noindent This function is trivially invariant under the action of $O(\Bbb{R}^n)$. If $\lambda_1, ..., \lambda_n$ are the eigenvalues of $A$, then:
$$
\opArcTan(A) = \sum_{i=1}^n\opArcTan(\lambda_i).
$$
\noindent For $r>0$, we define:
$$
\opSL_r(A) = \opArcTan(r^{-1}A).
$$
\noindent If $A$ is positive definite, then $SL_r$ is a strictly decreasing function of $r$. Moreover, $SL_\infty=0$ and $SL_0=n\pi/2$. Thus, for all $\theta\in]0,n\pi/2[$, there exists a unique $r>0$ such that:
$$
SL_r(A) = \theta.
$$
\noindent We define $R_\theta(A) = r$. $R_\theta$ is also invariant under the action of $O(n)$ on the space of positive definite, symmetric matrices.
\medskip
\noindent Let $M$ be an oriented Riemannian manifold of dimension $n+1$. Let $\Sigma=(S,i)$ be a strictly convex, immersed hypersurface in $M$. For $\theta\in]0,n\pi/2[$, we define $R_\theta(\Sigma)$ (the {\emph $\theta$-special Lagrangian curvature} of $\Sigma$) by:
$$
R_\theta(\Sigma) = R_\theta(A_\Sigma),
$$
\noindent where $A_\Sigma$ is the shape operator of $\Sigma$.
\goodbreak
\newsubhead{Properties of Special Lagrangian Curvature}
\noindent $R_\theta$ is an analytic homogenous function of order $1$. Importantly:
\proclaim{Lemma \nextprocno}
\noindent For all $\theta$, $R_\theta$ is a concave function over the set of positive definite symmetric matrices.
\endproclaim
\proclabel{LemmaConcavityOfSLCurvature}
\remark This property is necessary for the application of the Perron method. In the following proof, we explicitely determine the second derivative. However, a simpler, more geometric argument may also be employed. Indeed, the function $\opSL_r$ is concave. Moreover $\opSL_r(A)\geqslant\theta$ if and only if $R_\theta(A)\geqslant r$. It follows that $R_\theta^{-1}([r,+\infty[)$ is convex for all $r>0$, and the result follows since the Hessian of a homogeneous function is, up to a factor, the second second fundamental form of its level sets.
\medskip
\proof Define the function $\sigma$ over the space of symmetric matrices by:
$$
\sigma(A) = \opArg(\opDet(\opId + iA)).
$$
\noindent Trivially:
$$\matrix
D\sigma_A(M) \hfill&= \opTr(\mu^{-1}M),\hfill\cr
D^2\sigma_A(M,M) \hfill&= -2\opTr(\mu^{-1}AM\mu^{-1}M),\hfill\cr
\endmatrix$$
\noindent where $\mu = \opId + A^2$. Choose $\theta\in]0,n\pi/2[$. Define the function $r$ over the space of symmetric matrices such that:
$$
\sigma(r(A)^{-1}A) = \theta.
$$
\noindent Define $\mu_r$ and $\phi_r$ by:
$$
\mu_r = \opId + r^{-2}A^2,\qquad\phi_r = \opTr(\mu_r A).
$$
\noindent Using the chain rule and the formula for $D\sigma$ and $D^2\sigma$ yields:
$$
D^2r_A(M,M) = \frac{-2}{r\phi_r}\opTr(\mu_r^{-1}A\tilde{M}\mu_r^{-1}\tilde{M}),
$$
\noindent where:
$$
\tilde{M} = M - \frac{1}{\phi_r}\opTr(\mu_r^{-1}M)A.
$$
\noindent Thus, when $A$ is positive definite, for all $M$:
$$
D^2r_A(M,M) \leqslant 0,
$$
\noindent The result follows.\qed
\medskip
\noindent For $\theta>(n-1)\pi/2$, $R_\theta(A)$ approximates the smallest eigenvalue of $A$:
\proclaim{Lemma \nextprocno}
\noindent Let $\lambda_1(A)$ denote the smallest eigenvalue of the matrix $A$. For all $\theta\in ]0,n\pi/2[$, there exists $K_1$ such that:
$$
R_\theta(A) \geqslant K_1\lambda_1(A).
$$
\noindent For all $\theta>(n-1)\pi/2$, there exists $K_2$ such that:
$$
R_\theta(A) \leqslant K_2\lambda_1(A).
$$
\endproclaim
\proclabel{LemmaSLCIsLowestEval}
\remark Observe that the second relation is no longer valid in the case where $\theta\leqslant(n-1)\pi/2$.
\medskip
\remark In particular, when $\theta>(n-1)\pi/2$, constant special Lagrangian curvature yields uniform lower bounds on the principal curvatures.
\medskip
\proof Let $K_1=R_\theta(\opId)$. Then:
$$
R_\theta(A) = \lambda_1(A)R_\theta(A/\lambda_1(A)) \geqslant \lambda_1(A)R_\theta(\opId) = K_1\lambda_1(A).
$$
\noindent The first result follows. Let $\lambda_1=\lambda_1(A)$ and $r=R_\theta(A)$. Then:
$$
r^{-1}\lambda_1 \geqslant \opArcTan(r^{-1}\lambda_1) \geqslant \theta - (n-1)\pi/2.
$$
\noindent Thus, if $K_2:=(\theta-(n-1)\pi/2)^{-1}<\infty$, then:
$$
R_\theta(A) = r \leqslant K_2\lambda_1 = K_2\lambda_1(A).
$$
\noindent The second result follows.\qed
\medskip
\noindent The case where $\theta=(n-1)\pi/2$ is of particular interest. Here $R_\theta$ has the simplest form and the most interesting geometric properties. For example, when $n=2$:
$$
R_{\pi/2} = K^{1/2},
$$
\noindent where $K$ is the Gaussian curvature, and when $n=3$:
$$
R_\pi = (K/H)^{1/2},
$$
\noindent where $H$ is the mean curvature. However, the case where $\theta>(n-1)\pi/2$ is more regular. In the sequel, results for $\theta=(n-1)\pi/2$ are obtained by first treating this case, and then taking limits.
\medskip
\noindent It is important to note that, although $R_\theta$ is more appealing geometrically, $\opSL_r$ is analytically simpler. In the sequel, results for $R_\theta$ constant will often be proven for $\opSL_r$ constant, which is trivially equivalent.
\goodbreak
\newsubhead{Convex Conditions and Convex Immersions}
\noindent Let $M$ be an $(n+1)$-dimensional Riemannian manifold. let $TM$ and $UM\subseteq TM$ be the tangent and unitary bundles respectively over $M$. Let $\pi:TM\rightarrow M$ be the canonical projection.
\medskip
\noindent Let $\opSymm^+(\Bbb{R}^n)$ denote the space of symmetric, positive definite matrices over $\Bbb{R}^n$. Let $X$ be an open subset of $\opSymm^+(\Bbb{R}^n)$. We say that $X$ defines a homogeneous convex property if and only if:
\medskip
\myitem{(i)} $X$ is convex;
\medskip
\myitem{(ii)} $X$ is invariant under the action of $\opSO(n)$; and
\medskip
\myitem{(iii)} $X$ is homogeneous in the sense that, for all $\lambda\in[1,\infty[$, $\lambda X\subseteq X$.
\medskip
\noindent Let $K\subseteq M$ be a convex set. We say that $K$ posseses the property $X$ if and only if, for all $p\in\partial K$, and for every supporting normal $\msf{N}_p$ to $K$ at $p$, there exists a smooth hypersurfaces $\Sigma$ such that:
\medskip
\myitem{(i)} $\Sigma$ is an exterior tangent to $K$ at $p$;
\medskip
\myitem{(ii)} $\msf{N}_p$ is the outward pointing normal to $\Sigma$ at $p$; and
\medskip
\myitem{(iii)} if $A$ is the second fundamental form of $\Sigma$ at $p$, then $A\in X$.
\proclaim{Lemma \nextprocno}
\noindent Let $X$ define a homogeneous convex property. Let $K,K'\subseteq M$ be convex sets. If $K$ and $K'$ possess the property $X$, then so does $K\minter K'$.
\endproclaim
\proclabel{LemmaIntersectionIsStillEpsilonConvex}
\proof It suffices to check the condition at $p\in \partial K\minter\partial K'$. Let $\msf{N}_p$ and $\msf{N}'_p$ be supporting normals to $K$ and $K'$ respectively at $p$. If $\msf{N}_p=\msf{N}'_p$, then $K\minter K'$ possesses the property $X$ in the direction of $\msf{N}_p=\msf{N}'_p$. Likewise, if $\msf{N}_p=\msf{N}'_p$, then $K\minter K'$ is a single point. We thus assume that they are distinct and not colinear.
\medskip
\noindent Let $\pi:\opSymm^+(\Bbb{R}^n)\rightarrow\opSymm^+(\Bbb{R}^{n-1})$ be the projection defined by restriction to the subspace. Define $X'=\pi(X)$. Trivially, $X'$ is open, convex and homogeneous.
\medskip
\noindent Let $\Sigma_p$ and $\Sigma'_p$ be smooth convex hypersurfaces at $p$ as in the definition of possession of property $X$. $\Sigma_p$ and $\Sigma'_p$ are transverse at $p$. Let $A$ and $A'$ be the second fundamental forms at $p$ of $\Sigma_p$ and $\Sigma'_p$ respectively. Define $\Gamma = \Sigma_p\minter\Sigma'_p$. Near $p$, $\Gamma$ is a smooth submanifold. For $s,t> 0$ such that $s+t=1$, define $\msf{N}_{s,t}$ by:
$$
\msf{N}_{s,t} = s\msf{N}_p + t\msf{N}'_p.
$$
\noindent Let $A_\Gamma$ be the second fundamental form of $\Gamma$ at $p$. $A_\Gamma$ depends on a choice of normal vector to $\Gamma$ at $p$:
$$\matrix
A_\Gamma(\msf{N}_p)\hfill&= \pi(A) \hfill&\in X'\hfill\cr
A_\Gamma(\msf{N}_p')\hfill&= \pi(A') \hfill&\in X'\hfill\cr
\endmatrix$$
\noindent Thus, for all $s,t>0$ such that $s+t=1$:
$$
A_\Gamma(\frac{1}{\|\msf{N}_{s,t}\|}\msf{N}_{s,t}) = \frac{1}{\|\msf{N}_{s,t}\|}(sA_\Gamma(\msf{N}_p) + tA_\Gamma(\msf{N}'_p))\in X'.
$$
\noindent Thus, for all $s,t>0$ such that $s+t=1$, there exists an immersed hypersurface $\Sigma_{s,t}$ such that:
\medskip
\myitem{(i)} $\Sigma_{s,t}$ is an exterior tangent to $K\minter K'$ at $p$;
\medskip
\myitem{(ii)} $\msf{N}_{s,t}/\|\msf{N}_{s,t}\|$ is the outward pointing normal to $\Sigma_{s,t}$ at $p$; and
\medskip
\myitem{(ii)} if $A_{s,t}$ is the second fundamental form of $\Sigma_{s,t}$ at $p$, then $A_{s,t}\in X$.
\medskip
\noindent Since the set of supporting normals to $K\minter K'$ at $p$ is a convex set whose boundary is contained in the union of the sets of supporting normals to $K$ and $K'$ at $p$, the result follows.\qed
\medskip
\noindent Let $K\subseteq M$ be a convex set. For $\epsilon>0$, we say that $K$ is $\epsilon$-convex if and only if, for every $p\in\partial K$, for every supporting normal $\msf{N}_p$ to $K$ at $p$, and for every $0<\epsilon'<\epsilon$, there exists a smooth convex hypersurface $\Sigma$ such that:
\medskip
\myitem{(i)} $\Sigma$ is an exterior tangent to $K$ at $p$;
\medskip
\myitem{(ii)} $\msf{N}_p$ is the outward pointing normal to $\Sigma$; and
\medskip
\myitem{(iii)} the second fundamental form of $\Sigma$ at $p$ is bounded below by $\epsilon'\opId$.
\medskip
\noindent This is a homogeneous convex property, and is thus preserved under the intersection of convex sets.
\medskip
\noindent Choose $\theta\in[0,n\pi/2[$ and $r>0$. We say that $R_\theta(\partial K)\geqslant r$ in the weak sense if and only if, for every $p\in\partial K$, for every supporting normal $\msf{N}_p$ to $K$ at $p$, and for every $0<r'<r$, there exists a smooth convex hypersurface $\Sigma$ such that:
\medskip
\myitem{(i)} $\Sigma$ is an exterior tangent to $K$ at $p$;
\medskip
\myitem{(ii)} $\msf{N}_p$ is the outward pointing normal to $\Sigma$; and
\medskip
\myitem{(iii)} $R_\theta(\Sigma)(p)\geqslant r'$.
\medskip
\noindent Since is is also a homogeneous convex property, Lemma \procref{LemmaIntersectionIsStillEpsilonConvex} yields:
\proclaim{Lemma \nextprocno}
\noindent Choose $\theta\in[0,n\pi/2[$ and $r>0$. Let $K,K'\subseteq M$ be convex sets. If $R_\theta(\partial K),R_\theta(\partial K')\geqslant r$ in the weak sense, then $R_\theta(\partial K\minter K')\geqslant r$ in the weak sense.
\endproclaim
\proclabel{LemmaIntersectionStillHasCurvatureBoundedBelow}
\noindent If $K\subseteq M$ is a convex set, and $U\subseteq\partial K$ is an open subset of the boundary, let $\Cal{N}(U)$ denote the set of supporting normals to $K$ over $U$. Let $(N,\partial N)$ be a compact $n$-dimensional manifold with boundary. A convex immersion of $N$ into $M$ is a pair $(\varphi,\hat{\varphi})$ where:
\medskip
\myitem{(i)} $\varphi:N\rightarrow M$ and $\hat{\varphi}:N\rightarrow UM$ are $C^{0,1}$ mappings such that $\pi\circ\hat{\varphi}=\varphi$; and
\medskip
\myitem{(ii)} for every $p\in N$, there exists a convex set $K\subseteq M$ such that $\varphi(p)\in\partial K$ and neighbourhoods $U\subseteq N$ and $V\subseteq\partial K$ of $p$ and $\varphi(p)$ respectively such that $\hat{\varphi}$ restricts to a homeomorphism from $U$ to $\Cal{N}(U)$.
\medskip
\noindent In the sequel, we will denote the convex immersion simply by $\varphi$. $\epsilon$-convex immersions are defined in an analogous manner.
\proclaim{Lemma \nextprocno}
\medskip
\noindent Suppose that $n\geqslant 2$ and $\partial N\neq\emptyset$. Let $K\subseteq M$ be a convex subset. Let $\varphi:N\rightarrow M$ be a convex immersion. Suppose there exists an open subset $U\subseteq N$ and a point $p\in\partial K$ such that:
\medskip
\myitem{(i)} $\varphi(U)=\partial K\setminus\left\{p\right\}$; and
\medskip
\myitem{(ii)} $\varphi(\partial U)=\left\{p\right\}$.
\medskip
\noindent Then $\hat{\varphi}$ defines a homeomorphism between $N$ and $\Cal{N}(\partial K)$.
\endproclaim
\proclabel{LemmaRemoveableSingularity}
\proof Choose $p\in\partial U$. Let $U'\subseteq N$ and $K'\subseteq M$ be a neighbourhood of $p$ in $N$ and a convex subset of $M$ respectively as in the definition of convex immersions. The complement of $\varphi^{-1}(\left\{p\right\})$ in $U'$ has only one connected component. However:
$$
\partial(U\minter U') \subseteq \varphi^{-1}(\left\{p\right\})\minter U'.
$$
\noindent Since $U$ is not contained in $\varphi^{-1}(\left\{p\right\})$ it follows that:
$$
U'\setminus\varphi^{-1}(\left\{p\right\})\subseteq U.
$$
\noindent In particular $N=U\munion U'$. $\hat{\varphi}$ therefore defines a covering map from $N$ to $\Cal{N}(\partial K)$. Since $n\geqslant 2$ and the latter is homeomorphic to an $n$-dimensional sphere, $\hat{\varphi}$ is a homeomorphism, and the result follows.\qed
\medskip
\noindent Finally we recall the following Geometric Maximum Principle.
\proclaim{Lemma \nextprocno}
\noindent Let $M$ be a Riemannian manifold and let $\Sigma=(S,i)$ and $\Sigma'=(S',i')$ be $C^0$ convex, immersed hypersurfaces in $M$. For $\theta\in]0,n\pi/2[$, let $R_\theta$ and $R_\theta'$ be the $\theta$-special Lagrangian curvatures of $\Sigma$ and $\Sigma'$ respectively. If $p\in S$ and $p'\in S'$ are such that $q=i(p)=i'(p')$, and $\Sigma'$ is an interior tangent to $\Sigma$ at $q$, then:
$$
R_\theta(p) \leqslant R'_\theta(p').
$$
\endproclaim
\proclabel{LemmaGeomMaxPrinc}
\proof See \cite{SmiFCS}.\qed
\goodbreak
\newsubhead{Distance Functions}
\noindent Let $M$ be a Riemannian manifold. Let $p\in M$ be a point. Let $\Sigma\subseteq M$ be a strictly concave, smooth, immersed hypersurface passing through $p$. Let $d_\Sigma$ be the signed distance in $M$ to $\Sigma$.
\proclaim{Lemma \nextprocno}
\noindent There exists a neighbourhood, $U$ of $p$ in $M$ such that $d_\Sigma$ is concave in $U$.
\endproclaim
\proof For $t\in\Bbb{R}$, let $\Sigma_t=d_\Sigma^{-1}(\left\{t\right\})$. For $U$ a sufficiently small neighbourhood of $p$ and for $t$ small, the intersection of $\Sigma_t$ with $U$ is smooth and concave. For all $t$, let $\msf{N}_t$ and $II_t$ be the unit normal vector and the second fundamental form respectively of $\Sigma_t$. Then, for all $X$ tangent to $\Sigma_t$:
$$\matrix
\opHess(d_\Sigma)(X,X) \hfill&= A_t(X,X), \hfill\cr
\opHess(d_\Sigma)(X,\msf{N}) \hfill&= 0, \hfill\cr
\opHess(d_\Sigma)(\msf{N},\msf{N}) \hfill&= 0. \hfill\cr
\endmatrix$$
\noindent The result follows.\qed
\medskip
\noindent Let $U\subseteq M$ be an open subset. Let $\Sigma$ be a hypersurface in $M$. Let $\msf{N}$ and $II^\Sigma$ be the unit normal and the second fundamental form respectively of $\Sigma$. Let $\opHess^\Sigma$ be the Hessian for smooth functions defined on $\Sigma$. Trivially, we obtain:
\proclaim{Lemma \nextprocno}
\noindent Let $\phi:U\rightarrow\Bbb{R}$ be smooth. Then:
$$
\opHess^\Sigma(f) = \opHess(f) - df(\msf{N})II^\Sigma.
$$
\endproclaim
\proclabel{LemmaConcavityPreserved}
\proclaim{Corollary \nextprocno}
\noindent Suppose that $\Sigma$ is convex and $\langle\nabla\phi,\msf{N}\rangle\geqslant 0$. If $\phi$ is concave as a function over $U$, then it is also concave as a function over $\Sigma\minter U$.
\endproclaim
\proclabel{CorollaryConcaveFunction}
\remark Let $\Sigma$ be a convex surface with smooth boundary. Let $p\in\partial\Sigma$ and $\Sigma'$ be a concave surface tangent to $\partial\Sigma$ at $p$ such that $\Sigma$ locally lies in its exterior near $p$. Let $d_{\Sigma'}$ now be the distance function to $\Sigma'$. We see that $d_{\Sigma'}$ acts as a barrier for Laplacians derived from the Hessian of $\Sigma$. This will play a central r\^ole later.
\medskip
\noindent Let $d_p$ be the distance to $p$ in $M$. For all $r$, let $\Sigma_r$ be the sphere of radius $r$ about $p$. We recall:
\proclaim{Lemma \nextprocno}
\noindent If $X$ is tangent to $\Sigma_r$ and $\msf{N}$ is the unit exterior normal to $\Sigma_r$, then, near $p$:
$$\matrix
\opHess(d_p)(X,X) \hfill&=r^{-1}(1+O(r^2))\langle X,X\rangle,\hfill\cr
\opHess(d_p)(X,\msf{N}) \hfill&=0,\hfill\cr
\opHess(d_p)(\msf{N},\msf{N}) \hfill&=0.\hfill\cr
\endmatrix$$
\endproclaim
\noskipproclaim{Corollary \nextprocno}
$$
\opHess(d_p^2/2) = \langle\cdot,\cdot\rangle + O(r_p^2).
$$
\endnoskipproclaim
\proclabel{CorollaryApproximateHessianOfDistance}
\goodbreak
\newsubhead{Regularising Convex Sets}
\noindent We recall the definition of mollifiers:
\proclaim{Definition \nextprocno}
\noindent Let $M$ be a Riemannian manifold. A mollifier of $M$ is a smooth, positive function $\varphi:TM\rightarrow[0,+\infty[$ such that:
\medskip
\myitem{(i)} for all $p\in M$:
$$
\int_{T_pM}\varphi\opdVol_p = 1,
$$
\noindent where $\opdVol_p$ is the volume form of $T_pM$;
\medskip
\myitem{(ii)} $\varphi(v_p)=0$ for $\|v_p\|>1$; and
\medskip
\myitem{(iii)} $\varphi$ is preserved by parallel transport of $M$.
\endproclaim
\noindent We construct mollifiers as follows. Let $\psi:[0,\infty[\rightarrow[0,\infty[$ be a smooth, positive function such that:
$$
t\leqslant 1/2\Rightarrow \psi(t)=1,\qquad t\geqslant 1\Rightarrow\psi(t)=0.
$$
\noindent Let $\lambda>0$ be a positive constant and define $\varphi:TM\rightarrow[0,\infty[$ by:
$$
\varphi(v_p) = \lambda\psi(\|v_p\|).
$$
\noindent $\varphi$ is trivially preserved by parallel transport. If $\lambda$ is chosen such that the integral of $\varphi$ over any (and thus every) tangent space is equal to $1$, then $\varphi$ is a mollifier.
\medskip
\noindent If $\varphi$ is a mollifier, we define $(\varphi_\epsilon)_{\epsilon>0}:TM\rightarrow[0,+\infty[$ by:
$$
\varphi_\epsilon(v_p) = \epsilon^{-n}\varphi(\epsilon^{-1}v_p).
$$
\noindent Using mollifiers, we obtain:
\proclaim{Lemma \nextprocno}
\noindent Let $M$ be a Riemannian manifold. Choose $\theta\in[(n-1)\pi/2,n\pi/2[$ and $r>0$. Let $\Sigma\subseteq M$ be a compact, convex immersed hypersurface such that $R_\theta(\Sigma)\geqslant r$ in the weak sense. If $\theta>(n-1)\pi/2$, then, for all $\delta>0$ there exists a smooth, convex hypersurface $\Sigma'$ (which may be chosen arbitrarily close to $\Sigma$ in the $C^{0}$ sense) such that:
$$
R_\theta(\Sigma') \geqslant r-\delta.
$$
\noindent If $\theta = (n-1)\pi/2$, the the same result holds provided that the second fundamental form of $\Sigma$ is bounded below in the weak sense.
\endproclaim
\proclabel{LemmaSmoothingConvexSets}
\remark Mollification preserves homogeneous convex conditions up to a small error. This is the content of the proof.
\medskip
\proof Let $\varphi$ be a mollifier of $M$. Let $\opExp:TM\rightarrow M$ be the exponential map of $M$. We work locally and therefore assume that there exists a unique geodesic between any two points, $x,y\in M$. Let $\tau_{y,x}$ be parallel transport from $x$ to $y$ along this geodesic.
\medskip
\noindent Define $f:M\rightarrow\Bbb{R}$ by:
$$
f(p) = d(p,\Sigma).
$$
\noindent This function is convex in a small neighbourhood of $\Sigma$. We restrict to this neighbourhood for the rest of the proof. $f$ is a locally $C^{1,1}$ function away from $\Sigma$. In particular, $\opHess(f)$ is measurable and bounded in every compact subset of the complement of $\Sigma$. For $\epsilon>0$, define $f_\epsilon:M\rightarrow\Bbb{R}$ by:
$$
f_\epsilon(p) = \int_{T_pM}(f\circ\opExp)(V_p)\varphi_\epsilon(V_p)\opdVol_p.
$$
\noindent Trivially, $(f_\epsilon)_{\epsilon>0}\rightarrow f$ in the $C^1$ sense as $\epsilon\rightarrow 0$. It remains to show that the second derivative of $f_\epsilon$ has the desired properties for $\epsilon$ sufficiently small.
\medskip
\noindent We construct an approximation for the Hessian. For $\epsilon>0$, define $A_\epsilon\in\Gamma(\opSymm(TM))$ by:
$$
A_\epsilon(p) = \int_{T_pM}(\opExp_p^*\opHess f)\varphi_\epsilon\opdVol_p.
$$
\noindent For $t>0$, let $\Sigma_t$ be the level hypersurface of $f$ with value $t$. Let $\delta_1>0$ be small. There exists $T_0>0$ such that, for $t<T_0$, $R_\theta(\Sigma_t)>r-\delta_1$ (this is the stage that requires the supplementary condition when $\theta=(n-1)\pi/2)$. Let $p\in M$ be such that $0<f(p)<T_0$. Let $X_p\in T_pM$ be a unit vector orthogonal to $\nabla f$ at $p$. Let $\gamma:\Bbb{R}\rightarrow M$ be a unit speed geodesic such that $\partial_t\gamma(0)=X_p$. For all $t$, define the vector field $X_t$, such that, for all $V\in T_{\gamma(t)}M$:
$$
X_t(\opExp_{\gamma(t)}(V)) = D\opExp_{\gamma(t)}\cdot\partial_t\gamma(t).
$$
\noindent For $V$ in $T_pM$, we define $c_V:\Bbb{R}\rightarrow M$ by:
$$
c_V(t) = (\opExp\circ\tau_{\gamma(t),\gamma(0)})(V).
$$
\noindent By Taylor's Theorem:
$$\matrix
(f\circ c_V)(t) \hfill&= (f\circ c_V)(0) + df(\partial_t c_V)(0)t + \int_0^t (t-s)\partial_t^2(f\circ c_V)(s)ds\hfill\cr
&= (f\circ c_V)(0) + df(\partial_t c_V)(0)t\hfill\cr
&\qquad + \int_0^t (t-s)\opHess(f)(\partial_t c_V,\partial_t c_V)ds + \int_0^t (t-s)df(\nabla_{\partial_t c_V}\partial_t c_V)ds.\hfill\cr
\endmatrix$$
\noindent Let $\eta>0$ be small. Trivially, $c_V\rightarrow\gamma$ in the $C^\infty$ sense as $\|V\|\rightarrow 0$. Thus, since $\opHess(f)$ is bounded, there exists $\epsilon_0>0$ such that, for $\|V\|<\epsilon_0$ and for all $t\in]-r,r[$:
$$
(f\circ c_V)(t) \geqslant (f\circ c_V)(0) + df(\partial_t c_V)(0)t + \int_0^t (t-s)(\opHess(f)-\eta)_{c_V(s)}(X_s,X_s)ds.
$$
\noindent For $\epsilon<\epsilon_0$, since $\varphi\opdVol_p$ is invariant under parallel transport:
$$\matrix
(f_\epsilon\circ\gamma)(t)\hfill&= \int_{T_pM}\varphi_\epsilon(V)(f\circ\opExp\circ\tau_{\gamma(t),\gamma(0)})(V)\opdVol_p\hfill\cr
&=\int_{T_pM}\varphi_\epsilon(V)(f\circ c_V)(t)\opdVol_p\hfill\cr
&\geqslant \int_{T_pM}\varphi_\epsilon(V)((f\circ c_V)(0) + tdf(\partial_t c_V)(0))\opdVol_p\hfill\cr
&\qquad + \int_0^t (t-s)(\opHess(f)-\eta)_{c_V(s)}(X_s,X_s)ds\opdVol_p.\hfill\cr
\endmatrix$$
\noindent However:
$$\matrix
(f_\epsilon\circ\gamma)(0) \hfill&= \int_{T_pM}\varphi_\epsilon(V)(f\circ\opExp\circ\tau_{\gamma(t),\gamma(0)})(V)\opdVol_p|_{t=0}\hfill\cr
&= \int_{T_pM} \varphi_\epsilon(V)(f\circ c_V)(0)\opdVol_p\hfill\cr
\endmatrix$$
\noindent Moreover:
$$\matrix
df_\epsilon(\partial_t\gamma(0)) \hfill&= \partial_t\int_{T_pM}\varphi_\epsilon(V)(f\circ\opExp\circ\tau_{\gamma(t),\gamma(0)})(V)\opdVol_p|_{t=0}\hfill\cr
&=\int_{T_pM}\varphi_\epsilon(V)(\partial_t(f\circ c_V)(t)|_{t=0})\opdVol_p\hfill\cr
&=\int_{T_pM}\varphi_\epsilon(V)df(\partial_t c_V(0))\opdVol_p.\hfill\cr
\endmatrix$$
\noindent Finally, by definition:
$$
A_\epsilon(\gamma(s))(\partial_t \gamma,\partial_t \gamma) = \int_{T_pM}\varphi_\epsilon(V)\opHess(f-\eta)_{c_V(s)}(X_s,X_s)ds)\opdVol_p
.$$
\noindent Thus, for $\epsilon<\epsilon_0$:
$$
(f_\epsilon\circ\gamma)(t)\geqslant (f_\epsilon\circ\gamma)(0) + df_\epsilon(\partial_t\gamma(0))t + \int_0^t (t-s)(A_\epsilon-\eta)(\gamma(s))(\partial_t\gamma,\partial_t\gamma)ds.
$$
\noindent Consequently:
$$
\opHess(f_\epsilon)(p) \geqslant A_\epsilon - \eta.
$$
\noindent Let $E_p$ be the orthogonal complement of $\nabla f_\epsilon$ at $p$. Let $E$ be the distribution obtained by parallel transport of $E_p$ along geodesics leaving $p$. Let $A|_E$ be the restriction of $\opHess(f)$ to $E$. Since $f_\epsilon$ tends to $f$ in the $C^1$ sense, and since $\opHess(f)$ is bounded, for $q$ sufficiently close to $p$.
$$
R_\theta(A(q)|_E) \geqslant r - \delta/2.
$$
\noindent However, $R_\theta$ is a concave function. Thus, for $\epsilon$ sufficiently small:
$$
R_\theta(A_\epsilon(p)|_E) \geqslant r - \delta/2.
$$
\noindent And so:
$$
R_\theta(\|\nabla f_\epsilon\|^{-1}\opHess(f_\epsilon)|_E)\geqslant r -\delta.
$$
\noindent Since these estimates may be calculated locally uniformly, the result follows by taking an appropriate level subset of $f_\epsilon$, for $\epsilon$ sufficiently small.\qed
\goodbreak
\newhead{The Continuity Method}
\newsubhead{First Order Control}
\noindent Let $M$ be a Riemannian manifold. Let $\opExp$ denote the exponential mapping of $M$. Let $H$ be a smooth convex hypersurface. Let $\msf{N}_H$ be the outward pointing unit normal over $H$. Let $\Omega\subseteq H$ be an open set. We will say that a $C^{0,1}$ hypersurface $\Sigma$ is a graph over $\Omega$ if and only if there exists a $C^{0,1}$ function $f:\overline{\Omega}\rightarrow[0,+\infty[$ and a homeomorphism $\varphi:\overline{\Omega}\rightarrow\Sigma$ such that:
\medskip
\myitem{(i)} $f$ vanishes along $\partial\Omega$ (i.e. $\partial\Sigma=\partial\Omega$); and
\medskip
\myitem{(ii)} for all $p\in\Omega$:
$$
\varphi(p) = \opExp_p(f(p)\msf{N}_H(p)).
$$
\noindent We refer to $f$ as the graph function of $\Sigma$. Consider the family of graphs over $\Omega$. We define the partial order ``$>$'' on this family such that if $\Sigma$ and $\Sigma'$ are two graphs over $\Omega$ and $f$ and $f'$ are their respective graph functions, then:
$$
\Sigma > \Sigma' \Leftrightarrow f(p) > f'(p)\text{ for all }p\in\Omega.
$$
\noindent Since $\partial\Omega$ is smooth, for all $p\in\partial\Omega$, the set of supporting hyperplanes to $\partial\Omega$ at $p$ is parametrised by $\Bbb{R}$. Supporting hyperplanes may be locally considered as graphs over $\Omega$, and we obtain an analogous partial order on this set, which we also denote by $>$.
\medskip
\noindent Choose $\theta\in](n-1)\pi/2,n\pi/2[$. Suppose now that $R_\theta(H)=R_0$ where $R_0\geqslant 0$ is constant. Let $\hat{\Sigma}$ be a $C^{0,1}$ convex hypersurface which is a graph over $\Omega$. Suppose moreover that $R_\theta(\hat{\Sigma})\geqslant R_1$ in the weak sense. Let $R_0<r_1<r_2<...<r_\infty<R_1$ be a sequence of positive real numbers and let $(\Sigma_n)_{\ninn}$ be a sequence of graphs over $\Omega$ such that:
\medskip
\myitem{(i)} for all $\ninn$, $\Sigma_n$ is a smooth convex hypersurface such that $R_\theta(\Sigma_n)=r_n$;
\medskip
\myitem{(ii)} for all $\ninn$, $\Sigma_n<\hat{\Sigma}$; and
\medskip
\myitem{(iii)} for all $n>m$, $\Sigma_n>\Sigma_m$.
\proclaim{Lemma \nextprocno}
\noindent There exists a $C^{0,1}$ convex hypersurface with boundary $\Sigma_0$ which is $C^\infty$ in its interior such that:
\medskip
\myitem{(i)} $\Sigma_0$ is a graph over $\Omega$;
\medskip
\myitem{(ii)} $\hat{\Sigma}>\Sigma_0$; and
\medskip
\myitem{(iii)} The sequence of graph functions $(f_n)_\ninn$ converges to $f_0$ in the $C^{0,\alpha}$ sense over $\overline{\Omega}$ and in the $C^\infty_\oploc$ sense over $\Omega$.
\medskip
\noindent Moreover, if $\Sigma_0$ is smooth up to the boundary:
\medskip
\myitem{(iv)} for every $p\in\partial\Sigma$, $T_p\hat{\Sigma} > T_p\Sigma_0$.
\endproclaim
\proclabel{LemmaCompactnessOfConvexGraphs}
\proof For all $n$, let $f_n$ be the graph function of $\Sigma_n$. $(f_n)_\ninn$ is uniformly bounded above by the graph function of $\hat{\Sigma}$. Since $(f_n)_\ninn$ is strictly increasing, there exists $f_\infty$ to which this sequence converges pointwise. For all $n\in\Bbb{N}\munion\left\{\infty\right\}$, define $U_n$ by:
$$
U_n = \left\{ \opExp_p(t\msf{N}_H(p))\text{ s.t. }p\in\overline{\Omega}\text{ and }0\leqslant t\leqslant f_n(p)\right\}.
$$
\noindent Trivially, for all $i<j$, $U_i\subseteq U_j$, and:
$$
U_\infty = \munion_{i=1}^\infty U_n.
$$
\noindent Since $U_n$ is convex (away from $H$) for all $n$, so is $U_\infty$. Moreover, the supporting hyperplanes of $U_\infty$ are transverse to the normal geodesics leaving $H$. Indeed, let $p\in\partial U_\infty$ be a point where the supporting hyperplane is not transverse to the normal geodesic leaving $H$. Since $\hat{\Sigma}_0>\Sigma_n$ for all $n$, $p\notin\partial\Omega$. Let $q\in\Omega$ be such that $\opExp(f_\infty(q)\msf{N}(q))=p$.  Let $\gamma$ be the geodesic segment joining $q$ to $p$. $\gamma$ lies inside $U_\infty$. Moreover, it is tangent to $\partial U_\infty$ at $p$. Consequently, it lies in the boundary of $U_\infty$ and thus defines a continuous path in $\partial U_\infty$ from $p$ to $q$ which does not intersect $\partial\Omega$. This is absurd.
\medskip
\noindent Since $\partial U_\infty$ is the graph of $f_\infty$, it follows that $f_\infty$ is $C^{0,1}$ and that $(f_n)_\ninn$ converges to $f_\infty$ in the $C^{0,\alpha}$ sense. This proves $(i)$ and the first half of $(iii)$.
\medskip
\noindent Let $\epsilon>0$. Let $p\in\Omega$ be such that $d(p,\partial\Omega)>\epsilon$. For all $n\in\Bbb{N}\munion\left\{\infty\right\}$, define $p_n\in\Sigma_n$ by:
$$
p_n = \opExp(f_n(p)\msf{N}_H(p)).
$$
\noindent Trivially $(p_n)_\ninn$ converges to $p_\infty$. Choose $\epsilon>0$. For all $n$, let $B_n$ be the ball of radius $\epsilon$ about $p_n$ in $\Sigma_n$. By Theorem $1.4$ of \cite{SmiSLC}, for $\epsilon$ sufficiently small, there exists an immersed hypersurface $\Sigma'_\infty$ containing $p_\infty$ such that $(B_n,p_n)$ subconverges to $(\Sigma'_\infty,p_\infty)$ in the $C^\infty$ pointed Cheeger/Gromov sense. Trivially, $\Sigma'_\infty\subseteq\Sigma_\infty$. Thus every subsequence of $(f_n)_\ninn$ subconverges in the $C^\infty_\oploc$ sense to $f_\infty$. This proves the second half of $(iii)$.
\medskip
\noindent In order to prove $(ii)$, suppose the contrary. Then $\Sigma_\infty$ intersects $\hat{\Sigma}$ non trivially at an interior point, $p$, say. Since $\hat{\Sigma}\geqslant\Sigma_\infty$, $\Sigma_\infty$ is an interior tangent to $\hat{\Sigma}$ at this point, which is absurd by the Geometric Maximum Principal (Lemma \procref{LemmaGeomMaxPrinc}).
\medskip
\noindent $(iv)$ follows in a similar manner from the Geometric Maximum Principal, and this concludes the proof.\qed
\medskip
\noindent Let $\msf{N}_n$ be the unit normal vector field over $\Sigma_n$. Since $\partial\Omega$ is smooth, for all $n\in\Bbb{N}\munion\left\{\infty\right\}$, $\Sigma_n$ has only one supporting hyperplane at any point of $\partial\Omega$. $\msf{N}_n$ therefore extends continuously to $\partial\Omega$.
\medskip
\noindent By compactness, there exists $\epsilon>0$ such that, for all $p\in\partial\Omega$ and for all $q$ in $B_\epsilon(p)$, there exists a unique geodesic joining $p$ to $q$. For such $p$ and $q$, let $\tau_{q,p}$ denote parallel transport from $p$ to $q$ along this geodesic. The following uniform modulus of continuity is of subtle importance in the sequel:
\proclaim{Lemma \nextprocno}
\noindent There exists a continuous function $\delta:[0,\infty[\rightarrow[0,\infty[$ such that for all $n\in\Bbb{N}$, for all $p\in\partial\Omega$ and for all $q\in B_\epsilon(p)$, if $\msf{N}'$ is a supporting normal to $\Sigma_n$ at $p$, then:
$$
\|\tau_{q,p}\msf{N}_n(p) - \msf{N}'\| < \delta(d(p,q)).
$$
\endproclaim
\proclabel{LemmaUniformModulusOfContinuity}
\proof The normal to a convex set is continuous whereever it is uniquely defined. Likewise, for a sequence of convex sets converging towards a limit, the normal converges at any point in the limit where it is uniquely defined. Uniformity of convergence follows by compactness, and the result follows.\qed
\newsubhead{Constructing Barriers I}
\noindent Let $M$ be an $(n+1)$-dimensional manifold. Let $H\subseteq M$ be a smooth convex hypersurface such that:
$$
R_\theta(H) = R_0,
$$
\noindent where $R_0$ is constant. Let $\Omega\subseteq H$ be a bounded open subset of $H$ with smooth boundary. Let $\hat{\Sigma}\subseteq M$ be a convex hypersurface such that $\partial\hat{\Sigma}=\partial\Omega=:\Gamma$ and such that $R_\theta(\hat{\Sigma})\geqslant R_1$ in the weak sense. Suppose that $\Gamma$ is strictly convex as a subset of $M$ with respect to the outward pointing normal to $\Gamma$ in $\hat{\Sigma}$.
\medskip
\noindent For $\msf{N}$ a normal vector to $\Gamma$, let $A_\Gamma(\msf{N})$ be the second fundamental form of $\Gamma$ in the direction of $\msf{N}$. Thus, if $X$ and $Y$ are vector fields tangent to $\Gamma$:
$$
A_\Gamma(\msf{N})(X,Y) = -\langle \nabla_X Y, \msf{N}\rangle.
$$
\proclaim{Proposition \nextprocno}
\noindent Let $p\in\Gamma$. Let $\msf{N}_0$ be the outward pointing normal to $H$ at $p$. Let $\msf{N}_1$ be the outward pointing normal to $\Gamma$ in $\hat{\Sigma}$ at $p$. For $s,t\in [0,1]$ such that $s+t=1$, $A_\Gamma(s\msf{N}_0 + t\msf{N}_1)$ is strictly positive definite.
\endproclaim
\proclabel{PropositionConvexityCondition}
\remark Unit vectors colinear to vectors of the form $s\msf{N}_1 + t\msf{N}_1$ for $s,t\geqslant 0$ will be said to lie between the outward normals of $H$ and $\partial\hat{\Sigma}$.
\medskip
\remark In particular, $\partial\Omega$ is strictly convex as a subset of $M$ with respect to any normal lying between the normals of $H$ and $\partial\hat{\Sigma}$
\medskip
\proof By definition, $A_\Gamma(\msf{N}_0)$ and $A_\Gamma(\msf{N}_1)$ are both strictly positive definite. The result follows by convexity of the set of positive definite quadratic forms.\qed
\proclaim{Corollary \nextprocno}
\noindent The normal to $\hat{\Sigma}$ at $p$ points above $H$.
\endproclaim
\proof Otherwise, if $\hat{\msf{N}}$ is the outward pointing normal to $\hat{\Sigma}$ at $p$, then $-\hat{\msf{N}}$ lies between the normals of $H$ and $\partial\hat{\Sigma}$. $\Gamma$ is therefore strictly concave with respect to $\hat{\msf{N}}$. This is absurd, since $\hat{\Sigma}$ is strictly convex.\qed
\medskip
\noindent Let $p\in\partial\Omega$. Let $\msf{N}_p$ be a normal vector to $\Gamma$ at $p$ lying between the outward normals of $H$ and $\partial\hat{\Sigma}$. Let $\lambda_1,...,\lambda_{n-1}$ be the eigenvalues of $A_\Gamma(\msf{N}_p)$. We define $\mu(p,\msf{N}_p,r,\theta)\in]0,\infty]$ by:
$$
\mu(p,\msf{N}_p,r,\theta) = \sup\left\{ m > 0\text{ s.t. }SL_r(\lambda_1,...,\lambda_{n-1},m)<\theta\right\}.
$$
\noindent $\mu$ is continuous in $p$, $\msf{N}_p$, $r$ and $\theta$. Suppose that $\mu(p,\msf{N}_p,r,\theta)=+\infty$. We aim to construct a barrier for hypersurfaces of constant special Lagrangian curvature equal to $r$ whose boundary is $\Gamma$ and whose normal at $p$ is $\msf{N}_p$.
\medskip
\noindent Near $p$, let $\msf{N}'_H$ be the parallel transport of the upwards pointing normal of $H$ at $p$. Define the set $\hat{\Omega}\subseteq M$ near $p$ by:
$$
\hat{\Omega} = \left\{\opExp(t\msf{N}'_H(q))\text{ s.t. }t\in-]\epsilon,\epsilon[,q\in\Omega\right\}.
$$
\noindent $\hat{\Omega}$ may be considered as the solid vertical cylinder over $\Omega$. Define the real valued function $d_H$ and $d_\Omega$ over a neighbourhood of $p$ by:
$$
d_H(q) = d(q,H),\qquad d_\Omega(q) = d(q,\partial\hat{\Omega}),
$$
\noindent where $d_\Omega$ is chosen to be positive inside $\Omega$. Observe that $(\nabla d_H,\nabla d_\Omega)$ forms an orthonormal basis of the space of normal vectors to $\Gamma$ at $p$. Let $K>0$ be such that $\nabla d_H - K\nabla d_\Omega$ is parallel to $\msf{N}_p$. Define the real valued function $\Phi_0$ in a neighbourhood of $p$ by:
$$
\Phi_0 = d_H - Kd_\Omega.
$$
\noindent Since $\msf{N}_p$ lies between the normals to $H$ and $\partial\hat{\Sigma}$ at $p$, by Proposition \procref{PropositionConvexityCondition}, there exists a strictly convex hypersurface $H'$ which is a strict exterior tangent to $\Gamma$ at $p$ such that:
$$
T_pH' = T_p\partial\Omega\oplus\langle\msf{N}_p\rangle.
$$
\noindent Let $A_{H'}$ be the second fundamental form of $H'$. We may choose $H'$ such that $\|A_{H'}\|$ is arbitrarily small at $p$. We define $d_{H'}$ by:
$$
d_{H'}(q) = d(q,H').
$$
\noindent For any two functions $f$ and $g$ with non-colinear derivatives at $p$, define the the $(n-2)$-dimensional distribution $E(f,g)$ near $p$ by:
$$
E(f,g) = \langle \nabla f,\nabla g\rangle^\perp.
$$
\noindent Let $e_1,...,e_{n-1}$ be an orthonormal basis for $T_p\Gamma$ with respect to which $A_\Gamma(\msf{N}_p)$ is
diagonal. We extend this to a local frame in $TM$ such that, at $p$, for all $X$ and all $i$:
$$
\nabla_X e_i = -\opHess(d_{H'})(e_i,X)\nabla d_{H'} - (1+K^2)^{-1}\opHess(\Phi_0)(e_i,X)\nabla\Phi_0.
$$
\noindent Define the distribution $E$ near $p$ to be the span of $e_1,...,e_{n-1}$.
\proclaim{Proposition \nextprocno}
\noindent If $D$ represents the Grassmannian distance between two $(n-2)$-dimensional subspaces then:
$$
D(E,E(\Phi_0,d_{H'})) = O(d_p^2).
$$
\endproclaim
\proclabel{PropositionSubspacesCloseI}
\proof At $p$:
$$
\langle\nabla d_{H'},\nabla \Phi\rangle = 0.
$$
\noindent Thus, for every vector $X$ at $p$, and for all $i$, by definition of $e_i$:
$$\matrix
&\langle\nabla_X e_i,\nabla d_{H'}\rangle \hfill&= -\opHess(d_{H'})(e_i,X)\hfill\cr
& &= -\langle e_i,\nabla_X\nabla d_{H'}\rangle\hfill\cr
\Rightarrow\hfill&X\langle e_i,\nabla d_{H'}\rangle \hfill&= 0.\hfill\cr
\endmatrix$$
\noindent Likewise:
$$
X\langle e_i,\nabla \Phi\rangle = 0.
$$
\noindent The result now follows.\qed
\medskip
\noindent For any smooth function $f$ and any non-negative function $l$, we define $\opSL_r'(f,l,E)$ by:
$$
\opSL'_r(f,l,E) = \sum_{i=1}^{n-1}\opArcTan(\frac{1}{r\sqrt{1+l^2}}\lambda_i(f,E)),
$$
\noindent where $(\lambda_i(f,E))_{i\leqslant i\leqslant (n-1)}$ are the eigenvalues of the restriction of
$\opHess(f)$ to $E$.
\proclaim{Proposition \nextprocno}
\noindent Let $f$ be such that $f(p),\nabla f(p)=0$ and the restriction of $\opHess(f)$ to $H'$ is positive definite. There exists a function $x$ such that $x(p),\opHess(x)(p)=0$ and:
$$
\opSL'_r(\Phi_0 + x(d_{H'} - f),\|\nabla\Phi_0\|,E) \leqslant \theta - \pi/2  + O(d_p^2).
$$
\endproclaim
\proclabel{PropositionSecondOrderVariation}
\proof By definition of $\Phi_0$ and $\msf{N}_p$, at $p$:
$$
\opSL_r'(\Phi_0,\|\nabla\Phi_0\|,E)\leqslant\theta - \pi/2.
$$
\noindent The Hessian of $xf$ vanishes at $p$. Likewise, the Hessian of the second order term $xd_{H'}$ vanishes on $(\nabla d_{H'})^\perp$ at $p$ and thus so does its restriction to $E$. It follows that $x(d_{H'}-f)$ does not affect $\opSL'_r$ at $p$. Thus, for all $x$, at $p$:
$$
\opSL_r'(\Phi_0 + x(d_{H'} - f),\|\nabla\Phi_0\|,E)\leqslant\theta - \pi/2.
$$
\noindent Denote $l=\sqrt{1+\|\nabla\Phi_0\|^2}$. For $1\leqslant i\leqslant n-1$, define $\mu_i$ by:
$$
\mu_i = \frac{(rl)^{-1}}{1+(rl)^{-2}\lambda_i^2}.
$$
\noindent Define $A$ and $B$ by:
$$
A = \sum_{i=1}^{n-1}\mu_if_{;ii}, B=\sum_{i=1}^{n-1}\mu_id_{H';ii}.
$$
\noindent Since $f_{;ij}$ is positive definite, $A>0$. Likewise, since $H'$ is concave, $B<0$. Define the vectors $X$ and $Y$ at $p$ by:
$$\matrix
X \hfill&= \nabla\opSL'_r(\Phi,\|\nabla\Phi_0\|,E),\hfill\cr
Y \hfill&= \nabla\opSL'_r(\Phi + x(d_{H'} - f),\|\nabla\Phi_0\|,E).\hfill\cr
\endmatrix$$
\noindent Denote $P=x(d_{H'}-f)$. At $p$:
$$
\opHess(P) = \nabla x\otimes \nabla d_{H'} + \nabla d_{H'}\otimes\nabla x.
$$
\noindent At $p$, for all $i$, by definition, $\langle e_i,\nabla d_{H'}\rangle=0$. Likewise $\langle \nabla\Phi,\nabla d_{H'}\rangle=0$. Thus, recalling the formula for $\nabla e_i$:
$$\matrix
X\opHess(P)(e_i,e_j) \hfill&= (\nabla_X\opHess(P))(e_i,e_j) + \opHess(P)(\nabla_X e_i,e_j) + \opHess(P)(e_i,\nabla_X e_j)\hfill\cr
&= (\nabla_X\opHess(P))(e_i,e_j) - \opHess(d_{H'})(X,e_i)x_{;j} - \opHess(d_{H'})(X,e_j)x_{;i}.\hfill\cr
\endmatrix$$
\noindent Extend $(e_i)_{1\leqslant i\leqslant n-1}$ to a basis for $T_pM$ by defining:
$$
e_0 = \msf{N}_p,\qquad e_n = \nabla d_{H'}.
$$
\noindent Then, with respect to this basis:
$$
(Y-X)_k = -(A-B)x_{;k} - 2\sum_{i=i}^{n-1}\mu_i x_{;i} f_{;ik} + N_{ik}x_{;k},
$$
\noindent where $N=O(\delta)$. Consider the linear map, $M$, given by:
$$
(MV)_k = (A-B)V_k + 2\sum_{i=1}^{n-1}\mu_if_{;ik}V_i.
$$
\noindent Suppose that $MV=0$, then:
$$\matrix
&\sum_{k=1}^{n-1}(MV)_k\mu_k V_k \hfill&= 0\hfill\cr
\Rightarrow\hfill&\sum_{k=1}^{n-1}(A-B)\mu_k V_k^2 + 2\sum_{i,j=1}^{n-1}(\mu_i V_i)(\mu_j V_j)f_{;ik}\hfill&= 0.\hfill\cr
\endmatrix$$
\noindent Since $(A-B)>0$ and $f_{;ij}$ is positive definite, it follows that:
$$
V_k = 0\text{ for all }1\leqslant k\leqslant n-1.
$$
\noindent This in turn yields:
$$
(A-B) V_0 = (A-B) V_n = 0.
$$
\noindent And so $V=0$. It follows that $M$ is invertible. There therefore exists $x$ such that, at $p$:
$$
\nabla\opSL'_r(\Phi_0 + x(d_{H'} - f),\|\nabla\Phi_0\|) = 0.
$$
\noindent The result follows.\qed
\medskip
\noindent For $M>0$, we define $\Phi$ by:
$$
\Phi = \Phi_0 + x(d_{H'} - f) + Md_{H'}^2.
$$
\proclaim{Proposition \nextprocno}
\noindent If $D$ represents the Grassmannian distance between two $(n-2)$-dimensional subspaces then:
$$
D(E(\Phi_0,d_{H'}),E(\Phi,d_{H'})) = O(d_p^2) + O(d_{H'}).
$$
\endproclaim
\proof Since $xf$ is of order $3$ at $p$:
$$
\nabla\Phi = \nabla\Phi_0 + (x + 2Md_{H'})\nabla d_{H'} + O(d_p^2) + O(d_{H'}).
$$
\noindent Thus:
$$
\langle \nabla\Phi,\nabla d_{H'}\rangle = \langle \nabla\Phi_0 + O(d_p^2) + O(d_{H'}),\nabla d_{H'}\rangle,
$$
\noindent where $\langle\cdot,\cdot\rangle$ here represents the subspace generated by two vectors. The result follows.\qed
\medskip
\proclaim{Corollary \nextprocno}
\noindent If $D$ represents the Grassmannian distance between two $(n-2)$-dimensional subspaces then:
$$
D(E,E(\Phi,d_{H'})) = O(d_p^2) + O(d_{H'}).
$$
\endproclaim
\proclabel{CorollarySubspacesCloseII}
\proof This follows by the triangle inequality and Proposition \procref{PropositionSubspacesCloseI}.\qed
\proclaim{Proposition \nextprocno}
\noindent Suppose that $M\epsilon^2<1$, $0\leqslant d_{H'}<\epsilon^2$ and $0\leqslant d_p<\epsilon$. Then:
$$
\|\nabla\Phi\|^2 \geqslant \|\nabla\Phi_0\|^2 + O(\epsilon^2).
$$
\endproclaim
\proclabel{PropositionControlledDerivative}
\proof We examine each of the terms seperately. Trivially:
$$
\|\nabla x(d_{H'} - f)\|^2 \geqslant 0.
$$
\noindent And:
$$
\|\nabla Md_{H'}^2\|^2 = 4M^2d_{H'}^2.
$$
\noindent We now consider the interaction terms. Recalling that $0\leqslant d_{H'}<\epsilon^2$ and $0\leqslant d_p<\epsilon$:
$$
\langle\nabla\Phi_0,\nabla x(d_{H'} - f)\rangle = x\langle\nabla\Phi_0,\nabla d_{H'}\rangle + O(\epsilon^2).
$$
\noindent Recalling that $\langle\nabla\Phi_0,\nabla d_{H'}\rangle=O(\epsilon)$, we obtain:
$$
\langle\nabla\Phi_0,\nabla x(d_{H'} - f)\rangle = O(\epsilon^2).
$$
\noindent Likewise:
$$
\langle\nabla\Phi_0,\nabla Md_{H'}^2\rangle = Md_HO(\epsilon).
$$
\noindent Finally, since $M\epsilon^2<1$:
$$
\langle \nabla x(d_{H'} - f), \nabla Md_{H'}^2\rangle = 2Md_{H'}x + O(\epsilon^2).
$$
\noindent Combining these terms yields:
$$
\|\nabla\Phi\|^2 \geqslant \|\Phi_0\|^2 + Md_{H'}(Md_{H'} - O(\epsilon)) + O(\epsilon^2).
$$
\noindent However:
$$
Md_{H'}(Md_{H'} - O(\epsilon)) \geqslant -O(\epsilon^2).
$$
\noindent the result now follows.\qed
\proclaim{Proposition \nextprocno}
\noindent Let $\epsilon>0$. If $M\epsilon^2<1$, $d_{H'}<\epsilon^2$ and $d_p<\epsilon$, then:
$$
\opSL'_r(\Phi,\|\nabla\Phi\|,E(\nabla\Phi,\nabla d_{H'})) \leqslant \theta - \pi/2  + O(\epsilon^2).
$$
\endproclaim
\proclabel{PropositionSecondOrderVariationII}
\proof Define $\Phi_1$ by:
$$
\Phi_1 = \Phi_0 + x(d_{H'} - f).
$$
\noindent By Proposition \procref{PropositionSecondOrderVariation}:
$$
\opSL'_r(\Phi_1,\|\nabla\Phi_0\|,E) \leqslant \theta - \pi/2  + O(\epsilon^2).
$$
\noindent Since $\opHess(\Phi(1))=O(1)$, by Proposition \procref{PropositionControlledDerivative} and Corollary \procref{CorollarySubspacesCloseII}:
$$
\opSL'_r(\Phi_1,\|\nabla\Phi\|,E(\nabla\Phi,\nabla d_{H'})) \leqslant \theta - \pi/2  + O(\epsilon^2).
$$
\noindent Differentiating $Md_{H'}^2$ yields:
$$
\opHess(Md_{H'}^2) = 2M\nabla d_{H'}\otimes \nabla d_{H'} + 2M d_{H'}\opHess(d_{H'}).
$$
\noindent The first term vanishes along $(\nabla d_{H'})^\perp$, and the second term is negative, and thus does not affect the inequality either. The result now follows.\qed
\goodbreak
\newsubhead{Boundary Lower Bounds for the Normal}
\noindent Let $M$, $H$, $\hat{\Sigma}$ and $\Omega$ be as in the preceeding section. Let $\Sigma$ be a $C^0$ convex hypersurface such that:
\medskip
\myitem{(i)} $\Sigma$ lies between $\Omega$ and $\hat{\Sigma}$;
\medskip
\myitem{(ii)} the interior of $\Sigma$ is smooth; and
\medskip
\myitem{(iii)} $\partial\Sigma=\partial\Omega=\Gamma$.
\medskip
\noindent Let $\msf{N}_\Sigma$ be the exterior normal to $\Sigma$ (which is continuous).
\proclaim{Lemma \nextprocno}
\noindent There exists $K>0$ such that if $R_\theta(\Sigma)=r$, then $\mu(p,\msf{N}_\Sigma(p),r,\theta)<K$ for all $p\in\partial\Sigma=\Gamma$.
\endproclaim
\proclabel{LemmaBoundaryLowerNormalBounds}
\proof We assume the contrary and obtain a contradiction. By continuity, there exists $p\in\Gamma$ such that:
$$
\mu(p,\msf{N}_\Sigma(p),r,\theta) = +\infty.
$$
\noindent Define $d_H$, $d_\Omega$, $\Phi_0$ and $K$ as in the previous section. Since $\Sigma$ is convex and $\partial\Sigma=\Gamma$ is smooth, $\msf{N}_\Sigma$ is continuous at $p$. Since, by definition $\msf{N}_\Sigma(p)=\nabla(d_H - Kd_\Omega)(p)$, (c.f. Lemma \procref{LemmaUniformModulusOfContinuity}) there exists a continuous function $\delta:[0,\infty[\rightarrow[0,\infty[$ such that $\delta(0)=0$ and, along $\Sigma$:
$$
\|\pi(\nabla \Phi_0)(q)\| \leqslant \delta(d_p(q)),
$$
\noindent where $\pi$ is the orthogonal projection onto $T\Sigma$. Define $H'$, $d_{H'}$ as in the previous section. For $\epsilon>0$ small, define $U_\epsilon$ by:
$$
U_\epsilon = \left\{ q\in M\text{ s.t. }d_p(q)<\epsilon,d_{H'}(q)<\epsilon^2\right\}.
$$
\noindent Along $\partial\Sigma=\Gamma$, $\Phi_0=0$. Recall that any convex set is $C^{0,1}$. Thus, along $\partial U_\epsilon\minter\Sigma$, $d_\Sigma(q,\partial\Sigma)=O(d_{H'})=O(\epsilon^2)$, where $d_\Sigma$ is the Riemannian distance inside $\Sigma$, and so:
$$
\Phi(q) = \delta(\epsilon)O(\epsilon^2),
$$
\noindent along $\partial U_\epsilon\minter\Sigma$. Since $\Gamma$ is strictly convex and lies strictly inside $H'$, there exists a function $f$ such that:
\medskip
\myitem{(i)} $f(p),\nabla f(p)=0$ and the restriction of $\opHess(f)(p)$ to $H'$ is positive definite; and
\medskip
\myitem{(i)} $d_{H'} - f = O(d_p^3)$ along $\Gamma$.
\medskip
\noindent We define $\Phi$ as in the previous section. Along $\partial\Sigma\minter U_\epsilon=\Gamma\minter U_\epsilon$:
$$
\Phi(q) \geqslant Md_{H'}^2 - O(\epsilon^4).
$$
\noindent This is positive for sufficiently large $M$. Likewise, along $\partial U_\epsilon\minter\Sigma$:
$$
\Phi(q) \geqslant Md_{H'}^2 - \delta(\epsilon)O(\epsilon^2).
$$
\noindent There thus exists $K_1>0$ independant of $\epsilon$ such that, if $M=K_1\delta(\epsilon)\epsilon^2$, then $\Phi\geqslant 0$ along $\partial
U_\epsilon\minter\Sigma$.
\medskip
\noindent Let $A$ be the restriction of $\|\nabla\Phi\|^{-1}\opHess(\Phi)$ to $\nabla\Phi^\perp$. Let $\lambda_1\leqslant...\leqslant\lambda_n$ be the eigenvalues of A. Let $\lambda_1'\leqslant...\leqslant\lambda_n'$ be the eigenvalues of the restriction of $A$ to $\nabla\Phi^\perp\minter\nabla d_{H'}^\perp$. By the Minimax Principal, for $1\leqslant i\leqslant (n-1)$:
$$
\lambda_i<\lambda_i'.
$$
\noindent Thus, by Proposition \procref{PropositionSecondOrderVariationII}, there exists $K_2>0$, also independant of $\epsilon$ such that:
$$
\sum_{i=1}^{n-1}\opArcTan(r\lambda_i)\leqslant \theta - \pi/2 + K_2\epsilon^2.
$$
\noindent However:
$$
\lambda_n = O(M).
$$
\noindent There thus exists $K_3>0$, independant of $\epsilon$, such that:
$$\matrix
\opSL_r(A) \hfill&\leqslant \theta + (K_2\epsilon^2 - K_3M^{-1})\hfill\cr
&= \theta + \epsilon^2(K_2 - K_1K_3\delta(\epsilon)^{-1}).\hfill\cr
\endmatrix$$
\noindent Since $\delta(\epsilon)$ tends to $0$ as $\epsilon$ tends to $0$, there exists $\eta>0$ such that, for $\epsilon$ sufficiently small, throughout $U_\epsilon$:
$$
SL_r(A) \leqslant \theta - \eta < \theta.
$$
\noindent It follows that if $\Sigma_t =\Phi^{-1}(\left\{t\right\})$ for all $t$, then:
$$
\opSL_r(\Sigma_t\minter U_\epsilon) \leqslant \theta - \eta < \theta.
$$
\noindent At $p$, $\nabla^\Sigma\Phi=0$. Thus, reducing $\epsilon$ further if necessary, we may deform $\Phi$ slightly to $\Phi'$ (by subtracting a very small multiple of $d_\Omega$, for example) such that $\Phi'$ is non-negative along $\partial(\Sigma\minter U_\epsilon)$, $\Phi$ is strictly negative over a non trivial subset of $\Sigma\minter U_\epsilon$, and, if $\Sigma_t'=(\Phi')^{-1}(\left\{t\right\})$ for all $t$, then:
$$
SL_r(\Sigma_t'\minter U_\epsilon)\leqslant \theta - \eta/2 < \theta.
$$
\noindent Let $p\in\Sigma$ be the point where $\Phi$ is minimised. Let $t_0=\Phi'(p)$. Since $p$ lies in the interior of $\Sigma$, $\Sigma_{t_0}'$ is smooth at this point. $\Sigma_{t_0}'$ is an interior tangent to $\Sigma$ at $p$. In particular, $\Sigma_{t_0}'$ is convex near $p$ and:
$$
R_\theta(\Sigma_t'\minter U_\epsilon) < r,
$$
\noindent which is absurd, by the Geometric Maximum Principle (Lemma \procref{LemmaGeomMaxPrinc}), and the result follows.\qed
\medskip
\noindent Using limits yields:
\proclaim{Corollary \nextprocno}
\noindent There exists $K>0$ such that if $\Sigma$ is a smooth hypersurface lying between $\Omega$ and $\hat{\Sigma}$ such that $\partial\Sigma=\partial\Omega=\Gamma$ and $R_\theta(\Sigma)=r\in]R_0,R_1[$ is constant, then:
$$
\mu(p,\msf{N}_\Sigma(p),r,\theta) < K\text{ for all }p\in\Gamma.
$$
\endproclaim
\proclabel{CorollaryLowerNormalBound}
\goodbreak
\newsubhead{Constructing Barriers II}
\noindent Let $M$ be an $(n+1)$-dimensional manifold. Let $\delta_0\geqslant 0$ be small. Let $H\subseteq M$ be a smooth convex hypersurface such that:
$$
R_\theta(H) = R_0,
$$
\noindent where $R_0$ is constant. Let $\Omega\subseteq H$ be a bounded open subset of $H$. Let $\hat{\Sigma}\subseteq M$ be a convex hypersurface such that $\partial\hat{\Sigma}=\partial\Omega=:\Gamma$ and such that $R_\theta(\hat{\Sigma})>R_1$ in the weak sense. Suppose that $\Gamma$ is strictly convex as a subset of $M$ with respect to the outward pointing normal to $\Gamma$ in $\hat{\Sigma}$.
\medskip
\noindent Let $\Sigma$ be a smooth immersed convex hypersurface lying between $\Omega$ and $\hat{\Sigma}$ such that $\partial\Sigma=\Gamma$ and:
$$
R_\theta(\Sigma) = r \in [R_0,R_1].
$$
\noindent Let $\msf{N}$, $II$ and $A$ be the unit normal, the second fundamental form and the shape operator respectively of $\Sigma$. Define $B$ over $\Sigma$ such that:
$$
B^{ij}(\delta_{jk} + r^{-2} A^2_{jk}) = {\delta^i}_j.
$$
\noindent Define the operator $\Delta^B:C^\infty(\Sigma)\rightarrow C^\infty(\Sigma)$ by:
$$
\Delta^B f = B^{ij}f_{;ij},
$$
\noindent where $f_{;ij}$ is the Hessian of $f$ with respect to Levi-Civita covariant derivative of $\Sigma$. We aim to construct barriers for $\Delta^B$ at any point of $\partial\Sigma$.
\medskip
\noindent Let $p\in\partial\Sigma$ be a point. There exists a strictly convex hypersurface $\Sigma'$ tangent to $\partial\Sigma=\Gamma$ at $p$ such that $\Sigma$ lies in its interior and the normal to $\Sigma$ at $p$ also points into its interior.
\proclaim{Lemma \nextprocno}
\noindent Let $d_{\Sigma'}$ denote the distance to $\Sigma'$. Let $U$ be a neighbourhood of $p$ such that, throughout $\Sigma\minter U$:
$$
\langle \nabla d_{\Sigma'},\msf{N}\rangle \geqslant 0.
$$
\noindent Then, throughout $\Sigma\minter U$:
\medskip
\myitem{(i)} $d_{\Sigma'}\geqslant 0$, and
\medskip
\myitem{(ii)} $\Delta^B d_{\Sigma'}\leqslant 0$.
\endproclaim
\proclabel{LemmaConcaveFunction}
\remark $U$ depends only on the modulus of continuity for $\msf{N}$ near $p$.
\medskip
\proof See Corollary \procref{CorollaryConcaveFunction}.\qed
\proclaim{Lemma \nextprocno}
\noindent Let $\theta\in](n-1)\pi/2,n\pi/2[$ be an angle. Let $r\in]0,\infty[$. Suppose that $R_\theta(\Sigma)=r$. There exists $\epsilon,\delta>0$ which only depend on $\theta$ and $r$ such that, throughout $B_\delta(p)\minter\Sigma$:
$$
\Delta^B d_p^2\geqslant\epsilon.
$$
\endproclaim
\proclabel{LemmaLowerBoundOfLaplacian}
\proof By Lemma \procref{LemmaConcavityPreserved} and Corollary \procref{CorollaryApproximateHessianOfDistance}:
$$
(d_p^2)_{;ij} = 2\delta_{ij} - 2A_{ij}d_p \langle\nabla d_p,\msf{N}\rangle + O(d_p^2).
$$
\noindent By Lemma \procref{LemmaSLCIsLowestEval}, there exists $K_1>0$, which only depends on $\theta$ and $r$ such that:
$$
Tr(B)\geqslant \frac{1}{K_1}.
$$
\noindent Thus, throughout $\Sigma$:
$$
\Delta^B(d_p^2) \geqslant \frac{2}{K_1} - O(d_p).
$$
\noindent There exists $\delta$, which only depends on $K$ such that, for $d_p<\delta$, the error term is less than $\frac{1}{K}$ in magnitude. The result now follows.\qed
\newsubhead{Second Order Boundary Estimates}
\noindent Let $f$ be the signed distance to $\Sigma$. $f$ is a real valued function which is smooth in a nieghbourhood of $\Sigma$. By definition:
$$
\|\nabla f\|=1.
$$
\noindent For $X$, $Y$ tangent to $\Sigma$:
$$
\opHess(f)(X,Y) = II(X,Y),
$$
\noindent where $II$ is the second fundamental form of $\Sigma$. Let $p\in\partial\Sigma$. Let $X$ be a vector field over $\Bbb{H}^{n+1}$ which is tangent to $\partial\Sigma$ (but not necessarily tangent to $\Sigma$). Define $\varphi=Xf$. For any $Y$ tangent to $\Sigma$:
$$
Y\varphi = \opHess(f)(X,Y) + \langle\nabla f,\nabla_Y X\rangle =  II(X,Y) + \langle\nabla f,\nabla_Y X\rangle.
$$
\noindent Thus, a-priori bounds on $X$ and $\varphi$ yield a-proiri bounds on $II$.
\proclaim{Lemma \nextprocno}
\noindent For $X,Y\in T\Sigma$:
$$
(\nabla_\msf{N}\opHess(f))(X,Y) = \langle R_{\msf{N}X}\msf{N},Y\rangle - \langle A^2X,Y\rangle.
$$
\endproclaim
\proclabel{LemmaNormalDerivative}
\proof Define $\Phi:\Sigma\times]-\epsilon,\epsilon[\rightarrow M$ by:
$$
\Phi(p,t) = \opExp_p(t\msf{N}(p)).
$$
\noindent Pulling back through $\Phi$, we identify $M$ with $\Sigma\times]-\epsilon,\epsilon[$ and $\msf{N}$ with $\partial_t$. In particular, if $X$ is tangent to $\Sigma$, then $[X,\msf{N}]=0$. Trivially:
$$
\nabla f = \msf{N}.
$$
\noindent Thus:
$$
\opHess(f)(X,Y) = \langle\nabla_Y\msf{N},X\rangle = A(X,Y).
$$
\noindent Bearing in mind that $\nabla_\msf{N}\msf{N}=0$:
$$\matrix
(\nabla_N\opHess(f))(X,Y) \hfill&= \msf{N}\langle\nabla_X\msf{N},Y\rangle
- A(\nabla_\msf{N}X,Y) - A(X,\nabla_\msf{N}Y)\hfill\cr
&= \langle\nabla_\msf{N}\nabla_X\msf{N},Y\rangle + \langle\nabla_X\msf{N},\nabla_\msf{N}Y\rangle
- A(\nabla_\msf{N}X,Y) - A(X,\nabla_\msf{N}Y)\hfill\cr
&= \langle R_{\msf{N}X}\msf{N},Y\rangle + \langle \nabla_X\msf{N},\nabla_Y\msf{N}\rangle
- A(\nabla_X\msf{N},Y) - A(X,\nabla_Y\msf{N})\hfill\cr
&= \langle R_{\msf{N}X}\msf{N},Y\rangle - \langle A^2X,Y\rangle.\hfill\cr
\endmatrix$$
\noindent The result follows.\qed
\proclaim{Lemma \nextprocno}
\noindent There exists $K>0$, which only depends on $X$, $r$, $\theta$ and the structure of $M$ such that, throughout $\Sigma$:
$$
\left|\Delta^B\varphi\right|\leqslant K.
$$
\endproclaim
\proclabel{LemmaBoundedLaplacian}
\proof Let Latin indices represent directions in $TM$ and let Greek indices represent directions in $T\Sigma$. Let $\nu$ represent the exterior normal direction to $\Sigma$.
\medskip
\noindent Let $;$ denote covariant differentiation with respect to the Levi-Civita covariant derivative of $M$. Let $O(1)$ represent terms bounded in terms of $X$, $r$, $\theta$ or the structure of $M$. Recall that $A_{\alpha\beta}=f_{;\alpha\beta}$. This is symmetric in $\alpha$ and $\beta$. By definition of curvature:
$$
f_{;\alpha\beta k} = f_{;\alpha k\beta} + {R_{\beta k\alpha}}^lf_{;l} = f_{;k\alpha\beta} + O(1).
$$
\noindent Since $f_{;\nu k}=0$ for all $k$, for all $X,Y,Z$ tangent to $\Sigma$:
$$
(\nabla\opHess(f))(Y,Z;X) = (\nabla^\Sigma\opHess(f))(Y,Z;X).
$$
\noindent Thus, differentiating $R_\theta(A)=r$ along $\Sigma$ yields:
$$\matrix
&B^{\alpha\beta}f_{;\alpha\beta\gamma} \hfill&= 0\hfill\cr
\Rightarrow\hfill&B^{\alpha\beta}f_{;\gamma\alpha\beta}\hfill&= O(1).\hfill\cr
\endmatrix$$
\noindent We remark in passing that it is at this stage that the differential condition on $f$ (and therefore $\Sigma$) is used. We now consider the normal derivative. By Lemma \procref{LemmaNormalDerivative}:
$$\matrix
&f_{;\alpha\beta\nu} \hfill&= R_{\nu\alpha\nu\beta} -(A^2)_{\alpha\beta}\hfill\cr
\Rightarrow\hfill&f_{;\nu\alpha\beta} \hfill&= -(A^2)_{\alpha\beta}.\hfill\cr
\endmatrix$$
\noindent Thus:
$$
\left|B^{\alpha\beta}f_{;\nu\alpha\beta}\right| \leqslant nr^2.
$$
\noindent Thus, for all $k$:
$$
\left|B^{\alpha\beta}f_{;k \alpha\beta}\right| = O(1).
$$
\noindent We now consider $\varphi = X^k f_{;k}$:
$$
B^{\alpha\beta}\varphi_{;\alpha\beta} = B^{\alpha\beta}{X^k}_{;\alpha\beta}f_{;k}
+2B^{\alpha\beta}{X^k}_{;\alpha}f_{;k\beta}
+B^{\alpha\beta}X^k f_{;k\alpha\beta}.
$$
\noindent Since $\|B\|\leqslant 1$, the first term is controlled by a-priori bounds on $\nabla^2X$. The third term is controlled by a-priori bounds on $X$ and the preceeding discussion. We now control the second term. Recalling that $f_{;\nu k}=0$ for all $k$:
$$
B^{\alpha\beta}{X^k}_{;\alpha}f_{;k\beta} = B^{\alpha\beta}A_{\gamma\beta}{X^{\gamma}}_{;\alpha}.
$$
\noindent Since $\|BA\|\leqslant 1$, this term is controlled by a-priori bounds on $\nabla X$. Finally, since $\nabla f$ is the unit normal to $\Sigma$, by Lemma \procref{LemmaConcavityPreserved}:
$$
\opHess^\Sigma(\varphi)_{\alpha\beta} = \varphi_{;\alpha\beta} - A_{\alpha\beta}.
$$
\noindent Thus:
$$
\Delta^B\varphi = B^{\alpha\beta}(\varphi)_{;\alpha\beta} - B^{\alpha\beta}A_{\alpha\beta}.
$$
\noindent Since $\|AB\|\leqslant 1$, the result follows.\qed
\proclaim{Lemma \nextprocno}
\noindent There exists $K$, which only depends on $M$, $H$, $\hat{\Sigma}$, $r$, $\theta$ and the modulus of continuity of $\Sigma$ near $\partial\Sigma$ such that, along $\partial\Sigma$:
$$
\|A\|\leqslant K.
$$
\endproclaim
\proclabel{LemmaSecondOrderBoundaryEstimates}
\proof Let $p\in\partial\Omega$. The normal to $\Sigma$ at $p$ lies between the normals to $H$ and $\partial\hat{\Sigma}$ at $p$. There thus exists a convex hypersurface, $H'$, which is an exterior tangent to $\partial\Omega$ at $p$ and such that the normal to $\Sigma$ at $p$ points into $H'$. Define $d_{H'}$ by:
$$
d_{H'}(q) = d(q,H').
$$
\noindent By Lemma \procref{LemmaConcaveFunction}, there exists a neighbourhood $U_1$ of $p$, which only depends on $\hat{\Sigma}$ and the modulus of continuity of $\Sigma$ at the boundary, such that, throughout $\Sigma\minter U_1$:
$$
\Delta^Bd_{H'}\leqslant 0.
$$
\noindent Define $d_p$ by:
$$
d_p(q) = d(q,p).
$$
\noindent By Lemma \procref{LemmaLowerBoundOfLaplacian}, there exists $\epsilon>0$ and a neighbourhood $U_2$ of $p$ such that, throughout $\Sigma\minter U_2$:
$$
\Delta^Bd_p^2\geqslant\epsilon.
$$
\noindent Let $f$ be the perpindicular distance to $\Sigma$. Let $X$ be a vector field tangent to $\partial\Omega$. Consider the function $\varphi = Xf$. $\varphi$ vanishes along $\partial\Omega$. Since $\|\nabla f\|=1$, there exists $K_1>0$, which only depends on $X$ such that, throughout $\Sigma$:
$$
\left|\varphi\right| = \|Xf\| \leqslant K_1.
$$
\noindent By Lemma \procref{LemmaBoundedLaplacian}, there exists $K_2>0$ such that, throughout $\Sigma$:
$$
\left|\Delta^B\varphi\right| \leqslant K_2.
$$
\noindent Choosing $\delta>0$ such that $B_\delta(p)\subseteq U_1\minter U_2$, there exists $A_->0$ such that, throughout $B_\delta(p)\minter\Sigma$:
$$
\Delta^B(\varphi - Ad_p^2)\leqslant 0.
$$
\noindent There exists $B_->0$ such that:
\medskip
\myitem{(i)} $\Delta^B(\varphi + B_-d_{H'} - A_-d_p^2)\leqslant 0$ throughout $B_\delta(p)\minter\Sigma$; and
\medskip
\myitem{(ii)} $(\varphi + B_-d_{H'} - A_-d_p^2)\geqslant 0$ along $\partial(B_\delta(p)\minter\Sigma)$.
\medskip
\noindent Thus, by the maximum principal, throughout $B_\delta(p)\minter\Sigma$:
$$
\varphi \geqslant B_-d_{H'} - A_-d_p^2.
$$
\noindent Likewise, reducing $\delta$ if necessary, there exists $B_+$ and $A_+$ such that, throughout $B_\delta(p)\minter\Sigma$:
$$
\varphi \leqslant B_+d_{\Sigma'} - A_+d_p^2.
$$
\noindent We thus obtain a-priori bounds on $\nabla\varphi$ at $p$. Since $X$ is arbitrary, this yields a-priori bounds on $\opHess(f)(X,Y)$ for all pairs of vectors $X,Y\in T_p\Sigma$ where at least one of $X$ or $Y$ is tangent to $\partial\Sigma$. Since the second fundamental form of $\Sigma$ is the restriction to $T\Sigma$ of the hessian of $f$ (and since $\|\nabla f\|=1$), we obtain a-priori bounds on $A(X,Y)$ for all such pairs of vectors.
\medskip
\noindent Let $(e_1,...,e_{n-1})$ be an orthonormal basis of $T_p\partial\Sigma$ which diagonalises the restriction of $A$. Let $e_n$ be the inward pointing normal of $\partial\Sigma$ at $p$. With respect to this basis, there exists $0<\lambda_1<...<\lambda_{n-1}$ and $M>0$ such that:
$$
A = \pmatrix D\hfill& O(1)\hfill\cr O(1)&M\hfill\cr\endpmatrix,
$$
\noindent where $D=\opDiag(\lambda_1,...,\lambda_{n-1})$ is the diagonal matrix with entries $\lambda_1,...,\lambda_{n-1}$. Let $\lambda'_1,...,\lambda'_n$ by the eigenvalues of $A$. By Lemma $1.2$ of \cite{CaffNirSprIII}:
$$
\lambda'_i =\left\{ \matrix \lambda_i + o(1)\text{ if }1\leqslant i\leqslant n-1\hfill\cr
M(1 + o(M^{-1}))\text{ if }i=n.\hfill\cr\endmatrix\right.
$$
\noindent However, by Corollary \procref{CorollaryLowerNormalBound}, there exists $K$ such that $\mu(p,\msf{N}_p,r,\theta)<K$. $M$ therefore cannot become arbitrarily large. We thus obtain a-priori bounds on $M$, and the result now follows.\qed
\newsubhead{Second Order Interior Estimates}
\noindent Let $M$ be a Riemannian manifold. Let $K\subseteq M$ be a compact subset. Let $\nabla$ be the Levi-Civita covariant derivative over $M$ and let $R$ be the Riemann curvature tensor of $M$. Choose $\theta\in](n-1)\pi/2,n\pi/2[$ and $r>0$. Let $K\subseteq M$ be a compact subset. Let $\Sigma\subseteq M$ be a smooth, convex, immersed hypersurface contained in $K$ such that:
$$
R_\theta(\Sigma) = r.
$$
\noindent Let $A$ and $\msf{N}$ be the second fundamental form and the exterior normal of $\Sigma$ respectively.
\proclaim{Lemma \nextprocno}
\noindent Let $H$ be the mean curvature of $\Sigma$. There exists $C>0$ which only depends on $r$ and the norms of $R$ and $\nabla R$ over $K$ such that:
$$
\Delta^B H \geqslant -C(1+H) + \sum_{i,j=1}^n\frac{(\lambda_i-\lambda_j)r^2\lambda_i\lambda_j}{(1+\lambda_j^2)}.
$$
\endproclaim
\proclabel{LemmaLaplacianOfMeanCurvature}
\proof Choose $p\in\Sigma$ and let $e_1,...,e_n$ be an orthonormal basis of eigenvectors for $A$ in $T_p\Sigma$. Let $\lambda_1,...,\lambda_n$ be the corresponding eigenvalues of $r^{-1}A$. Let $;$ denote covariant differentiation with respect to the Levi-Civita covariant derivative of $\Sigma$. Let the index $\nu$ denote the direction normal to $\Sigma$. Differentiating the curvature condition twice yields:
$$\matrix
\sum_{i=1}^n\frac{1}{(1+\lambda_i^2)}A_{iik}\hfill&=0,\hfill\cr
\sum_{i=1}^n\frac{1}{(1+\lambda_i^2)}A_{iipq}\hfill&=\sum_{i,j=1}^n\frac{r^{-2}(\lambda_i+\lambda_j)}{(1+\lambda_i^2)(1+\lambda_j^2)}A_{ij;p}A_{ij;q}.\hfill\cr
\endmatrix$$
\noindent Let $R^\Sigma$ be the Riemann curvature tensor of $\Sigma$. By definition of curvature:
$$\matrix
A_{ij;k} - A_{ik;j} \hfill&= R_{kj\nu i},\hfill\cr
A_{ij;kl} - A_{ij;lk} \hfill&= R^\Sigma_{klij}\lambda_j + R^\Sigma_{klji}\lambda_i.\hfill\cr
\endmatrix$$
\noindent Thus, for all $i$ and $j$:
$$\matrix
A_{ii;jj} \hfill&= A_{ij;ij} + R_{ji\nu i;j}\hfill\cr
&= A_{ij;ji} + R_{ji\nu i;j} + R^\Sigma_{ijij}\lambda_j + R^\Sigma_{ijji}\lambda_i\hfill\cr
&= A_{jj;ii} + R_{ji\nu j;i} + R_{ji\nu i;j} + R^\Sigma_{ijij}\lambda_j + R^\Sigma_{ijji}\lambda_i.\hfill\cr
\endmatrix$$
\noindent Applying the second derivative of the curvature condition:
$$
\sum_{i,j=1}^n\frac{1}{(1+\lambda_j^2)}A_{jj;ii} \geqslant 0.
$$
\noindent For any $1$-form, $\xi$:
$$\matrix
\xi_{i;j} \hfill&= (\nabla\xi)_{ij} - A_{ij}\xi_\nu,\hfill\cr
\xi_{\nu;j} \hfill&= (\nabla_{\partial_j}\xi)(\msf{N}) + {A_j}^k\xi_k.\hfill\cr
\endmatrix$$
\noindent Thus:
$$\matrix
R_{ji\nu j;i} \hfill&= (\nabla R)_{ji\nu ji} - r\lambda_i(1-\delta_{ij})R_{j\nu\nu j} + r\lambda_i R_{jiij},\hfill\cr
R_{ji\nu i;j} \hfill&= (\nabla R)_{ji\nu ij} - r\lambda_j(1-\delta_{ij})R_{\nu i\nu i} + r\lambda_j R_{jiji}.\hfill\cr
\endmatrix$$
\noindent This yields:
$$\matrix
\sum_{i,j=1}^n\frac{1}{(1+\lambda_j^2)}(R_{ji\nu j;i} + R_{ji\nu i;j})\hfill&
\geqslant \sum_{i,j=1}^n\frac{-r\lambda_i(1-\delta_{ij})}{(1+\lambda_j^2)}R_{j\nu\nu j}
+ \sum_{i,j=1}^n\frac{-r\lambda_j(1-\delta_{ij})}{(1+\lambda_j^2)}R_{\nu i\nu i}\hfill\cr
&+\sum_{i,j=1}^n\frac{r(\lambda_i-\lambda_j)}{(1+\lambda_j^2)}R_{ijji} - C_1,\hfill\cr
\endmatrix$$
\noindent where $C_1$ only depends on the norm of $\nabla R$ over $K$. The first and third terms on the right hand side is bounded by a multiple of $H$ times the norm of $R$ over $K$. Likewise, the second term is bounded in terms of the norm of $R$ over $K$. Thus:
$$
\sum_{i,j=1}^n\frac{1}{(1+\lambda_j^2)}(R_{ji\nu j;i} + R_{ji\nu i;j}) \geqslant  - C(1+H),
$$
\noindent where $C$ only depends on $r$ and the norms of $R$ and $\nabla R$ over $K$. Finally, since $A$ is the shape operator of $\Sigma$:
$$
R^\Sigma_{ijij}\lambda_j + R^\Sigma_{ijji}\lambda_i = (\lambda_i - \lambda_j)R_{ijji} + r^2(\lambda_i-\lambda_j)\lambda_i\lambda_j.
$$
\noindent The result follows.\qed
\proclaim{Lemma \nextprocno}
\noindent There exists $D>0$, which only depends on $r$,$\theta$ and the norms of $R$ and $\nabla R$ over $K$ such that:
$$
H\geqslant D\Rightarrow \Delta^B H\geqslant 0.
$$
\endproclaim
\proclabel{LemmaMeanCurvatureIsSubharmonic}
\proof Symmetrising the inequality obtained in Lemma \procref{LemmaLaplacianOfMeanCurvature} yields:
$$
\Delta^B H \geqslant -C(1+H) + \sum_{i,j=1}^n F(\lambda_i,\lambda_j;r,C),
$$
\noindent where $F$ is given by:
$$
F(x,y;r,C) = \frac{r^2xy}{2(1+x^2)(1+y^2)}(x^3 + y^3 - x^2y - y^2x).
$$
\noindent Since $R_\theta(A)=r$, there exists $\epsilon>0$ which only depends on $r$ and $\theta$ such that
$\lambda_i\geqslant\epsilon$ for all $i$. We observe that, for all $x,t\geqslant 0$:
$$
F(x,y;r,C) \geqslant 0.
$$
\noindent Without loss of generality:
$$
\lambda_1 \geqslant H/n,\qquad \epsilon\leqslant \lambda_n \leqslant r\opTan(\theta/n).
$$
\noindent Consequently, $F(\lambda_1,\lambda_n;r,C)$ grows like $H^2$ as $H\rightarrow+\infty$. In particular, there exists $D>0$ such that, for $H\geqslant D$:
$$
F(\lambda_1,\lambda_n;r,C)\geqslant C(1+H).
$$
\noindent This is the desired value for $D$ and the result follows.\qed
\proclaim{Proposition \nextprocno}
\noindent There exists $K$, which only depends on $M$, $H$, $\hat{\Sigma}$, $r$, $\theta$ and the modulus of continuity of $\Sigma$ near $\partial\Sigma$ such that, throughout $\Sigma$:
$$
\|A\|\leqslant K.
$$
\endproclaim
\proclabel{PropositionInteriorSecondOrderBounds}
\proof By convexity, $\|A\|\leqslant H\leqslant n\|A\|$. If $H$ acheives its maximum on the boundary, then, by Lemma
\procref{LemmaSecondOrderBoundaryEstimates}:
$$
\|A\|\leqslant H\leqslant nK.
$$
\noindent If $H$ acheives its maximum in the interior, then, by Lemma \procref{LemmaMeanCurvatureIsSubharmonic} and the Maximum Principal:
$$
\|A\|\leqslant H\leqslant D.
$$
\noindent The result follows.\qed
\goodbreak
\newsubhead{The Dirichlet Problem I}
\noindent Let $M$ be an $(n+1)$-dimensional manifold. Choose $\theta\in[(n-1)\pi/2,n\pi/2[$. Let $H\subseteq M$ be a smooth convex hypersurface such that:
$$
R_\theta(H) = R_0.
$$
\noindent Where $R_0$ is constant. Let $\Omega\subseteq H$ be a bounded open subset of $H$. Let $\hat{\Sigma}\subseteq M$ be a convex hypersurface such that $\partial\hat{\Sigma}=\partial\Omega=:\Gamma$ and such that $R_\theta(\hat{\Sigma})\geqslant R_1$ in the weak sense, where $R_1\leqslant 1/\opTan^{-1}(\theta/n)$.
\proclaim{Proposition \nextprocno}
\noindent Let $f:\Omega\rightarrow\Bbb{R}$ be a smooth function. Define $\hat{f}:\Omega\rightarrow M$ by:
$$
\hat{f}(p) = \opExp(f(p)\msf{N}_H(p)),
$$
\noindent where $\msf{N}_H$ is the unit exterior normal over $H$. If $II$ denotes the second fundamental form of the graph of $f$, then:
$$
\hat{f}^* II = -B^{-1}\opHess(f)B^{-1} + R,
$$
\noindent where $B$ is a symmetric positive definite matrix, $R$ is a symmetric $2$-form, and $B$ and $R$ are functions only of $p$, $f$ and $\nabla f$.
\endproclaim
\proclabel{PropositionFormulaForSecondFF}
\proof Define $\Phi:\Omega\times[0,\infty[\rightarrow M$ by:
$$
\Phi(p,t) = \opExp(t\msf{N}_H(p)).
$$
\noindent Let $g_0$ be the Riemannian metric over $\Omega$. Since $M$ is a Hadamard manifold, $\Phi$ is a local diffeomorphism. Let $g$ be the Riemannian metric over $M$. $\hat{g}=\Phi^*g$ defines a Riemannian metric over $\Omega\times[0,\infty[$. With respect to $\hat{g}$, $\partial_t$ has unit length and is orthogonal to $T\Omega$. Let $M(p,t)$ be a symmetric matrix such that, for all vectors tangent to $T\Omega$:
$$
\hat{g}(X,Y) = g_0(M(p,t)X,Y).
$$
\noindent Let $\nabla_0$ be the Levi-Civita covariant derivative over $\Omega$. Let $\nabla$ be the Levi-Civita covariant derivative of $\hat{g}$. For $X$ tangent to $\Omega$, define $\hat{X}_f$ by:
$$
\hat{X}_f = (X,\langle\nabla^0f,X\rangle).
$$
\noindent Define the symmetric matrix $B:=B(p,f,\nabla f)$ such that, for all $X,Y\in T\Omega$:
$$
g_0(B X,B Y) = \hat{g}(\hat{X}_f,\hat{Y}_f).
$$
\noindent Define $\hat{\msf{N}}:=\hat{\msf{N}}(p,f,\nabla f)$ by:
$$
\hat{\msf{N}}_f = (-M^{-1}(p,f)\nabla^0 f, 1).
$$
\noindent $\hat{\msf{N}}_f$ is an outward normal to the graph of $f$. Define $\mu_f :=\mu_f(p,f,\nabla f)$ by:
$$
\mu_f= \|\hat{\msf{N}}_f\|.
$$
\noindent For all $X\in T\Omega$:
$$
\nabla_{\hat{X}_f}\hat{\msf{N}}_f = (-M^{-1}(p,f)\nabla^0_X \nabla^0(f), 0) + R_1(p,f,\nabla f)(X),
$$
\noindent where $R_1$ is a term which only depends on $p$, $f$ and $\nabla f$. Thus:
$$
\hat{g}(\nabla_{\hat{X}_f}\hat{\msf{N}}_f,\hat{Y}_f) = (-\nabla^0_X \nabla^0 f, 0) + R_2(p,f,\nabla f)(X,Y).
$$
\noindent It follows that:
$$
II = -\frac{1}{\mu_f}B^{-1}\opHess(f)B^{-1} + R_3.
$$
\noindent Where $\mu_f$, $B$ and $R_3$ only depend on $p$, $f$ and $\nabla f$. The result follows.\qed
\medskip
\noindent We now prove Theorems \procref{TheoremDirichletI}, \procref{TheoremDirichletIA} and \procref{TheoremDirichletIB}:
\medskip
{\noindent\bf Proof of Theorem \procref{TheoremDirichletI}:\ } Suppose that $\theta>(n-1)\pi/2$ and that $\hat{\Sigma}$ and $\Gamma$ are smooth. By Lemma \procref{LemmaSmoothingConvexSets}, the general case follows by approximation. Let $I\subseteq [R_0,R_1]$ be such that, for all $t\in I$, a solution exists which is smooth over $\overline{\Omega}$ and which lies strictly below $\hat{\Sigma}$ and whose supporting tangents along $\Gamma$ also lie strictly below those of $\hat{\Sigma}$. By definition $R_0\in I$. By Theorem $1.3$ of \cite{SmiSLC}, noting that $r>1/\opTan^{-1}(\theta/n)$, $I$ is open.
\medskip
\noindent Let $(r_n)_\ninn\in I$ be an increasing sequence converging to $r_0\in[R_0,R_1]$. For all $n$, let $\Sigma_n$ be a solution with $R_\theta(\Sigma_n)=r_n$. For all $n$, let $f_n:\overline{\Omega}\rightarrow\Bbb{R}$ be the function of which $\Sigma_n$ is the graph. By Lemma \procref{LemmaCompactnessOfConvexGraphs}, there exists $f_0$ to which $(f_n)_\ninn$ subconverges in the $C^{0,\alpha}$ sense for all $\alpha$ and whose graph lies below $\hat{\Sigma}$.
\medskip
\noindent Proposition \procref{PropositionInteriorSecondOrderBounds} yields uniform $C^2$ bounds on $(f_n)_\ninn$. For all $n$, $f_n$ satisfies an equation of the form:
$$
F(p,\phi,D\phi,D^2\phi;r,\theta) = 0.
$$
\noindent Since $f_n$ is uniformly bounded in the $C^2$ sense, $F$ is uniformly elliptic. By concavity of $R_\theta$ and Proposition \procref{PropositionFormulaForSecondFF}, $F$ is concave with respect to $D^2\phi$. Theorem $1$ of \cite{CaffNirSprII} therefore yields uniform $C^{2,\alpha}$ bounds on $(f_n)_\ninn$ for all $\alpha$. Repeated application of Schauder's estimates then yield uniform $C^k$ bounds on $(f_n)_\ninn$ for all $k$. It follows that $(f_n)_\ninn$ converges to $f_0$ in the $C^\infty$ sense over $\overline{\Omega}$.
\medskip
\noindent Thus, if $\Sigma_0$ is the graph of $f_0$, $\Sigma_0$ is smooth up to the boundary and $R_\theta(\Sigma_0)=r_0$. By part $(iv)$ of Lemma \procref{LemmaCompactnessOfConvexGraphs}, $\Sigma_0$ lies strictly below $\hat{\Sigma}$ and every supporting tangent to $\Sigma_0$ along $\Gamma$ also lies strictly below those of $\hat{\Sigma}$. $I$ is therefore closed, and existence for $r\in[R_0,R_1[$ follows by connectedness of the interval $[R_0,R_1]$. The case where $r=R_1$ is proven by taking limits, and the result follows.\qed
\medskip
{\noindent\bf Proof of Theorem \procref{TheoremDirichletIA}:\ }By Theorem \procref{TheoremDirichletI}, it remains to prove uniqueness. Choose $\theta\in[(n-1)\pi/2,n\pi/2[$ and $r\in[0,R_1]$. Let $\Sigma_1$ and $\Sigma_2$ be two distinct solutions such that:
$$
R_\theta(\Sigma_1),R_\theta(\Sigma_2) = r.
$$
\noindent Suppose that there exists $p\in\Sigma_2$ which lies strictly above $\Sigma_1$. By deforming $\Sigma_1$ slightly and moving it upwards by isometries of hyperbolic space, we obtain an immersed hypersurface, $\Sigma_1'$ and a point $p'\in\Sigma_1'$ such that $R_\theta(\Sigma_1)<r$ and $\Sigma_2$ is an interior tangent to $\Sigma_1'$ at $p'$. This contradicts the Geometric Maximum Principal. $\Sigma_2$ therefore lies below $\Sigma_1$. Likewise, $\Sigma_1$ lies below $\Sigma_2$ and they therefore coincide. The result follows.\qed
\medskip
{\noindent\bf Proof of Theorem \procref{TheoremDirichletIB}:\ } Let $U'$ be the intersection of all horoballs in $\Bbb{H}^{n+1}$ containing $\Omega$. Let $U$ be the intersection of $U'$ with one of the connected components of $\Bbb{H}^{n+1}\setminus H$. Define $\hat{\Sigma} = \partial U'$. $\hat{\Sigma}$ satisfies the hypotheses of Theorem \procref{TheoremDirichletIA} and the result now follows.\qed
\goodbreak
\newhead{The Perron Method}
\newsubhead{Extended Normals}
\noindent Let $M$ be a Hadamard manifold of sectional curvature bounded above by $-1$. Let $UM$ be the unitary bundle of $M$. Let $\Sigma=(S,i)$ be a smooth convex immersed hypersurface in $M$. Let $\msf{N}_\Sigma$ be the outward pointing normal over $\Sigma$. Let $\Omega$ be an open subset of $\Sigma$. Let $\msf{N}_{\partial\Omega}$ be the outward pointing normal of $\partial\Omega$ in $\Sigma$. Define $N\Omega$ and $N\partial\Omega$ by:
$$\matrix
N\Omega \hfill&= \left\{\msf{N}(p)\text{ s.t. }p\in\Omega\right\},\hfill\cr
N\partial\Omega \hfill&=\left\{\msf{V}_p\text{ s.t. }p\in\partial\Omega\text{ \& }\langle V_p,\msf{N}_{\partial\Omega}(p)\rangle,\langle V_p,\msf{N}_\Sigma(p)\rangle\geqslant 0\right\}.\hfill\cr
\endmatrix$$
\noindent We define $\hat{N}\Omega$ by:
$$
\hat{N}\Omega = N\Omega\munion N\partial\Omega.
$$
\noindent We call $\hat{N}\Omega$ the extended normal of $\Omega$. $\Omega$ embeds naturally as an open subset of $\hat{N}\Omega$. Moreover, $i$ extends naturally to an immersion $\hat{\mathi}:\hat{N}\Omega\rightarrow UM$. We define $\Phi:\hat{N}\Omega\times[0,\infty[\rightarrow M$ by:
$$
\Phi(p,t) = \opExp(t\hat{\mathi}(p)).
$$
\noindent Since $M$ is a Hadamard manifold and $\Sigma$ is convex, for every $p\in\hat{N}\Omega$ there exists a neighbourhood $U$ of $p$ in $\hat{N}\Omega$ such that the restriction of $\Phi$ to $U\times]0,\infty[$ is a homeomorphism onto its image. We refer to $\Phi$ as the end of $\Omega$. The differential structure of $M$ pulls back through $\Phi$ to a differential structure over $\hat{N}\Omega\times]0,\infty[$, which we also refer to as the end of $\Omega$ when there is no ambiguity. We denote it by $\Cal{E}(\Omega)$. $\Cal{E}(\Omega)$ is foliated by the geodesics normal to $\hat{\Omega}$. We refer to this foliation and the resulting vector field as the vertical foliation and vector field respectively.
\medskip
\noindent The boundary of $\partial\Cal{E}(\Omega)$ divides into two parts, which we denote by $\partial_w\Cal{E}(\Omega)$ and $\partial_f\Cal{E}(\Omega)$ and define as follows:
$$\matrix
\partial_w\Cal{E}(\Omega) \hfill&= \left\{(p,t)\text{ s.t. }p\in\partial\hat{N}\Omega\text{ \& }t\in[0,\infty[\right\},\hfill\cr
\partial_f\Cal{E}(\Omega) \hfill&= \left\{(p,0)\text{ s.t. }p\in\hat{N}\Omega\right\}.\hfill\cr
\endmatrix$$
\noindent We refer to $\partial_w\Cal{E}(\Omega)$ and $\partial_f\Cal{E}(\Omega)$ as the wall and the floor respectively of $\Cal{E}(\Omega)$. Trivially, $\partial_f\Cal{E}(\Omega)$ is identified with $\hat{N}\Omega$. The experienced reader will be aware that $\Cal{E}(\Omega)$ also has an ideal boundary at infinity. This will not concern us.
\medskip
\noindent Let $U\subseteq\Cal{E}(\Omega)$ be an open set. We say that $U$ is convex if and only if the shortest geodesic in $\Cal{E}(\Omega)$ joining any two points in $U$ also lies in $U$. We say that $U$ is boundary convex if and only if, for every boundary point $p\in\partial U$, there exists $r>0$ such that $B_r(p)\minter U$ is convex.
\medskip
\noindent If $U\subseteq\Cal{E}(\Omega)$ is boundary convex, let $\delta_U$ be the distance function to $U$ in $\Cal{E}(\Omega)$.
\proclaim{Proposition \nextprocno}
\noindent Let $p\in\partial U$ and $r>0$ be such that:
$$
d(p,\partial\Cal{E}(\Omega))>2r.
$$
\noindent Then, $\delta_U$ is convex in $B_r(p)$. In particular, $U\minter B_r(p)$ is convex.
\endproclaim
\proof For any $q\in B_r(p)$, the shortest geodesic joining $q$ to $U$ does not intersect $\partial\Cal{E}(\Omega)$. The first assertion now follows from the fact that the distance to a convex set in a Hadamard manifold is a convex function. Since $B_r(p)$ is convex, the second assertion follows.\qed
\medskip
\noindent Let $(U_n)_\ninn$ be a sequence of compact boundary convex subsets of $\Cal{E}(\Omega)$. We say that $(U_n)_\ninn$ is nested if and only if for all $i>j$, $U_i\subseteq U_j$. By classical point set topology, there exists a compact subset $U_0\subseteq\Cal{E}(\Omega)$ such that $(U_n)_\ninn$ converges to $U_0$ in the Haussdorf sense.
\proclaim{Proposition \nextprocno}
\noindent $U_0\setminus\partial\Cal{E}(\Omega)$ is boundary convex away from $\partial\Cal{E}(\Omega)$.
\endproclaim
\proof Choose $p\in\partial U_0\setminus\partial\Cal{E}(\Omega)$. There exists $r>0$ such that $d(p,\partial\Cal{E}(\Omega))>3r$. There exists $(p_n)_\ninn\in\Cal{E}(\Omega)$ which converges to $p$ such that, for all $n$:
$$
p_n\in\partial U_n.
$$
\noindent and:
$$
d(p_n,\partial\Cal{E}(\Omega))>2r.
$$
\noindent For all $n\in\Bbb{N}\munion\left\{0\right\}$, define $\delta_n=\delta_{U_n}$. By the preceeding proposition, the restriction of $\delta_n$ to $B_r(p_n)$ is convex. Taking limits, it follows that $\delta_0$ is convex, and so $U_0\minter B_r(p)$ is convex.\qed
\proclaim{Proposition \nextprocno}
\noindent Choose $p_0\in U_0\setminus\partial\Cal{E}(\Omega)$. Let $(p_n)_\ninn$ be a sequence converging to $p_0$ such that $p_n\in\partial U_n$ for all $n$. For all $n$, let $\msf{N}_n$ be a supporting normal to $U_n$ at $p_n$. If
$\msf{N}_n$ converges to $\msf{N}_0$, then $\msf{N}_0$ is a supporting normal to $U_0$ at $p_0$.
\endproclaim
\proof This follows from the analogous result in a Hadamard manifold.\qed
\medskip
\noindent We now say that $U$ is boundary $\epsilon$-convex if and only if, for every boundary point $p\in\partial U$, there exists $r>0$ such that $B_r(p)\minter U$ is $\epsilon$-convex.
\proclaim{Lemma \nextprocno}
\noindent Let $M$ be a Hadamard manifold of sectional curvature bounded above by $-1$. Choose $\epsilon>0$. Let $\Sigma\subseteq M$ be a convex hypersurface whose second fundamental form is bounded below by $\epsilon\opId$ in the weak sense. Let $d_\Sigma$ be the distance to $\Sigma$ in $M$. Let $\gamma:\Bbb{R}\rightarrow M$ be a curve lying on the outside of $\Sigma$ of geodesic curvature is less than $\epsilon$. For $d_{\Sigma}<\epsilon^{-1}$, $d_{\Sigma}^2$ is a convex function of $\gamma$.
\endproclaim
\proclabel{ConvexityOfDistanceFunction}
\proof Choose $t\in\Bbb{R}$. Let $p\in\Sigma$ be the closest point to $\gamma(t)$. Let $\epsilon'<\epsilon$ be greater than the geodesic curvature of $\gamma$ near $t$. Let $d$ be the distance to $\Sigma$ at $\gamma(t)$. Let $\eta:[0,d]\rightarrow M$ be the shortest geodesic segment from $p$ to $\gamma(t)$. By definition of $\Sigma$, there exists a strictly convex hypersurface, $\Sigma'$, which is an exterior tangent to $\Sigma$ at $p$ whose second fundamental form equals $\epsilon'\opId$ at $p$.
\medskip
\noindent For $s>0$, let $\Sigma'_s$ be the hypersurface at constant distance $s$ from $\Sigma'$. Let $A'_s$ be the second fundamental form of $\Sigma'_s$ at $\eta(s)$. By Lemma \procref{LemmaNormalDerivative}:
$$
\nabla_{\partial_s}A_s = W_s - A_s^2,
$$
\noindent where $W_s$ is such that, for all $X$ tangent to $\Sigma_s$:
$$
W_s(X,X) = \langle R_{\partial_sX}\partial_s,X\rangle.
$$
\noindent For all $s$, by definition of $M$:
$$
W_s \geqslant \opId.
$$
\noindent Thus:
$$
A_s \geqslant \opTanh(s+\hat{\epsilon})\opId,
$$
\noindent where $\hat{\epsilon}'$ is given by:
$$
\opTanh(\hat{\epsilon}') = \epsilon'.
$$
\noindent Let $d_{\Sigma'}$ be the distance to $\Sigma'$. Then, along $\eta$:
$$
\opHess(d_\Sigma') = \opTanh(s+\hat{\epsilon}')\opId^\perp,
$$
\noindent where:
$$
\opId^\perp = (\opId - \nabla d_{\Sigma'}\otimes \nabla d_{\Sigma'}).
$$
\noindent Thus:
$$\matrix
&\opHess({d_\Sigma'}^2) \hfill&= 2d_\Sigma'\opTanh(s+\hat{\epsilon}')\opId^\perp + 2\nabla d_{\Sigma'}\otimes \nabla d_{\Sigma'}\hfill\cr
\Rightarrow&\opHess({d_\Sigma'}^2)(X,X) \hfill&\geqslant 2\opMin(d_\Sigma'\opTanh(s+\hat{\epsilon}'),1)\|X\|^2.\hfill\cr
\endmatrix$$
\noindent Thus, along $\gamma$ at $t$:
$$\matrix
(\partial_t^2 (d_{\Sigma'}\circ\gamma)^2)\hfill&\geqslant 2\opMin(d\opTanh(d+\hat{\epsilon}'),1) - 2d\langle\nabla d_{\Sigma'},\nabla_{\partial_t\gamma}\partial_t\gamma\rangle\hfill\cr
&\geqslant 2\opMin(d\opTanh(d+\hat{\epsilon}'),1) - 2d\epsilon'\hfill\cr
&\geqslant 2\opMin(d\opTanh(d+\hat{\epsilon}') - d\epsilon',0)\hfill\cr
&=0.\hfill\cr
\endmatrix$$
\noindent It follows that $(d_{\Sigma'}\circ\gamma)^2$ is convex at $t$. Since:
$$
(d_{\Sigma'}\circ\gamma)^2 \geqslant (d_{\Sigma}\circ\gamma)^2,
$$
\noindent and since both functions are equal at $t$, the result now follows.\qed
\proclaim{Proposition \nextprocno}
\noindent If $(U_n)_\ninn$ is $\epsilon$-boundary convex for all $n$, then $U_0$ is also $\epsilon$-boundary convex away from $\partial\Cal{\Omega}$.
\endproclaim
\proclabel{PropositionConvexityInTheLimit}
\proof This follows by a similar reasoning as before, this time using Lemma \procref{ConvexityOfDistanceFunction} instead of the convexity of the distance from a geodesic to a convex set.\qed
\goodbreak
\newsubhead{Graphs Over Extended Normals}
\noindent We extend the notion of graphs to extended normals. Let $\Sigma=(S,i)$ be a convex immersed submanifold. We say that $\Sigma$ is a graph over $\hat{N}\Omega$ if there exists:
\medskip
\myitem{(i)} a relatively compact open subset $\Omega_\Sigma\subseteq\hat{N}\Omega$ such that $\Omega\subseteq\Omega_\Sigma$;
\medskip
\myitem{(ii)} a homeomorphism $\alpha:\Sigma\rightarrow\Omega_\Sigma$; and
\medskip
\myitem{(iii)} a continuous function $f:\Omega_\Sigma\rightarrow[0,\infty[$,
\medskip
\noindent such that $f$ vanishes along $\partial\Omega_\Sigma$, and for all $p\in\Sigma$:
$$
i(p) = \opExp(f(p)(\msf{N}\circ\alpha)(p)).
$$
\noindent We call $f$ and $\Omega_\Sigma$ the graph function and the graph domain respectively of $\Sigma$. We define $U_f$ by:
$$
U_f = \left\{(p,t)\text{ s.t. }p\in\Omega_\Sigma\text{ \& }t\leqslant f(p)\right\}.
$$
\noindent By definition of $\Sigma$, $U_f$ is a boundary convex subset of $\Cal{E}(\Omega_\Sigma)$.
\medskip
\noindent Let $\Sigma$ and $\Sigma'$ be two graphs over $\hat{N}\Omega$. Let $f$, $f'$ and $\Omega$, $\Omega'$ be their respective graph functions and graph domains. We define the partial order ``$\geqslant$'' over the space of graphs over $\hat{N}\Omega$ such that $\Sigma\geqslant\Sigma'$ if and only if:
$$
U_{f'} \subseteq U_f.
$$
\noindent In other words, if and only if $\Omega'\subseteq\Omega$ and:
$$
f|_{\Omega'} \geqslant f'.
$$
\noindent If $\Sigma'\leqslant\Sigma$, then we say that $\Sigma'$ lies below $\Sigma$.
\medskip
\noindent If $\Sigma$ is a graph over $\hat{N}\Omega$, we define $\opVol(\Sigma)$ to be the volume of $U_{\Sigma}$. By compactness, this is finite. Trivially:
$$
\Sigma'\leqslant\Sigma \Rightarrow \opVol(\Sigma')\leqslant\opVol(\Sigma).
$$
\noindent Moreover equality holds in the above relation if and only if $\Sigma=\Sigma'$.
\proclaim{Lemma \nextprocno}
\noindent Let $\Sigma_1>\Sigma_2>...$ be a decreasing sequence of $\epsilon$-convex immersed hypersurfaces which are graphs over $\hat{N}\Omega$. For all $i$, let $f_i$ be the graph function of $\Sigma_i$. There exists an $\epsilon$-convex immersed hypersurface $\Sigma_0$ such that:
\medskip
\myitem{(i)} for all $i$, $\Sigma_i>\Sigma_0$;
\medskip
\myitem{(ii)} $\Sigma_0$ is a graph over $\hat{N}\Omega$; and
\medskip
\myitem{(iii)} if $f_0$ is the graph function of $\Sigma_0$ over $\hat{N}\Omega$, then $f_0$ is $C^{0,1}_\oploc$, and $(f_n)_\ninn$ converges to $f_0$ in the $C^{0,\alpha}_\oploc$ sense for all $\alpha$.
\endproclaim
\proclabel{LemmaCompactnessOfConvexSubsetsOfHyperbolicEnds}
\remark Even without $\epsilon$-convexity, the graph function of the limit would still be $C^{0,1}$ over $\Omega$ and the graph functions would also converge accordingly over this set. The $\epsilon$-convexity is required to ensure that the limit function is also $C^{0,1}$ over $\Omega_0\setminus\Omega$.
\medskip
\proof Trivially $(f_n)_\ninn$ is a decreasing sequence. There thus exists $f_0$ to which this sequence converges pointwise. For all $n\in\Bbb{N}\munion\left\{0\right\}$, denote $U_n := U_{f_n}$. Trivially, $(U_n)_\ninn$ is a nested sequence. $(U_n)_\ninn$ therefore converges to $U_0$ in the Haussdorf sense. Since $U_n$ is a graph over $\hat{N}\Omega$ for every $n$, Proposition \procref{PropositionConvexityInTheLimit} may be modified to show that $U_0$ is boundary $\epsilon$-convex at every $p\in\partial U_0$ which does not lie in $\partial_w\Cal{E}(\Omega)$.
\medskip
\noindent For $p\in\hat{N}\Omega$ and for $n\in\Bbb{N}\munion\left\{0\right\}$, define $\hat{p}_n\in\Cal{E}(\Omega)$ by:
$$
\hat{p}_n = (p,f_n(p)).
$$
\noindent For $\epsilon>0$, define $\Omega_\epsilon$ by:
$$
\Omega_\epsilon = \left\{p\in\hat{N}\Omega\text{ s.t. }f_0(p)\geqslant\epsilon\right\}\munion\Omega.
$$
\noindent There exists $r_\epsilon>0$ such that, for every $p\in\Omega_\epsilon$ and for all $n$:
$$
d(\hat{p}_n,\partial_w\Cal{E}(\Omega)) > 2r_\epsilon.
$$
\noindent For all $n$, the supporting tangents to $\Sigma_n$ over $\Omega_\epsilon$ are uniformly bounded away from the vertical vector field. Indeed, suppose the contrary, then there exists $p\in\overline{\Omega}_\epsilon$ such that the vertical vector at $\hat{p}_0$ is tangent to $\Sigma_0$ at $\hat{p}_0$. However, by continuity, the geodesic segment joining $p$ to $\hat{p}$ lies in $U_0$. This is a absurd, since $\Sigma_0$ is $\epsilon$-convex at $p$.
\medskip
\noindent Let $\gamma:I\rightarrow\Omega_\epsilon$ be a rectifiable curve. For all $n\in\Bbb{N}\munion\left\{0\right\}$, let $\gamma_n:I\rightarrow\Sigma_n$ be the lift of $\gamma$. Since $(f_n)_\ninn$ is unformly bounded over $\Omega_\epsilon$, and since its slope is uniformly bounded, there exists $B_\epsilon$, independant of $\gamma$, such that, for all $n\in\Bbb{N}$:
$$
\opLength(\gamma_n) \leqslant B_\epsilon\opLength(\gamma).
$$
\noindent It follows that the sequence $(\hat{f}_n)_\ninn:=(p,f_n)_\ninn$ is uniformly Lipschitz over $\Omega_\epsilon$. Thus so is $(f_n)_\ninn$. Consequently $f_0$ is $C^{0,1}$ over $\Omega_\epsilon$ and $(f_n)_\ninn$ converges to $f_0$ in the $C^{0,\alpha}$ sense over $\Omega_\epsilon$ for all $\alpha$.
\medskip
\noindent Define $\Omega_0'$ by:
$$
\Omega_0' = f^{-1}(]0,\infty[)\munion\Omega.
$$
\noindent Let $\Omega_0$ be the connected component of $\Omega_0'$ containing $\Omega$. Then $f_0$ is $C^{0,1}_\oploc$ over $\Omega_0$ and $(f_n)_\ninn$ converges to $f_0$ in the $C^{0,\alpha}_\oploc$ sense over $\Omega_0$ for all $\alpha$. Moreover, $f_0$ extends to a continuous function over the closure of $\Omega_0$ which vanishes on $\partial\Omega_0$. The result follows.\qed
\medskip
\noindent The following lemma describes an important property of convex graphs which will be referred to as ``fatness'' in the sequel:
\proclaim{Lemma \nextprocno}
\noindent Let $\Sigma$ be an $\epsilon$-convex immersed hypersurface which is a graph over $\hat{N}\Omega$. Let $p\in\Sigma$ be an interior point. There exists $\eta>0$ and a supporting normal $\msf{N}_p$ to $\Sigma$ at $p$ such that, for any other supporting normal $\msf{N}_p'$ to $\Sigma$ at $p$:
$$
\langle \msf{N}_p,\msf{N}'_p\rangle \geqslant \eta.
$$
\endproclaim
\proclabel{LemmaConvexGraphsAreFat}
\proof Let $\Cal{N}_p$ denote the set of supporting normals to $\Sigma$ at $p$. $\Cal{N}_p$ is a convex subset of the sphere of unit vectors over $p$. Let $V_p$ be the vertical vector at $p$. Since $\Sigma$ is a strictly convex graph, there exists $\eta_1>0$ such that, for every supporting normal $\msf{N}_p\in\Cal{N}_p$:
$$
\langle \msf{N}_p, V_p\rangle \geqslant \eta_1.
$$
\noindent $\Cal{N}_p$ is thus strictly contained in the hemisphere about $V_p$. We denote this hemisphere by $H$. If $V_p\in\Cal{N}_p$, then the result follows with $\msf{N}_p=V_p$. Suppose therefore that $V_p\notin\Cal{N}_p$. By convexity, there exists a totally geodesic subsphere $S$, orthogonal to $\partial H$ such that $V_p\in S$ and $\Cal{N}_p$ lies strictly to one side of $S$ in $H$. Let $S'$ be obtained by rotating $S$ about $S\minter H$ until it meets $\Cal{N}_p$. Choose $\msf{N}_p\in S'\minter\Cal{N}_p$. $\msf{N}_p$ has the desired properties, and the result follows.\qed
\goodbreak
\newsubhead{The Dirichlet Problem II}
\noindent Let $M$ be an $(n+1)$-dimensional Hadamard manifold of sectional curvature bounded above by $-1$. Let $\Sigma$ be a smooth convex immersed hypersurface in $M$. Let $\Omega\subseteq\Sigma$ be an open subset and let $\hat{\Sigma}$ be a convex immersed hypersurface in $M$ which is a graph over the extended normal of $\Omega$.
\proclaim{Proposition \nextprocno}
\noindent For all $p\in M$, for every normal vector $\msf{N}_p$ over $p$, for all $\theta\in[(n-1)\pi/2,n\pi/2[$, and for all sufficiently small $\epsilon>0$, there exists $\delta>0$ and an immersed hypersurface $\Sigma$ of radius $\delta$ about $p$ which is normal to $\msf{N}_p$ at $p$ such that:
$$
R_\theta(\Sigma) = \epsilon\opTan^{-1}(\theta/n)\text{ and\ }\|A_\Sigma\|\leqslant 2\epsilon,
$$
\noindent where $A_\Sigma$ is the second fundamental form of $\Sigma$.
\endproclaim
\proclabel{PropositionAdaptedDiscs}
\remark Such disks will be refered to as $\delta$-adapted disks. They are important for the use of the Perron method in the proof of Theorem \procref{TheoremDirichletII}.
\medskip
\proof We use the Implicit Function Theorem for elliptic operators. Let $\Sigma_0$ be an immersed hypersurface in $M$ which is normal to $\msf{N}_p$ at $p$ such that:
$$
A_0 = \epsilon\opId,
$$
\noindent where $A_0$ is the second fundamental form of $\Sigma_0$.
\medskip
\noindent Let $\msf{N}_0$ be the normal vector field over $\Sigma_0$. Let $f:\Sigma_0\rightarrow\Bbb{R}$ be a smooth function representation an infinitesimal normal deformation of $\Sigma_0$. Then by Lemma $3.1$ of \cite{SmiSLC}:
$$
D\opSL_\theta\cdot f = -\Delta^B f + gf,
$$
\noindent where $g$ is a bounded function. For $\epsilon>0$ sufficiently small, $A_0$ is bounded above and below over $B_\delta(p)$. Moreover by Lemma $3.1$ of \cite{SmiSLC}, since the sectional curvature of $M$ is bounded above by $-1$, for $\epsilon$ sufficiently small $g>0$. There thus exists $K>0$ such that, for all smooth $f$ of compact support:
$$
\langle D\opSL_\theta\cdot f,f\rangle_{L^2} \geqslant K\|f\|^2_{L^2}.
$$
\noindent Thus, if $G:C^\infty(B_\epsilon(p))\rightarrow C_0^\infty(B_\epsilon(p))$ is the Green's operator of $D\opSL_\theta$, then:
$$
\|G\| \leqslant K.
$$
\noindent The radius in $W_{2,p}$ over which the inverse of $\opSL_\theta$ is defined is determined by the norms of $G$, $D\opSL_\theta$ and $D^2\opSL_\theta$. It is thus uniformly bounded below as the radius, $\delta$, tends to $0$. However, the $W_{2,p}$ distance between $\opSL_\theta(\Sigma_0)$ and the constant function tends to $0$ as $\delta$ tends to $0$. Thus, for $\delta$ sufficiently small, the Implicit Function Theorem yields an immersed hypersurface of constant special Lagrangian curvature. This reasoning can be adapted to ensure that the resulting hypersurface passes through $p$, is normal to $\msf{N}_p$ at $p$ and has second fundamental form colinear to $\opId$ at $p$. The result follows by reducing $\delta$ further if necessary.\qed
\medskip
\noindent We now prove Theorem \procref{TheoremDirichletII}:
\medskip
{\bf\noindent Proof of Theorem \procref{TheoremDirichletII}:\ }We suppose that $\theta>(n-1)\pi/2$. The case $\theta=(n-1)\pi/2$ follows by approximation. We consider first the case where $r>R_0$. By Lemma \procref{LemmaSLCIsLowestEval}, there exists $\epsilon_0>0$ which only depends on $r$ and $\theta$ such that $\hat{\Sigma}$ is $\epsilon_0$-convex. Let $\Sigma'$ be an $\epsilon$-convex immersion in $M$ which is a graph over $\hat{\Omega}$ such that $\Sigma'\leqslant\hat{\Sigma}$ and $R_\theta(\Sigma')\geqslant r$ in the weak sense. Let $p\in\Sigma'$ be an interior point. Let $\msf{N}_p$ be a supporting normal to $\Sigma'$ at $p$. By Lemma \procref{LemmaConvexGraphsAreFat} (fatness), $\msf{N}_p$ may be chosen such that for any other supporting normal $\msf{N}'_p$ to $\Sigma'$ at $p$:
$$
\langle\msf{N}_p,\msf{N}'_p\rangle \geqslant \epsilon_1,
$$
\noindent for some $\epsilon_1>0$. Let $0<\delta\ll\epsilon_0$ be small. Let $(\Sigma_\delta,\partial\Sigma_\delta)$ be a $\delta$-adapted disk with normal $\msf{N}_p$. By $\epsilon_0$-convexity, $\partial\Sigma_\delta$ lies above $\Sigma'$ for $\delta$ sufficiently small. Let $(\Sigma_{\delta,t})$ be a family of inward deformations of $\Sigma_\delta$ (in the direction opposite to $\msf{N}_p$) such that $\partial\Sigma_{\delta,t}$ lies above $\Sigma'$ for all $t$. By making $\Sigma_\delta$ smaller if necessary, we may assume that, for sufficiently small $t$, $\Sigma_{\delta,t}$ still has constant $\theta$-special Lagrangian curvature.
\medskip
\noindent Since $\Sigma'$ is strictly convex, for sufficiently small $t$, there exists a non-trivial open subset $\Omega_t\subseteq\Sigma_{\delta,t}$ which is relatively compact with respect to $\Sigma_{\delta,t}$ and such that a portion of $\Sigma'$ is a graph over $\Omega$. We denote this portion by $\Sigma_{\delta,t}'$.
\medskip
\noindent By continuity of the supporting normal to a convex set, there exists $t_0>0$ such that, for $t<t_0$, the pair $(\Omega_t,\Sigma_{\delta,t}')$ satifies the hypotheses of Theorem \procref{TheoremDirichletI}. There thus exists a convex immersion $\Sigma^r_{\delta,t}$ such that:
\medskip
\myitem{(i)} $\Sigma^r_{\delta,t}$ is a graph over $\Omega_t$;
\medskip
\myitem{(ii)} $\Sigma^r_{\delta,t}$ lies beneath $\Sigma_{\delta,t}'$ as a graph over $\Omega_t$;
\medskip
\myitem{(iii)} $\Sigma^r_{\delta,t}$ is smooth away from the boundary; and
\medskip
\myitem{(iv)} $R_\theta(\Sigma^r_{\delta,t}) = r$.
\medskip
\noindent We define $\Sigma''_t$ by replacing the portion $\Sigma'_{\delta,t}$ of $\Sigma'$ with $\Sigma^r_{\delta,t}$. By Lemma \procref{LemmaIntersectionStillHasCurvatureBoundedBelow}, $R_\theta(\Sigma'')\geqslant r$ in the weak sense. Moreover, $\Sigma''_t$ can be chosen to vary continuously with $t$.
\medskip
\noindent Suppose that, for $t<t_0$, $\partial\Omega_t$ does not intersect $\partial\Sigma'=\partial\Omega$. Then, for all $t$, $\Sigma'_{\delta,t}$ is a strict graph over $\hat{N}\Omega$ which lies below $\hat{\Sigma}$. Indeed, suppose first that there exists $t_1<t_0$ and $p\in\Sigma^r_{\delta,t}$ which also lies in $\Omega$. Let $t_1$ be the first such time. Since $\partial\Omega_t$ does not intersect $\partial\Omega$, $p$ is an interior point. Since $t_1$ is the first intersection time, $\Sigma^r_{\delta,t}$ is an exterior tangent to $\Omega$ at this point. However, at $p$:
$$
R_\theta(\Sigma^r_{\delta,t}) = r > R_0.
$$
\noindent This is absurd, by the Geometric Maximum Principal.
\medskip
\noindent Suppose now that $\Sigma'_{\delta,t}$ is not a graph over $\hat{N}\Omega$. By continuity, there exists in interior point $p\in\Sigma'_{\delta,t}$ such that the vertical vector at $p$ is tangent to $\Sigma'_{\delta,t}$ at $p$. By continuity, the vertical geodesic segment joinging $\hat{N}\Omega$ to $p$ lies below $\Sigma'_{\delta,t}$. It follows that the vertical vector at $p$ is an interior tangent to $\Sigma'_{\delta,t}$ at $p$. This is impossible by strict convexity.
\medskip
\noindent We denote by $A$ the above described operation for obtaining new convex immersions out of old ones. Let $\Sigma_1$ and $\Sigma_2$ be two convex immersions which are graphs over $\hat{N}\Omega$ such that $\Sigma_1,\Sigma_2\leqslant \hat{\Sigma}$ and $R_\theta(\Sigma_1),R_\theta(\Sigma_2)\geqslant r$ in the weak sense. Let $f_1$ and $f_2$ be their respective graph functions. Define $f_{1,2}$ by:
$$
f_{1,2} = \opMin(f_1,f_2),
$$
\noindent Let $\Sigma_{1,2}$ be the graph of $f_{1,2}$. Trivially, $\Sigma_{1,2}\leqslant\hat{\Sigma}$. Moreover, by Lemma \procref{LemmaIntersectionStillHasCurvatureBoundedBelow}, $R_\theta(\Sigma_{1,2})\geqslant r$ in the weak sense. We denote this operation for obtaining new convex immersions out of old ones by $B$.
\medskip
\noindent Let $\Cal{F}$ be the family of all convex immersions which may be obtained from $\hat{\Sigma}$ by a finite combination of the operations $A$ and $B$. Define $V_0\geqslant 0$ by:
$$
V_0 = \minf\left\{ \opVol(\Sigma)\text{ s.t. }\Sigma\in\Cal{F}\right\}.
$$
\noindent There exists a sequence $\Sigma_1>\Sigma_2>...$ of strictly convex immersions in $\Cal{F}$ such that:
$$
\opVol(\Sigma_n)_\ninn\rightarrow V_0.
$$
\noindent For all $n\in\Bbb{N}$, let $f_n$ and $\hat{\Omega}_n$ be the graph function and graph domain of $\Sigma_n$ respectively. $(f_n)_\ninn$ is a decreasing sequence. By Proposition \procref{LemmaCompactnessOfConvexSubsetsOfHyperbolicEnds}, there exists $f_0:\hat{\Omega}_0\rightarrow[0,\infty[$ such that:
\medskip
\myitem{(i)} $f_0$ is continuous over the closure of $\hat{\Omega}_0$;
\medskip
\myitem{(ii)} $f_0$ vanishes along $\partial\hat{\Omega}_0$;
\medskip
\myitem{(iii)} $f_0$ is $C^{0,1}_\oploc$ inside $\hat{\Omega}_0$;
\medskip
\myitem{(iv)} $(f_n)_\ninn$ converges to $f_0$ in the $C^{0,\alpha}_\oploc$ sense over $\hat{\Omega}_0$ for all $\alpha$; and
\medskip
\myitem{(v)} if $\Sigma_0$ is the graph of $f_0$, then $\Sigma_0$ is $\epsilon$-convex.
\medskip
\noindent Let $p\in\Omega_0$ be an interior point. Let $\msf{N}_p$ be a supporting normal to $\Sigma_0$ at $\hat{p}$ chosen such that, for any other supporting normal $\msf{N}_p'$ at $\hat{p}$:
$$
\langle\msf{N}_p,\msf{N}'_p\rangle \geqslant \epsilon_1,
$$
\noindent for some $\epsilon_1>0$. For all $n$, let $d_n$ be the restriction to $\Sigma_n$ of the length metric of $\Cal{E}(\Omega)$. The construction outlined at the beginning of the proof may be carried out uniformly near $\hat{p}$ for all $n$. We thus obtain $t>0$ and for all $n$:
\medskip
\myitem{(i)} $\Omega_{t,n}$;
\medskip
\myitem{(ii)} $\Sigma_{n,\delta,t}$; and
\medskip
\myitem{(iii)} $\Sigma^r_{n,\delta,t}$,
\medskip
\noindent such that, for all $n$:
\medskip
\myitem{(i)} $\Sigma_{n,\delta,t}$ and $\Sigma^r_{n,\delta,t}$ are graphs above $\Omega_{t,n}$;
\medskip
\myitem{(ii)} $\Sigma^r_{n,\delta,t}$ lies below $\Sigma_{n,\delta,t}$; and
\medskip
\myitem{(iii)} $\Sigma^r_{n,\delta,t}$ has radius at least $\epsilon_2$ about $\hat{p}_n$ with respect to $d_n$ for some fixed $\epsilon_2>0$, where $\hat{p}_n$ is the point in $\hat{\Sigma}_{n,\delta,t}$ lying above $p_0$.
\medskip
\noindent For all $n$, we define $\Sigma'_n$ by replacing the portion $\Sigma_{n,\delta,t}$ of $\Sigma_n$ with $\Sigma^r_{n,\delta,t}$. For all $n$, $\Sigma'_n\in\Cal{F}$ and $\Sigma'_n\leqslant\Sigma_n$. Let $\Sigma'_0$ be the limit of $(\Sigma'_n)_\ninn$. Trivially:
$$
\Sigma'_0\leqslant\Sigma_0.
$$
\noindent We assert that $\Sigma'_0=\Sigma_0$. Indeed, otherwise, $\Sigma'_0\neq\Sigma_0$, in which case:
$$
\opVol(\Sigma'_0) < \opVol(\Sigma_0),
$$
\noindent which is absurd. By Theorem $1.4$ of \cite{SmiSLC}, it follows that $\Sigma_0$ is smooth over a radius of $\epsilon_2$ about $\hat{p}_0$. Since $p\in\hat{\Omega}_0$ is abitrary, it follows that $\Sigma_0$ is smooth over the whole of $\hat{\Omega}_0$. Moreover, $R_\theta(\Sigma_0)=r$, and the result follows for $r>R_0$.
\medskip
\noindent Let $\Sigma^r$ be the hypersurface obtained in this manner such that $R_\theta(\Sigma^r)=r$. Then, for all $r>r'$:
$$
\Sigma^r > \Sigma^{r'}.
$$
\noindent Thus, taking the limit as $r$ tends to $R_0$ yields the desired solution when $r=R_0$. The result follows.\qed
\goodbreak
\newhead{The Perron Method II}
\newsubhead{Pseudo-Immersions}
\noindent In order to prove Theorem \procref{TheoremDirichletIII}, we require a compactification of the space of convex immersions when there is no ambient end. To this end, we define pseudo-immersions.
\medskip
\noindent Let $M$ be an $(n+1)$-dimensional Hadamard manifold. Let $TM$ and $UM\subseteq TM$ be the tangent and unitary bundles respectively over $M$. Let $\pi:TM\rightarrow M$ be the canonical projection. Let $N$ be a compact $n$-dimensional manifold without boundary. A pseudo-immersion of $N$ into $M$ is a pair $(\varphi,\hat{\varphi})$ where:
\medskip
\myitem{(i)} $\varphi:N\rightarrow M$ is a $C^{0,1}$ mapping; and
\medskip
\myitem{(ii)} $\hat{\varphi}:N\rightarrow UM$ is an injective $C^{0,1}$ mapping,
\medskip
\noindent such that:
$$
\pi\circ\hat{\varphi} = \varphi.
$$
\noindent In the sequel, we will denote such a pair simply by $\varphi$. Since $\varphi$ is Lipschitz, the path metric and the volume of $M$ pull back to a (possibly degenerate) path metric and volume form over $N$, which we denote by $d_\varphi$ and $\opdVol_\varphi$ respectively. Likewise, the path metric of $UM$ pulls back to a path metric over $N$, which we denote by $\hat{d}_\varphi$. Since $\hat{\varphi}$ is injective, $\hat{d}_\varphi$ is non-degenerate. For $p\in N$, we denote the balls of radius $r$ in $N$ about $p$ with respect to $d_\varphi$ and $\hat{d}_\varphi$ by $B_r(p;N)$ and $\hat{B}_r(p;N)$ respectively. We denote these simply by $B_r(p)$ and $\hat{B}_r(p)$ respectively when there is no ambiguity concerning the ambient manifold.
\medskip
\noindent We say that a sequence $(\varphi_n,\hat{\varphi}_n)_\ninn$ converges to $(\varphi_0,\hat{\varphi}_0)$ if and only if $(\varphi_n)_\ninn$ and $(\hat{\varphi}_n)_\ninn$ converge to $\varphi_0$ and $\hat{\varphi}_0$ respectively in the $C^{0,\alpha}$ sense over $N$ for all $\alpha$.
\medskip
\noindent For $r>0$, if $\varphi$ is a pseudo-immersion and $p\in N$, we say that $\varphi$ is convex over a radius of $r>0$ at $p$, if and only if there exists a convex set $K\subseteq M$ such that:
\medskip
\myitem{(i)} $\varphi(p)\in\partial K$;
\medskip
\myitem{(ii)} $\hat{\varphi}(p)$ is normal to $K$ at $\varphi(p)$; and
\medskip
\myitem{(iii)} $\varphi(\hat{B}_r(p))\subseteq K$.
\medskip
\noindent We say that a pseudo-immersion, $\varphi$, is convex if and only if there exists $r>0$ such that $\varphi$ is convex over a radius $r$ at every point of $N$. For $\epsilon>0$, we define $\epsilon$-convexity in an analogous manner. Trivially, every convex immersion is also a convex pseudo-immersion and every $\epsilon$-convex immersion is also an $\epsilon$-convex pseudo-immersion.
\medskip
\noindent For a convex pseudo-immersion $\varphi$, we define the mapping $\Phi:N\times[0,\infty[\rightarrow M$ by:
$$
\Phi(p,t) = \opExp(t\hat{\varphi}(p)).
$$
\proclaim{Proposition \nextprocno}
\noindent For every $p\in N$, there exists a neighbourhood $p\in U\subset N$ such that the restriction of $\Phi$ to $U\times [0,\infty[$ is injective.
\endproclaim
\remark By conservation of the domain, the restriction of $\Phi$ to this set is then a homeomorphism onto its image.
\medskip
\remark In this case, we refer to $\Phi$ as the end of $\varphi$. We furnish the manifold $N\times]0,\infty[$ with the differential structure of $M$ pulled back through $\Phi$. We also refer to the resulting manifold as the end of $\phi$, and we denote it by $\Cal{E}(\varphi)$.
\medskip
\proof Let $r>0$ be such that $\varphi$ is convex over a radius $r$ about every point in $N$. Choose $q\in \hat{B}_{r/2}(p)$. Let $K_p$ and $K_q$ be convex sets as in the definition of convexity at $p$ and $q$. Let $K=K_p\minter K_q$. $K$ is also convex. Thus, if $\gamma$ is the geodesic segment joining $p$ to $q$, then $\gamma$ lies in $K$ and thus makes an optuse angle with any normal vector to $K$ at $p$ and $q$. Consequently, the half geodesics leaving $\varphi(p)$ and $\varphi(q)$ in the respective directions of $\hat{\varphi}(p)$ and $\hat{\varphi}(q)$ never intersect. Since $q\in \hat{B}_{r/2}(p)$ is arbitrary, the result follows.\qed
\medskip
\noindent Let $\varphi,\varphi':N\rightarrow M$ be two convex pseudo-immersions. We say that $\varphi'$ is a graph over $\varphi$ if and only if there exists a $C^{0,1}$ function $f:N\rightarrow [0,\infty[$ such that for all $p\in N$:
$$
\varphi'(p) = \opExp_{\varphi(p)}(f(p)\hat{\varphi}(p)).
$$
\noindent We observe that if $\varphi'$ is a graph over $\varphi$, then $\opdVol_{\varphi'}\geqslant\opdVol_\varphi$, with equality if and only if $\varphi'=\varphi$. We thus define the partial order ``$\leqslant$'' on the set of convex pseudo-immersions such that $\varphi\leqslant\varphi'$ if and only if $\varphi'$ is a graph over $\varphi$.
\proclaim{Proposition \nextprocno}
\noindent Let $\varphi:N\rightarrow M$ be an $\epsilon$-convex pseudo-immersion. Let $\Sigma\subseteq M$ be a convex immersed hypersurface such that:
\medskip
\myitem{(i)} the second fundamental form of $\Sigma$ is bounded above by $\epsilon$ in the weak sense;
\medskip
\myitem{(ii)} $\varphi(p)\in\Sigma$ and $\hat{\varphi}(p)$ is normal to $\Sigma$ at $\varphi(p)$; and
\medskip
\myitem{(iii)} $\Sigma$ has radius at most $\epsilon^{-1}$ about $\varphi(p)$.
\medskip
\noindent Then $\Sigma$ lifts to an immersed hypersurface in $\Cal{E}(\varphi)$.
\endproclaim
\proclabel{PropositionConvexCurvesLift}
\proof This follows from Lemma \procref{ConvexityOfDistanceFunction}.\qed
\proclaim{Proposition \nextprocno}
\noindent Let $K\subseteq M$ be compact. Choose $\epsilon>0$ and let $\varphi:N\rightarrow M$ be an $\epsilon$-convex pseudo-immersion such that $\varphi(N)\subseteq K$. There exists $r$, which only depends on $K$ and $\epsilon$ such that $\varphi$ is $\epsilon$-convex over a radius of $r$.
\endproclaim
\proclabel{PropositionUniformityOfEpsilonConvexity}
\remark This lemma makes $\epsilon$-convexity uniform over sequences ensuring that this property is preserved when limits are taken.
\medskip
\proof Choose $p\in N$. Let $\Sigma\subseteq M$ be a convex immersed hypersurface normal to $\hat{\varphi}(p)$ at $p$ such that the norm of its second fundamental form is bounded above by $\epsilon$. For $r>0$, let $\Sigma_r$ be the ball of radius $r$ in $\Sigma$ about $p$. By Proposition \procref{PropositionConvexCurvesLift}, $\Sigma_{r_1}$ lifts to an immersed hypersurface $\hat{\Sigma}$ in $\Cal{E}(\varphi)$ for some $r_1>0$. There exists $r_2>0$ be such that $\hat{B}_{r_2}(p,\Sigma)\subseteq \Sigma_{r_1}$. There exists a neighbourhood, $U$ of $\hat{\varphi}(p)\in UM$ such that every geodesic passing through $U$ intersects $\Sigma_{r_1}$ transversely. Consequently, the slope of $\hat{\Sigma}_{r_1}$ as a graph over $\varphi$ is uniformly bounded for $\varphi(q)\in U$. There therefore exists $K_1,r_3>0$ such that, if $\gamma$ is a curve in $\hat{B}_{r_3}(p,N)$ and $\hat{\gamma}$ is the curve in $\hat{\Sigma}$ lying above $\gamma$, then:
$$
\hat{l}(\gamma)/K_1 \leqslant \hat{l}(\hat{\gamma}) \leqslant K_1\hat{l}(\gamma),
$$
\noindent where $\hat{l}$ denotes length with respect to $\hat{d}$. Thus there exists $r_4>0$ such that a subset of $\hat{\Sigma}_r$ is a graph over $\hat{B}_{r_4}(p,N)$. In otherwords, for all $q\in\hat{B}_{r_4}(p,N)$, every half-geodesic leaving $\varphi(q)$ in the direction of $\hat{\varphi}(q)$ intersects $\hat{\Sigma}_{r_1}$ non-trivially.
\medskip
\noindent Let $\Omega$ be a convex set such that $\Sigma_{r_1}\subseteq\partial\Omega$. Define $d_\Omega:M\rightarrow[0,\infty[$ by:
$$
d_\Omega(q) = d(q,\Omega).
$$
\noindent $d_\Omega$ is a convex function over $M$. Choose $q\in B_{r_4}(p;N)$ and suppose that $\varphi(q)\notin\Omega$. Since the half-geodesic leaving $\varphi(q)$ in the direction of $\hat{\varphi}(q)$ intersects $\Omega$ non-trivially, and since $d_\Omega$ is a convex function, at $q$:
$$
\langle \hat{\varphi}(q),\nabla d_\Omega(q)\rangle < 0.
$$
\noindent For sufficiently small $r_4$, this is not possible and we therefore obtain the desired value for $r$. Since this construction may be carried out uniformly for $\varphi(p)\in K$, the result follows.\qed
\proclaim{Lemma \nextprocno}
\noindent Let $\varphi:N\rightarrow M$ be a smooth strictly convex immersion. Let $(\varphi_n)_\ninn:N\rightarrow M$ be $\epsilon$-convex pseudo-immersions such that:
\medskip
\myitem{(i)} for all $n$, $\varphi_n\leqslant\varphi$;
\medskip
\myitem{(ii)} there exists $p\in N$ such that $(f_n(p))_\ninn$ is bounded; and
\medskip
\myitem{(iii)} there exists $\delta>0$ and $N\in\Bbb{N}$ such that for $n\geqslant N$:
$$
f_n \geqslant \delta.
$$
\noindent Then, there exists an $\epsilon$-convex immersion $\varphi_0:N\rightarrow M$ such that:
\medskip
\myitem{(i)} $\varphi_0\leqslant\varphi$; and
\medskip
\myitem{(ii)} $(\varphi_n)_\ninn$ subconverges to $\varphi_0$.
\medskip
\noindent Moreover, $(\hat{d}_n)_{n\in\Bbb{N}\munion\left\{0\right\}}$ is uniformly equivalent to $d$ over $N$.
\endproclaim
\proclabel{LemmaCompactnessOfConvexPseudoImmersions}
\proof Since $\varphi$ is smooth, $d_\varphi$ and $\hat{d}_\varphi$ are equivalent. Let $\gamma:I\rightarrow N$ be a curve. By convexity, for all $n$:
$$
\opLength(\varphi\circ\gamma) \geqslant \opLength(\varphi_n\circ\gamma).
$$
\noindent Thus, for all $n$, $\varphi_n$ is $1$-Lipschitz. Moreover, for all $p\in N$:
$$
f_n(p) = d(p,\varphi_n(p)).
$$
\noindent Thus, for all $n$ and for all $p,q\in N$:
$$
\left|f_n(p) - f_n(q)\right| \leqslant d(p,q) + d(\varphi_n(p),\varphi_n(q)) \leqslant 2d(p,q).
$$
\noindent Thus $f_n$ is $2$-Lipschitz for all $n$. Since $(f_n)_\ninn=(d(\varphi_n(p),\varphi(p)))_\ninn$ is uniformly bounded at one point, there exist $C^{0,1}$ functions $\varphi_0:N\rightarrow M$ and $f_0:N\rightarrow [0,\infty[$ to which $(\varphi_n)_\ninn$ and $(f_n)_\ninn$ respectively converge in the $C^{0,\alpha}$ sense over $N$ for all $\alpha$. By condition $(iv)$, $f_0>0$ over $N$.
\medskip
\noindent For all $n$, and for all $p\in N$:
$$
\hat{\varphi}_n(p) = \frac{1}{f_n(p)}\opExp^{-1}_{\varphi_n(p)}(\varphi(p)).
$$
\noindent There thus exists $\hat{\varphi}_0$ to which $(\hat{\varphi}_n)_\ninn$ converges in the $C^{0,\alpha}$ sense over $N$ for all $\alpha$. By Proposition \procref{PropositionUniformityOfEpsilonConvexity}, $\varphi_0$ is $\epsilon$-convex.
\medskip
\noindent Let $\gamma$ be a curve in $N$. Define $\tilde{\gamma}_n$ by:
$$
\tilde{\gamma}_n = (f_n\hat{\varphi}_n)(\gamma_n(t)).
$$
\noindent Since $f_n$ and $\hat{\varphi}_n$ are uniformly bounded in the $C^{0,1}$ sense, there exists $B>0$ such that, for all $n$:
$$
\opLength(\tilde{\gamma}_n) \leqslant B\opLength(\gamma)_n.
$$
\noindent Conversely, since the derivative of the exponential mapping is bounded over any compact subset of $TM$, and since $\opExp(\tilde{\gamma}_n) = \gamma_n$, by increasing $B$ if necessary, we obtain, for all $n$:
$$
\opLength(\gamma_n)\leqslant B\opLength(\tilde{\gamma}_n).
$$
\noindent For all $n$, and for all $p\in N$, define $\eta_{n,p}$ to be the geodesic leaving $\varphi_n(p)$ in the direction of $\hat{\varphi}_n(p)$. Since $\varphi$ is strictly convex and is a graph over $\varphi_n$, there exists $\epsilon>0$ such that, for all $p\in N$ and for all $n$:
$$
\langle \partial_t\eta_{n,p},\hat{\varphi}(p)\rangle \geqslant \epsilon.
$$
\noindent Indeed, let $B$ be a small ball lying on the outside of $\varphi$ and tangent to $\varphi$ at $p$. Let $N$ be such that, for $n\geqslant N$, $f_n(p)>0$. Since $\varphi$ is smooth, moving $B$ inwards slightly and intersecting with the interior of $\varphi$ yields a convex set $K_p$ lying in the end of $\varphi_n$ for all $n\geqslant N$. If $\gamma_n(p)$ is the geodesic leaving $\varphi_n(p)$ in the direction of $\hat{\varphi}_n(p)$, then there exists $\epsilon_1>0$ such that:
$$
\gamma_n(p)([f_n(p)-\epsilon,f_n(p)]) \subseteq K.
$$
\noindent This yields a sequence of geodesic segments with length uniformly bounded below. The assertion now follows, since otherwise, these segments would converge to a segment tangent to $\varphi$ at $p$, which is impossible, by the strict convexity of $\varphi$.
\medskip
\noindent Thus, the derivative of the projection onto $UM$ is uniformly bounded below along $\tilde{\gamma}$. So, by increasing $B$ again if necessary, we obtain, for all $n$:
$$
\frac{1}{B}\opLength(\gamma;\hat{d}_{\varphi_n}) \leqslant \opLength(\gamma;\hat{d}_\varphi) \leqslant B\opLength(\gamma;\hat{d}_{\varphi_n}).
$$
\noindent The result follows.\qed
\goodbreak
\newsubhead{The Isotopy Problem}
\noindent We now prove Theorem \procref{TheoremDirichletIII}:
\medskip
{\bf\noindent Proof of Theorem \procref{TheoremDirichletIII}:\ }Suppose first that $\theta>(n-1)\pi/2$. Let $\varphi:N\rightarrow M$ be the immersion. We may assume that $R_\theta(\varphi)\geqslant r$ in the weak sense. Indeed, let $\hat{\varphi}:N\rightarrow UM$ be the exterior normal over $N$. For $t\geqslant 0$, define $\varphi_t$ by:
$$
\varphi_t(p) = \opExp_{\varphi(p)}(t\hat{\varphi}(p)).
$$
\noindent Since the sectional curvature of $M$ is bounded above by $-1$, for all $\epsilon>0$, there exists $T>0$ such that for $t\geqslant T$, $\varphi_t$ is $(1-\epsilon)$-convex. In particular, for $\epsilon$ sufficiently small, $R_\theta(\varphi_t)\geqslant r$. We may thus replace $\varphi$ with $\varphi_t$ for sufficiently large $t$.
\medskip
\noindent By Lemma \procref{LemmaSmoothingConvexSets}, we may assume that $\varphi$ is smooth. Let $\varphi':N\rightarrow M$ be a convex immersion such that $\varphi\geqslant\varphi'$, and $R_\theta(\varphi')\geqslant r$ in the weak sense. Choose $p\in N$. By Propositions \procref{PropositionConvexCurvesLift} and \procref{PropositionAdaptedDiscs}, we may construct an adapted disk $(\Sigma,\partial\Sigma)$ at $p$ which is normal to $\hat{\varphi}(p)$ and which lifts to $\Cal{E}(\varphi')$. By choosing the norm of the second fundamental form of $\Sigma$ sufficiently small, we may assume that $\Sigma$ has negative curvature.
\medskip
\noindent Let $(\Sigma_t)_{t\in[0,\epsilon[}$ be a family obtained by moving $\Sigma$ downwards (in the direction opposite to $\hat{\varphi}'(p)$). For sufficiently small $\epsilon$, $\Sigma_t$ can be chosen to be adapted for all $t$. Moreover, the norm of the second fundamental form of $\Sigma$ may be chosen sufficiently small so that the intersection of $\varphi'(N)$ with $\Sigma$ is $\eta$-convex, for some $\eta>0$. Finally, we assume that $\partial\Sigma_t$ lies in $\Cal{E}(\varphi')$ for all $t$.
\medskip
\noindent Let $(\Omega_t)_{t\in [0,\epsilon[}$ be the continuous family of connected open subsets of $\Sigma_t$ defined such that $\Omega_0=\left\{p\right\}$ and $\partial\Omega_t = \varphi'(N)\minter\Sigma_t$. Let $\Sigma'_t$ be the portion of $\varphi'(N)$ lying above $\Omega_t$. We claim that, for all $t$, $\Omega_t$ is a convex open set with non-trivial interior. Indeed, suppose that $\Omega_t$ degenerates. By strict convexity, this is only possible if $\Omega_{t_0}$ is a single point for some $t_0>0$. By Lemma \procref{LemmaRemoveableSingularity}, $\varphi'(N)$ is the boundary of a convex set, and is therefore homotopically trivial, which contradicts the hypotheses. The assertion follows.
\medskip
\noindent We now claim that, for all $t$, $\Sigma'_t$ is a graph over the extended normal of $\Omega_t$. Indeed, suppose the contrary. By continuity and strict convexity, there exists $t_0>0$ such that, either, the graph of $\Sigma'_{t_0}$ is vertical over $\Omega_{t_0}$ at some interior point, or the outward normal of $\Sigma'_{t_0}$ points vertically downwards at some point on the boundary. The former case is excluded by strict convexity of $\Sigma'_t$. In the latter case, $\Omega_{t_0}$ is a single point, in which case $\varphi'(N)$ is the boundary of a convex set, which contradicts the hypotheses as before. The assertion follows.
\medskip
\noindent Choose $0<t<\epsilon$. By Theorem \procref{TheoremDirichletII}, there exists $\Sigma'_{t,r}$ which is smooth up to the boundary, and which is a graph over the extended normal of $\Omega_t$ lying between $\Omega_t$ and $\Sigma'_t$ such that:
$$
R_\theta(\Sigma'_{t,r}) = r.
$$
\noindent We define $\varphi''$ by replacing $\Sigma'_t$ in $\varphi'$ with $\Sigma'_{t,r}$. $\varphi''$ is a convex immersion and $\varphi''\leqslant\varphi'$. By Lemma \procref{LemmaRemoveableSingularity}, $R_\theta(\varphi'')\geqslant r$ in the weak sense. Moreover, by examening the proof of Theorem \procref{TheoremDirichletII}, if $\Sigma'_{t,r}$ is chosen to be the maximal solution (in the sense that its graph function is maximal), then $\varphi''$ is obtained from $\varphi'$ by isotopic deformation. In particular, this implies as before that $\varphi$ is a graph over $\varphi''$.
\medskip
\noindent Let $\Cal{F}$ be the family of convex immersions in $M$ obtained by a finite number of iterations of the operation described above. By Lemma $4.1$ of \cite{GuanSpruckI}, if $\varphi_1$ and $\varphi_2$ are two convex immersions in $\Cal{F}$, then there exists a third convex immersion $\varphi_{1,2}$ in $\Cal{F}$ such that $\varphi_1,\varphi_2\geqslant\varphi_{1,2}$. For $\varphi'\in\Cal{F}$, let $\opVol(\varphi')$ denote the volume between $\varphi'$ and $\varphi$ in the end of $\varphi'$. Define $V_0$ by:
$$
V_0 = \msup\left\{\opVol(\varphi')\text{ s.t. }\varphi'\in\Cal{F}\right\}.
$$
\noindent There exists a sequence $\varphi_1\geqslant\varphi_2\geqslant ...$ in $\Cal{F}$ such that:
$$
(\opVol(\varphi_n))_\ninn \rightarrow V_0.
$$
\noindent For all $n$, define $d_n$ by:
$$
d_n = \minf\left\{ f_n(p)\text{ s.t. }p\in N\right\}.
$$
\noindent We claim that $(d_n)_\ninn$ is bounded. Indeed, suppose the contrary. Since the sectional curvature of $M$ is bounded above by $-1$, by convexity:
$$
\opDiam(\varphi_n)\leqslant (\opLog(\opSinh(d_n)))^{-1}\opDiam(\varphi).
$$
\noindent Thus, as $(d_n)_\ninn\rightarrow\infty$, $\opDiam(\varphi_n)_\ninn\rightarrow 0$. This contradicts the hypotheses on $N$, and the assertion follows. In particular, $V_0$ is finite.
\medskip
\noindent Thus, by Lemma \procref{LemmaCompactnessOfConvexPseudoImmersions}, there exists an $\epsilon$-convex pseudo-immersion $\varphi_0:N\rightarrow M$ such that $\varphi_0<\varphi$ to which $(\varphi_n)_\ninn$ subconverges. Since $\varphi_0$ maximises volume, by an analogous reasoning to that used in the proof of Theorem \procref{TheoremDirichletII}, $\varphi_0$ is smooth and:
$$
R_\theta(\varphi_0) = r.
$$
\noindent By construction, $\varphi_0$ is isotopic to $\varphi$ and $(i)$ follows.
\medskip
\noindent Suppose now that $\theta =(n-1)\pi/2$. Let $(\theta_n)_\ninn$ be a decreasing sequence converging to $\theta$ and let $(r_n)_\ninn$ be a sequence converging to $r$. For all $n$, let $\varphi_n:N\rightarrow M$ be a smooth immersion such that $\varphi_n\leqslant\varphi$ and:
$$
R_{\theta_n}(\varphi_n) = r_n.
$$
\noindent By Theorem $1.4$ of \cite{SmiSLC}, there exists a (possibly degenerate) immersion $\varphi_0:N\rightarrow M$ to which $(\varphi_n)_\ninn$ subconverges. In the degenerate case, the image of $\varphi_0$ is a bundle of $(n-1)$-dimensional spheres over a complete geodesic. By compactness, it follows that $N=S^{n-1}\times S^1$, which contradicts the hypotheses. $\varphi_0$ is therefore not degenerate, and so:
$$
R_\theta(\varphi_0) = r.
$$
\noindent $(ii)$ follows, and this concludes the proof.\qed
\goodbreak
\newhead{Bibliography}
{\leftskip = 5ex \parindent = -5ex
\leavevmode\hbox to 4ex{\hfil \cite{Almgren}}\hskip 1ex{Almgren F. J. Jr., {\sl Plateau's problem. An invitation to varifold geometry.}, Student Mathematical Library, {\bf 13}, American Mathematical Society, Providence, RI, (2001)}
\medskip
\leavevmode\hbox to 4ex{\hfil \cite{AndBarbBegZegh}}\hskip 1ex{Andersson L., Barbot T., B\'eguin F., Zeghib A., Cosmological time versus CMC time in spacetimes of constant curvature, arXiv:math/0701452}
\medskip
\leavevmode\hbox to 4ex{\hfil \cite{Cabezas}}\hskip 1ex{Cabezas-Rivas E., Miquel V., Volume-preserving mean curvature flow in the hyperbolic cpace, {\sl Indiana Univ. Math. J.} {\bf 56} (2007), no.5, 2061--2086}
\medskip
\leavevmode\hbox to 4ex{\hfil \cite{CaffNirSprI}}\hskip 1ex{Caffarelli L., Nirenberg L., Spruck J., The Dirichlet problem for nonlinear second-Order elliptic equations. I. Monge Amp\`ere equation, {\sl Comm. Pure Appl. Math.} {\bf 37} (1984), no. 3, 369--402}
\medskip
\leavevmode\hbox to 4ex{\hfil \cite{CaffNirSprII}}\hskip 1ex{Caffarelli L., Kohn J. J., Nirenberg L., Spruck J., The Dirichlet problem for nonlinear second-order elliptic equations. II. Complex Monge Amp\`ere, and uniformly elliptic, equations, {\sl Comm. Pure Appl. Math.} {\bf 38} (1985), no. 2, 209--252}
\medskip
\leavevmode\hbox to 4ex{\hfil \cite{CaffNirSprIII}}\hskip 1ex{Caffarelli L., Nirenberg L., Spruck J., The Dirichlet problem for nonlinear second-order elliptic equations. III. Functions of the eigenvalues of the Hessian, {\sl Acta Math.} {\bf 155} (1985), no. 3-4, 261--301}
\medskip
\leavevmode\hbox to 4ex{\hfil \cite{CaffNirSprV}}\hskip 1ex{Caffarelli L., Nirenberg L., Spruck J., Nonlinear second-order elliptic equations. V. The Dirichlet problem for Weingarten hypersurfaces, {\sl Comm. Pure Appl. Math.} {\bf 41} (1988), no. 1, 47--70}
\medskip
\leavevmode\hbox to 4ex{\hfil \cite{ChenYau}}\hskip 1ex{Cheng S. Y., Yau S. T., On the regularity of the Monge-Amp\`ere equation\break $\opDet(\partial^2u\partial x_i\partial x_j)=F(x,u)$, \sl{Comm. Pure Appl. Math.} {\bf 30} (1977), no. 1, 41--68}
\medskip
\leavevmode\hbox to 4ex{\hfil \cite{GallKapMard}}\hskip 1ex{Gallo D., Kapovich M., Marden A., The monodromy groups of Schwarzian equations on closed Riemann surfaces, {\sl Ann. of Math.} {\bf 2} 151 (2000), no. 2, 625--704}
\medskip
\leavevmode\hbox to 4ex{\hfil \cite{Guan}}\hskip 1ex{Guan B., The Dirichlet problem for Monge-Amp\`ere equations in non-convex domains and spacelike hypersurfaces of constant Gauss curvature, {\sl Trans. Amer. Math. Soc.} {\bf 350} (1998), 4955--4971}
\medskip
\leavevmode\hbox to 4ex{\hfil \cite{GuanSpruckI}}\hskip 1ex{Guan B., Spruck J., The existence of hypersurfaces of constant Gauss curvature with prescribed boundary, {\sl J. Differential Geom.} {\bf 62} (2002), no. 2, 259--287}
\medskip
\leavevmode\hbox to 4ex{\hfil \cite{GuanSpruckII}}\hskip 1ex{Guan B., Spruck J., Szapiel M., Hypersurfaces of constant curvature in Hyperbolic space I, {\sl J. Geom. Anal}}
\medskip
\leavevmode\hbox to 4ex{\hfil \cite{GuanSpruckIII}}\hskip 1ex{Guan B., Spruck J., Hypersurfaces of constant curvature in Hyperbolic space II, arXiv:0810.1781}
\medskip
\leavevmode\hbox to 4ex{\hfil \cite{Gut}}\hskip 1ex{Guti\'errez C., {\sl The Monge-Amp\`ere equation}, Progress in Nonlinear Differential Equations and Their Applications, {\bf 44}, Birkh\"a user, Boston, (2001)}
\medskip
\leavevmode\hbox to 4ex{\hfil \cite{HarvLaws}}\hskip 1ex{Harvey R., Lawson H. B. Jr., Calibrated geometries, {\sl Acta Math.} {\bf 148} (1982), 47--157}
\medskip
\leavevmode\hbox to 4ex{\hfil \cite{Huisken}}\hskip 1ex{Huisken G., Contracting convex hypersurfaces in Riemannian manifolds by their mean curvature, {\sl Invent. Math.} {\bf 84}  (1986), no. 3, 463--480}
\medskip
\leavevmode\hbox to 4ex{\hfil \cite{LabI}}\hskip 1ex{Labourie F., Un lemme de Morse pour les surfaces convexes (French), {\sl Invent. Math.} {\bf 141} (2000), no. 2, 239--297}
\medskip
\leavevmode\hbox to 4ex{\hfil \cite{LabII}}\hskip 1ex{Labourie F., Probl\`eme de Minkowski et surfaces \`a courbure constante dans les vari\'et\'es hyperboliques (French), {\sl Bull. Soc. Math. France} {\bf 119} (1991), no. 3, 307--325}
\medskip
\leavevmode\hbox to 4ex{\hfil \cite{LoftI}}\hskip 1ex{Loftin J. C., Affine spheres and convex $\Bbb{RP}\sp n$-manifolds, {\sl Amer. J. Math.} {\bf 123} (2001), no. 2, 255--274}
\medskip
\leavevmode\hbox to 4ex{\hfil \cite{LoftII}}\hskip 1ex{Loftin J. C., Riemannian metrics on locally projectively flat manifolds, {\sl Amer. J. Math.} {\bf 124} (2002), no. 3, 595--609}
\medskip
\leavevmode\hbox to 4ex{\hfil \cite{MazzPac}}\hskip 1ex{Mazzeo R, Pacard P., Constant curvature foliations on asymptotically hyperbolic spaces, arXiv:0710.2298}
\medskip
\leavevmode\hbox to 4ex{\hfil \cite{RosSpruck}}\hskip 1ex{Rosenberg H., Spruck J., On the existence of convex hypersurfaces of constant Gauss curvature in hyperbolic space, {\sl J. Differential Geom.} {\bf 40} (1994), no. 2, 379--409}
\medskip
\leavevmode\hbox to 4ex{\hfil \cite{SchKra}}\hskip 1ex{Schlenker J. M., Krasnov K., The Weil-Petersson metric and the renormalized volume of hyperbolic $3$-manifolds, arXiv:0907.2590}
\medskip
\leavevmode\hbox to 4ex{\hfil \cite{SmiFCS}}\hskip 1ex{Smith G., Moduli of Flat Conformal Structures of Hyperbolic Type, arXiv:0804.0744}
\medskip
\leavevmode\hbox to 4ex{\hfil \cite{SmiSLC}}\hskip 1ex{Smith G., Special Lagrangian curvature, arXiv:math/0506230}
\medskip
\leavevmode\hbox to 4ex{\hfil \cite{SmiGKM}}\hskip 1ex{Smith G., Equivariant plateau problems, {\sl Geom. Dedicata} {\bf 140} (2009), 95--135}
\medskip
\leavevmode\hbox to 4ex{\hfil \cite{SmiHPP}}\hskip 1ex{Smith G., Hyperbolic Plateau problems, arXiv:math/0506231}
\medskip
\leavevmode\hbox to 4ex{\hfil \cite{SmiPKS}}\hskip 1ex{Smith G., Pointed $k$-surfaces, {\sl Bull. Soc. Math. France} {\bf 134} (2006), no. 4, 509--557}
\medskip
\leavevmode\hbox to 4ex{\hfil \cite{SmiPHD}}\hskip 1ex{Smith G., Probl\`emes elliptiques pour des sous-vari\'et\'es Riemanniennes, Th\`ese, Orsay, 2004}
\par}

\enddocument